\documentclass [11pt,oneside]{article}

\reversemarginpar \textwidth=14cm \textheight=19cm

\usepackage{amssymb} \usepackage{amsfonts} \usepackage{amsmath}
\usepackage{amsthm} \usepackage{epsfig} 

\newcommand{\mettifig}[1]{\epsfig{file=#1}}

\newtheorem{lemma}{Lemma}[section] 
\newtheorem{teo}[lemma]{Theorem}
\newtheorem{rem}[lemma]{Remark} 
\newtheorem{prop}[lemma]{Proposition}
\newtheorem{cor}[lemma]{Corollary} 
\newtheorem{conj}[lemma]{Conjecture}
 
\newtheorem{defn}[lemma]{Definition}

\newcommand{\matN}{\ensuremath {\mathbb{N}}}
\newcommand{\matR} {\ensuremath {\mathbb{R}}}
\newcommand{\matQ} {\ensuremath {\mathbb{Q}}}
\newcommand{\matZ} {\ensuremath {\mathbb{Z}}}

\newcommand{\matP} {\ensuremath {\mathbb{P}}}
\newcommand{\matH} {\ensuremath {\mathbb{H}}}
\newcommand{\matRP} {\ensuremath {\mathbb{RP}}}

\newcommand{\calL} {\ensuremath {\mathcal{L}}}
\newcommand{\calS} {\ensuremath {\mathcal{S}}}
\newcommand{\calM} {\ensuremath {\mathcal{M}}}

\newcommand{\calB} {\ensuremath {\mathcal{B}}}

\newcommand{\calT} {\ensuremath {\mathcal{T}}}

\newcommand{\calH}{\ensuremath {\mathcal{H}}}

\newcommand{\calG}{\ensuremath {\mathcal{G}}}

\newcommand{\nota} [1] {\caption{\footnotesize{#1}}}
\newcommand{\matr} [4] {\left(\begin{array}{@{}c@{\ }c@{}} #1 & #2 \\ #3 & #4 \\ \end{array} \right)}

\newfont{\Got}{eufm10 scaled 1200}
\newcommand{\permu}{{\hbox{\Got S}}}

\newcommand{\GL}{{\rm GL}}
\newcommand{\SL}{{\rm SL}}

\font\titsc=cmcsc10 scaled 1200

\newcommand{\basis}{\mathfrak{b}}
\newcommand{\Tsolid}{\hbox{\textit{\textbf{T}}}}            
\newcommand{\dimo}[1]{\vspace{2pt}\noindent\textit{Proof of \ref{#1}}.\ }
\newcommand{\finedimo}{{\hfill\hbox{$\square$}\vspace{2pt}}}

\newcommand{\timtil}{
\begin{picture}(12,8)
\put(2,0){$\times$}
\put(3,5){$\scriptstyle\sim$}
\end{picture}
}

\newcommand{\timforsetil}{
\begin{picture}(16,8)
\put(4,0){$\times$}
\put(2,5){${\scriptscriptstyle(}{\scriptstyle\sim}{\scriptscriptstyle)}$}
\end{picture}
}

\author{Bruno \titsc{Martelli}
\and Carlo \titsc{Petronio}}

\title{Complexity of geometric three-manifolds}

\begin{document}

\maketitle

\begin{abstract}
  We compute for all orientable irreducible geometric $3$-manifolds
  certain complexity functions that approximate from
  above Matve\-ev's natural complexity, known to be equal to the
  minimal number of tetrahedra in a triangulation. We can show that
  the upper bounds on Matveev's complexity implied by our computations
  are sharp for thousands of manifolds, and we conjecture they are for
  infinitely many, including all Seifert manifolds.  Our
  computations and estimates apply to all the Dehn fillings of
  $M6_1^3$ (the complement of the three-component chain-link,
  conjectured to be the smallest triply cusped hyperbolic manifold),
  whence to infinitely many among the smallest closed hyperbolic
  manifolds.  Our computations
  are based on the machinery of the
  decomposition into `bricks' of irreducible manifolds, developed in a
  previous paper. As an application of our results we completely
  describe the geometry of all 3-manifolds of complexity up to 9.
  
  \vspace{4pt}

\end{abstract}


\tableofcontents

\section*{Introduction}
The complexity $c(M)$ of a closed orientable
$3$-manifold $M$ was defined in~\cite{Mat} as the minimal number of
vertices of a simple spine of $M$. In the same paper it was shown that
$c$ is additive under connected sum and that, if $M$ is irreducible,
$c(M)$ equals the minimal number of tetrahedra in a triangulation of
$M$, unless $M$ is $S^3$, $\matRP^3$, or $L_{3,1}$. This means that
indeed $c(M)$ is a very natural measure of how complicated $M$ is.
Despite this fact, $c(M)$ is only known for a finite number of $M$'s.
In this paper we give upper bounds for the complexity of a wide class
of $M$'s, including many of the smallest known
 hyperbolic manifolds, all Seifert manifolds, and all torus
bundles over the circle. To be more precise we will introduce a
sequence $(c_n)_{n=0}^\infty$ of functions with values in
$\matN\cup\{\infty\}$ such that, given $M$, the sequence
$(c_n(M))_{n=0}^\infty$ is non-increasing and has constant value
$c(M)$ for $n\gg 0$. We will then explicitly compute $c_n$ for
$n\leqslant 9$ on all \emph{geometric} manifolds, \emph{i.e.}~on those 
carrying one of the eight locally homogeneous model geometries~\cite{Sco}.
The upper estimate on $c(M)$
given by the knowledge of $c_9(M)$ is known to be sharp for many
$M$'s, and we conjecture it to be
sharp for all torus bundles and all Seifert manifolds, except for certain
very special ones fibred over the sphere with three exceptional fibres
(a better estimate $c^*$ is defined for them).
It also allows to compute the complexity of 16 of the 
smallest known closed orientable hyperbolic manifolds.

This paper is based on the theory developed in~\cite{MaPe} and
reviewed below in Section~\ref{review:section}. We roughly sketch this
theory here to give the reader an idea of what $c_n$ is.
We define a \emph{marking} for a torus $T$ to be an embedding in $T$ of a $\theta$-graph
such that $T\setminus\theta$ is a disc, where the embedding is viewed up to isotopy.
We then call {\em manifold with marked boundary} an orientable, connected and 
compact 3-manifold with (possibly empty) boundary consisting of tori,
with a marking fixed on each boundary torus.
Using simple polyhedra we have defined in~\cite{MaPe} 
a complexity $c(M)$ for every manifold $M$ with marked boundary.
When $M$ is closed $c(M)$ equals Matveev's
complexity, and the natural properties of $c$ extend
from the closed to the marked-boundary case.
Moreover, given two manifolds with marked boundary $M$ and $M'$, we can
\emph{assemble} them by taking two tori $T\subset \partial M$ and $T'\subset \partial M'$
and gluing them via a map that matches their markings (for the sake
of simplicity we do not mention here the operation of
self-assembling, see Section~\ref{review:section}).
The complexity $c$ turns out to be subadditive under assemblings, so
we have called \emph{sharp} an assembling under which $c$ is
additive. Finally, we have called \emph{brick} a manifold with marked boundary 
that cannot be split as a
sharp assembling: it is then quite easy to show that every closed irreducible manifold 
can be sharply split
into bricks.  Using a computer, we have found that up to complexity
9 there are precisely 11
bricks with boundary, and 19 closed ones (that cannot be assembled
at all).

Adding to the 19 closed bricks the assemblings of the 11 bounded ones
we get all the $1901$ closed
irreducible manifolds of complexity at most $9$, and infinitely many
other manifolds of bigger complexity.
We can then define $c_n(M)$ as the smallest sum of complexities over the
expressions of $M$ as an assembling of bricks having complexity up to
$n$, with $c_n(M)=\infty$ if no such expression exists. 
The precise computation of $c_n$
for $n\leqslant 9$ on Seifert and hyperbolic manifolds is too
complicated to be reproduced here (see
Theorems~\ref{Seifert:bundle:teo} and~\ref{hyp:Dehn:fill:teo} below),
but there is one aspect that makes it remarkable and we can state now.
Of course $c_9(M)$ gives a better or equal upper estimate for $c(M)$
than $c_0(M),\ldots,c_8(M)$, but from the computation one sees that,
at least for most Seifert manifolds and all torus bundles, the
non-increasing sequence $(c_n(M))_{n=0}^9$ actually stabilizes very
early (namely at $n=1$ for lens spaces and torus bundles, and at $n=3$ for 
all other Seifert manifolds, except some atoroidal ones). 
This seems to suggest that
this stabilized value actually equals $c(M)$ for these $M$'s.

One of the key ingredients of our results is the complete classification
established in~\cite{MaPe:chain}
of all non-hyperbolic manifolds obtained by Dehn surgery
on the chain-link with three components (denoted by $6^3_1$
in~\cite{rolfsen}).
These manifolds are relevant to the theory because
two of our bricks are homeomorphic to the complement of $6^3_1$ 
(itself denoted by $M6^3_1$ in~\cite{Ca-Hi-We}),
with appropriate markings.

\vspace{.5cm}

This paper is structured as follows.  Since the 
precise definition of $c_n$ requires some details from~\cite{MaPe},
we have deferred it to Section~\ref{review:section}, 
collecting in Section~\ref{statement:section} all our statements, organized so
to be readable even without a complete understanding of what $c_n$ is.
Section~\ref{statement:section} also contains the definition
of two complexity functions (one for coprime pairs of integers
and one for invertible $2\times 2$ integer matrices) necessary to state our results.
Section~\ref{review:section} then contains the exact definitions of brick and $c_n$,
together with the preliminary constructions and results needed for these
definitions. In Section~\ref{slope:triod:section} we discuss the spaces of
slopes and $\theta$-graphs on a given torus, and we prove various properties
of the auxiliary complexity functions introduced in Section~\ref{statement:section}.
In Section~\ref{atoroidal:Seifert:torus:section} we start setting up the
combinatorial machine that we employ to compute $c_n$, beginning with the atoroidal case.
Our main result here is that this case essentially reduces to the understanding
of the Dehn fillings of the bricks, rather than more elaborate assemblings.
This fact is used in Section~\ref{atoroidal:Seifert:torus:section} itself to
compute $c_n$ on atoroidal Seifert manifolds, and later in 
Section~\ref{hyperbolic:complexity:section} for the Dehn fillings of
our hyperbolic bricks.
Section~\ref{Seifert:section} concludes the paper, with the computation of
$c_n$ for non-hyperbolic but geometric $3$-manifolds. This computation
is based on the extension to manifolds with marked boundary of the notion of Seifert
fibration and Euler number. In Section~\ref{tables:section} we provide a
description of all the closed, irreducible, and orientable
manifolds of complexity up to $9$.

\section{Complexity computations}\label{statement:section}
As briefly mentioned in the introduction, our theory of decomposition
of irreducible $3$-manifolds into bricks allows to define a non-increasing
sequence $\big(c_n\big)_{n=0}^\infty$ of approximated complexity functions such that,
for any given $M$, we have $c_n(M)=c(M)$ for $n\gg 0$.  Deferring the precise
definition of $c_n$ to Section~\ref{review:section}, we state here all the results
we can show about it.  We need to begin with some terminology and some preliminaries.

\paragraph{Terminology for $3$-manifolds}
All the manifolds we will consider are compact, connected, $3$-dimensional, and
orientable (but mostly unoriented).
We call \emph{lens
  space} any genus-1 manifold distinct from $\matRP^3$, and
\emph{Seifert manifold} any total space of a Seifert fibration over a
surface, possibly with boundary. We say that a Seifert manifold is
\emph{genuine} when it is not $S^3$, $\matRP^3$, a lens space, $T\times
I$, the solid torus $\Tsolid$, or a torus bundle over $S^1$.
For $A\in\SL_2(\matZ)$ we denote by $T_A$ the total space of the torus
bundle over $S^1$ with monodromy $A$.
Later we will employ the labels of the
Callahan-Hildebrand-Weeks census~\cite{Ca-Hi-We} for cusped
hyperbolic manifolds.  Namely, we denote by $Mi_j^k$ the $j$-th
manifold among those having an ideal triangulation with $i$ tetrahedra
and $k$ cusps, ordered by increasing volume.

We explain now the notation we use for the parameters of a
Seifert manifold. Let $F$ be surface, $t$ be an integer, and
$\{(p_i,q_i)\}_{i=1}^k$ be coprime pairs of integers with $|p_i|\geqslant 2$.
We denote by $\big(F,(p_1,q_1),\ldots,(p_k,q_k),t\big)$ the oriented
Seifert manifold obtained from either $F\times S^1$
or $F\timtil S^1$ (depending on orientability of $F$)
by drilling $k+1$ trivially fibred solid tori and Dehn filling
the resulting boundary components along the slopes
$p_1a_1+q_1b_1,\ldots,p_ka_k+q_kb_k,a_{k+1}+tb_{k+1}$, where
the $a_i$'s are contained in a section of the bundle, the
$b_i$'s are fibres, and each $(a_i,b_i)$ is a positive basis in homology.
Recall~\cite{FoMa} that for $i\leqslant k$ the $i$-th filling gives rise to
an exceptional fibre with ``orbital type'' $(|p_i|,r_i)$, where
$0<r_i<|p_i|$ and $q_ir_i \equiv p_i / |p_i|\ ({\rm mod}\ |p_i|)$, but we will never use
the orbital type, only the ``filling type'' $(p_i,q_i)$.  In addition, the
Euler number of the Seifert fibration is given by $t+\sum q_i/p_i$, but we
will never use the Euler number itself, only the ``twisting number'' $t$.

The Seifert parameters are called \emph{normalized} if $p_i>q_i>0$ for all $i$.
Any set of Seifert parameters can be promoted to a normalized one
by replacing each $(p_i,q_i)$ by $(|p_i|,q'_i)$, where
$0<q_i'<|p_i|$, $q'_i\equiv q_i p_i / |p_i|\ ({\rm mod}\ |p_i|)$, and
$t$ by $t+\sum(q_i/p_i-q'_i/|p_i|)$. The next result, deduced 
from~\cite{FoMa} and often needed below, shows that 
normalized parameters with an extra assumption correspond 
almost bijectively to genuine Seifert manifolds:

\begin{prop}\label{Seifert:generalities:prop} 
Let $M$ be a closed and oriented Seifert manifold.
\begin{itemize}
\item If $M$ is genuine and atoroidal then it can be expressed uniquely as 
$$\big(S^2,(p_1,q_1),(p_2,q_2),(p_3,q_3),t\big)$$
with $p_i>q_i>0$ for all $i$, and $t\in\matZ$. The only such expressions
that do not define atoroidal genuine Seifert manifolds
are the following torus bundles:
    \begin{eqnarray*}
      \big(S^2,(2,1),(3,1),(6,1),-1\big) &=& T_{\tiny{\matr 01{-1}1}} \\
      \big(S^2,(2,1),(4,1),(4,1),-1\big) &=& T_{\tiny{\matr 01{-1}0}} \\
      \big(S^2,(3,1),(3,1),(3,1),-1\big) &=& T_{\tiny{\matr 01{-1}{-1}}}.
    \end{eqnarray*}
\item If $M$ is genuine and toroidal then it can be expressed uniquely as 
$$\big(F, (p_1, q_1), \ldots, (p_k, q_k), t\big)$$
with $k-\chi(F)>0$, $p_i>q_i>0$ for all $i$, and $t\in\matZ$. 
The only such expressions
that do not define toroidal genuine Seifert manifolds are those with 
$F=S^2$ and $k=3$, and the following one:
$$\big(S^2,(2,1),(2,1),(2,1),(2,1),-2\big) = T_{\tiny{\matr {-1}00{-1}}}.$$
\item The expression of $-M$ is obtained from that of $M$ 
by replacing $t$ by $-t-k$ and each $q_i$ by $p_i-q_i$.
\end{itemize}
\end{prop}    

\begin{rem}\label{k:e:rem}
\emph{When a closed Seifert manifold 
$\big(F, (p_1, q_1), \ldots, (p_k, q_k), t\big)$
is viewed as an orientable but unoriented manifold, we can assume by the last
point of the previous proposition that $t\geqslant -k/2$, but this expression
can fail to be unique when the value of $t$ is precisely $-k/2$.}
\end{rem}

In addition to the standard terminology, we define now
$\calM^*$ to be the class of all genuine atoroidal Seifert manifolds 
of type $\big(S^2,(2,1),(n,1),(m,1),-1\big)$, that is
\begin{equation}\label{Mstar:list:eqn}
\calM^*=\big\{(S^2,(2,1),(n,1),(m,1),-1):\ 2\leqslant n\leqslant m,\ 
(n,m)\ne(3,6),(4,4)\big\}.
\end{equation}

\paragraph{Coprime pairs and integer matrices}
Our statements involve notions of complexity for coprime pairs of
integers and for invertible $2\times 2$ integer matrices. These
notions first appeared in~\cite{Mat:book}
and~\cite{Anisov}.

\begin{defn}\label{p:q:complexity:defn}
  {\em We recursively define an $\matN$-valued function
    $|\cdot\,,\cdot\,|$ on pairs of coprime integers by setting
    $|\pm 1,0|=|0,\pm 1|=|1,1|=0$ and
    $$|p,q|=1+ \left\{ \begin{array}{ll}
        |p,q-p| & {\rm if}\ q>p>0, \\
        |p-q,q| & {\rm if}\ p>q>0. \\
\end{array}\right.
$$
If $p>0>q$ we define $|p,q|=|p,-q|+1$.
If $p<0\neq q$ we define $|p,q|$ as $|-p,-q|$.}
\end{defn}

A direct method for computing $|p,q|$ will be described in 
Proposition~\ref{lines:and:|p,q|:prop}.
The following proposition will also be proved in
Section~\ref{slope:triod:section}.  It implies in particular that
$|p,q|$ is a well-defined function for a lens space $L_{p,q}$.

\begin{prop}\label{good:|p,q|:for:lens:prop}
\begin{enumerate}
\item $|p,q|=|q,p|$ for all coprime $p,q$;
\item If $0<q<p$ then $|p,p-q|=|p,q|$;
\item If $0<q,q'<p$ and $q\cdot q'\equiv \pm 1\ ({\rm mod}\ p)$ then
  $|p,q|=|p,q'|$.
\end{enumerate}
\end{prop}

Let us now consider the following matrices:
$$
S_1 = \matr{1}{-1}{0}{-1},\quad S_2 = \matr{-1}{0}{-1}{1},\quad S_3
= \matr{0}{1}{1}{0},\quad J = \matr{-1}001.
$$
The following proposition will be proved in
Section~\ref{slope:triod:section}.

\begin{prop} \label{matrix:decomposition:prop}
  Every $A\in\GL_2(\matZ)$ can be written in a unique way as
  $$A=\varepsilon\cdot S_{i_0}\cdot J\cdot S_{i_1} \cdot J\cdots
  S_{i_{n-1}}\cdot J\cdot S_{i_n}\cdot S_1^m$$
  with $n\geqslant 0$,
  $\varepsilon\in\{\pm1\}$, $i_0,i_n\in\{1,2,3\}$,
  $i_1,\ldots,i_{n-1}\in\{1,2\}$, and $m\in\{0,1\}$.
\end{prop}

\begin{rem}
  \emph{Since ${\rm det}(A)=(-1)^{m+1}$, uniqueness of $m$ in this expression is obvious.
Existence, and uniqueness of the other
    indices, will be harder to show.}
\end{rem}

\begin{defn}\label{matrix:complexity:defn}
  {\em For $A\in\GL_2(\matZ)$ we define $|A|\in\matN$ as the number
    $n$ of $J$'s appearing in the expression of $A$ just described.
    Moreover we define $||A||$ to be the minimum of $|BAB^{-1}|$ as
    $B$ varies in $\GL_2(\matZ)$.}
\end{defn}

\begin{rem} \label{easy:equalities:rem}
  {\em Since $J^{-1}=J$ and $S_i^{-1}=S_i$ for $i=1,2,3$, one readily
    sees from the definition that $|A^{-1}|=|A|$ when ${\rm
      det}(A)=-1$. The same equality can actually be checked also
    for ${\rm det}(A)=+1$, using also the relations
    $$S_1\cdot S_i=S_{(1\,3)(i)}\cdot S_3,\quad S_3\cdot J=-J\cdot
    S_3,\quad S_3\cdot S_i=S_{(1\,2)(i)}\cdot S_3,\quad S_3\cdot
    S_i=S_{(1\, 3)(i)}\cdot S_1.$$
    In particular we have
    $||A^{-1}||=||A||$ for all $A$.} 
\end{rem}

\begin{rem}
\emph{Given a certain $A\in\GL_2(\matZ)$, it is in general very difficult to
compute $|A|$ and $\|A\|$ using the definition directly.
A practical method to compute them will be given in Corollary~\ref{||:cor}.}
\end{rem}

\paragraph{Finiteness of $c_n$} 
We now state our main results, starting from finiteness of our
approximated complexity functions. Recall that a special list
$\calM^*$ of Seifert manifolds was introduced above in~(\ref{Mstar:list:eqn}).

\begin{teo}\label{finiteness:teo}
  Let $M$ be a closed, orientable, and irreducible $3$-manifold. If
  $M\not\in\calM^*$ we have
  $$c_0(M)\geqslant c_1(M)=c_2(M)\geqslant c_3(M)=\ldots
  =c_7(M)\geqslant c_8(M) \geqslant c_9(M).$$
  Moreover:
\begin{itemize}
\item $c_0(M)<+\infty$ if and only if $M\in\{S^3, \matRP^3, L_{3,1}\}$
  or $M=T_A$ with $||A||\leqslant 1$;
\item $c_1(M)<+\infty$ if and only if $M$ is either $S^3$, or
  $\matRP^3$, or a lens space, or a torus bundle over $S^1$;
\item $c_3(M)<+\infty$ if and only if $M$ splits along some tori into Seifert
  manifolds;
\item $c_8(M)<+\infty$ if and only if $M$ splits along some tori into Seifert 
  manifolds and/or copies of $M2_2^1$;
\item $c_9(M)<+\infty$ if and only if $M$ splits along some tori into Seifert 
  manifolds and/or copies of $M6_1^3$.
\end{itemize}
\end{teo}
Since a solid torus is a Seifert manifold, every Dehn filling of $M6_1^3$ has finite $c_9$. 
The manifold $M2_2^1$ is the figure-8 knot sister, and it is a filling of $M6_1^3$.

\paragraph{Computation of $c_n$: Seifert case}
We will now be much more precise on the values of the $c_n$'s on
torus bundles, on Seifert manifolds, and on hyperbolic manifolds.
Concerning non-geometric manifolds, we only mention here that the general computation of
$c_n$ would make an essential use
of the notion of complexity $|A|$ introduced above for a matrix
$A$ (see Lemma~\ref{c3:graph:lem}).

We start with torus bundles and Seifert manifolds, using the ``almost unique''
normalized parameters described in Proposition~\ref{Seifert:generalities:prop}
and Remark~\ref{k:e:rem}, 
and leaving to the reader to check that the formula for $c_n$ is
unchanged when a Seifert manifold admits two normalized expressions.

\begin{teo}\label{Seifert:bundle:teo}
  Let $M$ be a closed, orientable, and irreducible $3$-manifold.
\begin{itemize}
\item If $M\in\{S^3, \matRP^3, L_{3,1}\}$ then $c_n(M)=0$ for
  $n=0,\ldots, 9$;
\item If $M=T_A$ then:
\begin{itemize}
\item If $\|A\|\leqslant 1$ then $c_n(M)=6$ for all $n$; 
\item If $\|A\|\geqslant 2$ then $c_0(M)=+\infty$ and $c_n(M)=\|A\|+5$
  for $n=1,\ldots, 9$;
\end{itemize}
\item If $M$ is a lens space $L_{p,q}$ different from $L_{3,1}$ then
  $c_0(M)=+\infty$ and $c_n(M)=|p,q|-2$ for $n=1,\ldots, 9$;
\item If $M$ is a genuine Seifert manifold $\big(F, (p_1, q_1),
  \ldots, (p_k, q_k), t\big)$ with $k-\chi(F)>0$, $p_i>q_i>0$ for all $i$ and 
  $t\geqslant -k/2$, but $M$ is neither a member of $\calM^*$
  nor of the form $(S^2,(2,1),(3,1),(p,q),-1)$ with $p/q>5$, then
  $c_n(M)=+\infty$ for $n=0,1,2$ and
  $$c_n(M)=\max\{0,t-1+\chi(F)\}-6(\chi(F)-1) +
  \sum_{i=1}^k(|p_i,q_i|+2)$$
  for $n=3,\ldots,9$;
\item If $M=(S^2,(2,1),(3,1),(p,q),-1)$ with $p/q>5$ and $p/q\not\in\matZ$, then
  $c_n(M)=+\infty$ for $n=0,1,2$, $c_n(M)=|p,q|+3$ for $n=3,\ldots,7$,
  and $c_n(M)=|p,q|+2$ for $n=8,9$.
\end{itemize}
\end{teo}

\begin{rem}\label{Mstar:disclaimer:rem}
\emph{The only Seifert manifolds to which the previous result 
does not apply are the members of
$\calM^*$. Another estimate $c^*$, defined in Section~\ref{review:section},
seems to be more appropriate for them.}
\end{rem}

\paragraph{Fillings of $M6_1^3$}
By Theorem~\ref{finiteness:teo}, an atoroidal $M$ with
$c_9(M)<+\infty$ is either a small Seifert
manifold or a Dehn filling of $M6_1^3$.
In~\cite{MaPe:chain} we have completely classified all non-hyperbolic 
Dehn fillings of $M6_1^3$. 
We provide here only a simplified version of this result.
To state it, we recall again that $M6_1^3$ is the complement of a link in $S^3$,
and we fix on each boundary component the natural meridian-longitude homology basis.
Using this basis, we identify the set of slopes on each component
to $\matQ\cup\{\infty\}$.

\begin{teo}\label{short:chain:fill:teo}
Let $M$ be a (closed or partial) Dehn filling of $M6_1^3$.
Then $M$ is hyperbolic
except if one of the following occurs:
\begin{itemize} 
\item one of the filling coefficients belongs to
$\{\infty,-3,-2,-1,0\}$;
\item a pair of filling coefficients belongs to
$\{\{1,1\},\{-4,-1/2\},\{-3/2,-5/2\}\}$;
\item the triple of filling coefficients belongs to the following list:
$$
\begin{array}{c}
\{-5,-5,-1/2\},\ \{-4,-4,-2/3\},\ \{-4,-3/2,-3/2\},\\ \{-4,-1/3,1\},\ 
\{-8/3,-3/2,-3/2\},\ \{-5/2,-5/2,-4/3\},\\ \{-5/2,-5/3,-5/3\},\ \{-7/3,-7/3,-3/2\},\
\{1,2,2\},\\ \{1,2,3\},\ \{1,2,4\},\ \{1,2,5\},\ \{1,3,3\},\ \{2,2,2\}.
\end{array}
$$
\end{itemize}
\end{teo}

Since $M2_2^1$ is the $\{-1,4\}$ filling of $M6_1^3$, we also deduce:
\begin{cor}\label{B7:fill:cor}
  A $p/q$-Dehn filling of $M2^1_2$ is hyperbolic unless
  $$p/q\in\{\infty, -3, -2, -1, -1/2, -1/3, 0, 1\}.$$
\end{cor}

\paragraph{Computation of $c_n$: hyperbolic case}
Now we turn to the exact computation of $c_8$ and $c_9$ on the hyperbolic Dehn
fillings of $M6^3_1$. To do so we need to introduce a
certain numerical complexity $h$ for triples of rational numbers.
Namely, we define:
$$h(p/q,r/s,t/u)=g+|p+2q,q|+|r+2s,s|+|t+2u,u|$$
where 
$$g=
\left\{\begin{array}{l} 6 \qquad {\rm if\ }
    p/q,r/s,t/u {\rm\ are\ all\ different\ from\ } 1 {\rm\ and\ greater\ than\ }-2,\\
    4\qquad {\rm if\ } p/q=1,\ r/s\in\{-4,-5\}, {\rm\ and\ } t/u<-1 {\rm\ up\ to\
    permutation},\\
    2\qquad {\rm if\ } p/q=1,\ r/s\in\{-4,-5\}, {\rm\ and\ } t/u>-1 {\rm\ up\ to\
    permutation},\\
    5 \qquad {\rm in\ all\ other\ cases.}
\end{array}\right.$$

\begin{teo}\label{hyp:Dehn:fill:teo}
\begin{itemize}
\item If $M$ is a hyperbolic Dehn filling of $M2^1_2$
  then $c_8(M)$ is the minimum of
  $7+|p,-q|$ over all realizations of $M$ as a $p/q$-filling of
  $M2^1_2$;
\item If $M$ is a hyperbolic Dehn filling of $M6^3_1$ 
  then $c_9(M)$ is the minimum of
  $h(p/q,r/s,t/u)$ over all realizations of $M$ as a
  $(p/q,r/s,t/u)$-filling of $M6^3_1$.
\end{itemize}
\end{teo}

The somewhat implicit form of the previous statement is due to the
fact that very little is known about manifolds arising as different
Dehn fillings of the same cusped hyperbolic manifold
(see~\cite{cosmetic} for a careful discussion of the one-cusped case).
We have however shown in~\cite{MaPe:chain} that the following repetitions 
occur among the fillings of $M6^3_1$ and, after
extensive experiments carried out with SnapPea,
we put forward the conjecture that this is actually the complete list
of repetitions (up to permutation of the filling coefficients):

\begin{prop}\label{repetitions:prop}
Let $N_{p/q,r/s,t/u}$ be the closed Dehn filling of $N=M6^3_1$
with coefficients $\{p/q,r/s,t/u\}$. Then, for all $p/q$ and $r/s$ in $\matQ\cup\infty$,
\begin{eqnarray}
N_{-\frac 32,\frac pq,\frac rs} & = & N_{-4,\frac{q}{p+2q}-1,-\frac rs-3}\label{511right:eqn} \\
N_{-\frac 32,\frac pq,\frac rs} & = & N_{-\frac 32,-\frac q{p+2q}-2,-\frac
s{r+2s}-2}\label{511left:eqn} \\
N_{-\frac 52,\frac pq,\frac rs} & = & N_{-\frac 52,-\frac q{p+2q}-1,-\frac
s{r+s}-2}\label{512:eqn} \\
N_{-\frac 12,\frac pq,\frac rs} & = & N_{-\frac 12,-\frac pq-4,-\frac rs-4}\label{513:eqn} \\
N_{1,2,\frac pq} & = & N_{1,2,-\frac pq+2}\label{514:eqn}   \\
N_{1,-4,\frac pq} & = & N_{1,-4,\frac qp}. \label{515:eqn} 
\end{eqnarray}
\end{prop}

\paragraph{Hyperbolic manifolds up to complexity 10} 
Combining Theorems~\ref{short:chain:fill:teo} and~\ref{hyp:Dehn:fill:teo} 
it is now quite easy to give the complete list of all closed hyperbolic
manifolds of complexity up to $9$, and a partial list of those of complexity $10$.
The first list is obtained by enumerating the triples $\{p/q,r/s,t/u\}$
such that $h(p/q,r/s,t/u)\leqslant 9$, discarding those listed in
Theorem~\ref{short:chain:fill:teo}, and taking into account 
relations~(\ref{511right:eqn}) to~(\ref{515:eqn}) to eliminate duplicates.
This process leads to the following four triples:
$$\{-4,-3/2,1\},\qquad \{-4,1,2\},\qquad \{-5,-1/2,1\},\qquad \{-3/2,-3/2,1\},$$
which indeed give distinct manifolds, all 
having complexity $9$. More importantly, these manifolds
are precisely the four smallest known closed orientable hyperbolic manifolds, listed
in Table~\ref{hyp_nine:table}. (Both in this table and in the next one, Num is the position in the list
of known closed orientable hyperbolic manifolds, ordered by increasing volume:
the first 10 are in~\cite{Ho-We}, a longer list is distributed with
SnapPea).

\begin{table}
\begin{center} \begin{tabular}{c|c|c|c}
Slopes & Volume & Num & Homology \\ \hline
$\{-4,-3/2,1\}$ & 
$0.942707362777$ & $1$ & $\matZ_5 + \matZ_5$ \\ 
$\{-4,1,2\}$ & $0.981368828892$ & $2$ & $\matZ_5$ \\ 
$\{-5,-1/2,1\}$ & $1.01494160641$ & $3$ & $\matZ_3 + \matZ_6$ \\ 
$\{-3/2,-3/2,1\}$ & $1.26370923866$ & $4$ & $\matZ_5 + \matZ_5$ \\
\end{tabular} \end{center}
\nota{The closed orientable hyperbolic manifolds of complexity up to $9$.}\label{hyp_nine:table}
\end{table}

The partial list of closed hyperbolic manifolds having complexity $10$ is obtained
by the same method, except that we must consider the triples 
$\{p/q,r/s,t/u\}$ such that
$h(p/q,r/s,t/u)=10$, and discard those giving 
one of the manifolds
of Table~\ref{hyp_nine:table}. The result is described in Table~\ref{hyp_ten:table}.
The fact that all the triples in the tables contain slope $1$ implies that
all the manifolds listed are actually obtained by Dehn surgery on the
Whitehead link.

\begin{table}
\begin{center} \begin{tabular}{c|c|c|c}
Slopes & Volume & Num & Homology \\ \hline
$\{-5,1,2\}$ & $1.28448530047$ & $5$ & $\matZ_6$ \\ 
$\{-5,1/2,1\}$ & $1.41406104417$ & $7$ & $\matZ_6$ \\ 
$\{-4,1,3\}$ & $1.41406104417$ & $8$ & $\matZ_{10}$ \\ 
$\{-4,-4/3,1\}$ & $1.42361190029$ & $9$ & $\matZ_{35}$ \\ 
$\{-4,-5/2,1\}$ & $1.54356891147$ & $13$ & $\matZ_{35}$ \\ 
$\{-4,-5/3,1\}$ & $1.58316666062$ & $15$ & $\matZ_{40}$ \\ 
$\{-5,-3/2,1\}$ & $1.5886466393$ & $17$ & $\matZ_{30}$ \\ 
$\{-3/2,-1/2,1\}$ & $1.64960971581$ & $19$ & $\matZ_{15}$ \\ 
$\{-5,-1/3,1\}$ & $1.83193118835$ & $23$ & $\matZ_2+\matZ_{12}$ \\ 
$\{-5/3,-3/2,1\}$ & $1.88541472555$ & $31$ & $\matZ_{40}$ \\ 
$\{-5,-2/3,1\}$ & $1.91084379309$ & $34$ & $\matZ_{30}$ \\ 
$\{-3/2,-4/3,1\}$ & $1.95370831542$ & $39$ & $\matZ_{35}$ \\
\end{tabular} \end{center}
\nota{A partial list of closed orientable hyperbolic manifolds of complexity $10$.}\label{hyp_ten:table}
\end{table}

\paragraph{Geometric census} Using Theorems~\ref{Seifert:bundle:teo} 
and~\ref{hyp:Dehn:fill:teo} we have been able to refine the main result
of~\cite{MaPe} and determine the exact list and the geometry of all manifolds
up to complexity $9$. Our census is summarized in Table~\ref{geom_compl:table},
explained in greater detail in Section~\ref{tables:section}, and publicly 
available from~\cite{website}. 
The three manifolds of complexity $0$, \emph{i.e.}~$S^3$, $\matRP^3$, and $L_{3,1}$,
do not appear in Table~\ref{geom_compl:table}.
The number of manifolds wth complexity
$9$ is $1155$: the wrong number $1156$ in~\cite{MaPe} was the result
of a list containing the same graph manifold twice.

\newcommand{\zero}{$\cdot$}
\begin{table}  \begin{center} \begin{tabular}{r|cccccccccc}
complexity & 1 & 2 & 3 & 4 & 5 & 6 & 7 & 8 & 9 \\ \hline
lens & $2$ & $3$ & $6$ & $10$ & $20$ & $36$ & $72$ & $136$ & $272$ \\ 
other elliptic & \zero & $1$ & $1$ & $4$ & $11$ & $25$ & $45$ & $78$ & $142$ \\ 
flat  & \zero & \zero & \zero & \zero & \zero & $6$ & \zero & \zero & \zero \\ 
Nil  & \zero & \zero & \zero & \zero & \zero & $7$ & $10$ & $14$ & $15$ \\ 
$\mathbb{H}^2 \times S^1$  & \zero & \zero & \zero & \zero & \zero & \zero & \zero & $2$ & \zero \\ 
${\rm SL}_2$  & \zero & \zero & \zero & \zero & \zero & \zero & $39$ & $162$ & $514$ \\ 
Sol & \zero  & \zero & \zero & \zero & \zero & \zero & $5$ & $9$ & $23$ \\ 
hyperbolic  & \zero & \zero & \zero & \zero & \zero & \zero & \zero & \zero & $4$ \\ 
non-geometric  & \zero & \zero & \zero & \zero & \zero & \zero & $4$ & $35$ & $185$ \\ \hline
total  & $2$ & $4$ & $7$ & $14$ & $31$ & $74$ & $175$ & $436$ & $1155$ \\
\end{tabular} \end{center}
\nota{Number of manifolds per complexity and geometry}\label{geom_compl:table}
\end{table} 

\section{Decomposition of manifolds into bricks}\label{review:section}
In this section we provide the details needed below of the theory
developed in~\cite{MaPe} (see also~\cite{MaPe:nonori} for the
non-orientable version). In the whole paper our manifolds
will be 3-dimensional, compact, orientable, and connected by default,
and they will be viewed up to homeomorphism.

\begin{figure}
\begin{center}
\mettifig{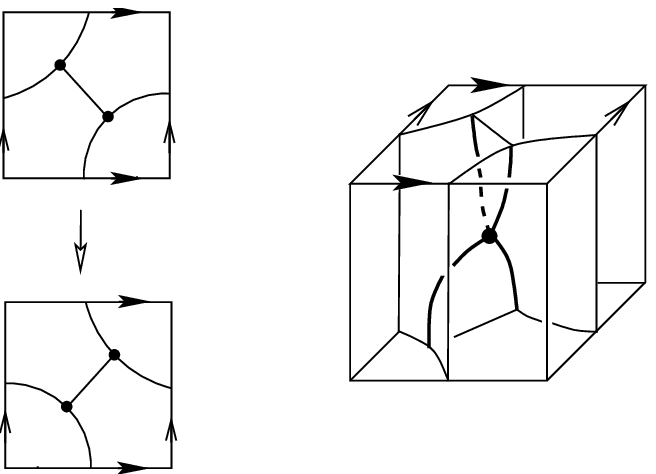} 
\nota{A flip on a $\theta$-graph in the torus, and a skeleton for $B_3$.} \label{flip:b3:fig}
\end{center}
\end{figure}

\paragraph{Manifolds with marked boundary}
As mentioned in the introduction, our theory deals with 
manifolds $M$ bounded by tori and endowed with a fixed \emph{marking}
of the boundary, \emph{i.e.}~a collection of $\theta$-graphs in $\partial M$
whose complement in $\partial M$ consists of discs. Note that
$\partial M$ can be empty.
As we will see in Section~\ref{slope:triod:section}, the same manifold can
be marked in infinitely many distinct ways. The easiest way to change
the marking of a boundary torus is to perform
a \emph{flip}, as shown in Fig.~\ref{flip:b3:fig}-left.

Now we describe two fundamental operations on the set of manifolds with marked boundary.
The first one, already defined in the introduction, is binary:
if $M$ and $M'$ are two such objects, take two tori $T\subset\partial
M, T'\subset\partial M'$ marked by $\theta\subset T, \theta'\subset
T'$ and a homeomorphism $\psi:T\to T'$ such that $\psi(\theta)=\theta'$.
Gluing $M$ and $M'$ along $\psi$ we get a new 3-manifold with marked
boundary. We call this operation (and its result) an \emph{assembling}.
Note that, although there are infinitely many non-isotopic maps between
two tori, only finitely many of them send one marking to the other, so
there is a finite number of non-equivalent assemblings of $M$ and $M'$.

Let us describe the second operation.
Let $M$ be a manifold with marked boundary, and $T,T'$ be 
distinct components of $\partial M$, with markings $\theta\subset T$ and 
$\theta'\subset T'$.
Let $\psi : T \to T'$ be a homeomorphism such
that $\psi(\theta)$ equals either $\theta'$ or a $\theta$-graph obtained
from $\theta'$ via a flip. The new 
manifold with marked boundary obtained from $M$ 
by gluing $T$ and $T'$ along $\psi$ is called a \emph{self-assembling}
of $M$. (There is a technical reason for loosening the natural requirement that 
$\psi(\theta)=\theta'$, mentioned later.)
There is only a finite number of non-equivalent self-assemblings of a given $M$.

\paragraph{Skeleta} 
A compact 2-dimensional
polyhedron $P$ is said to be \emph{simple} if the link of every point in $P$ is contained in
the 1-skeleton $K$ of the tetrahedron. A point, a compact graph, a
compact surface are thus simple. 
A point having the whole of $K$ as a link is 
called a \emph{vertex}.
The set $V(P)$ of the vertices of $P$ consists of isolated points, so
it is finite. Points, graphs and surfaces of course do not contain vertices.

A sub-polyhedron $P$ of a 3-manifold $M$ with marked boundary is
called a {\em skeleton} 
of $M$ if $P\cup\partial M$ is simple, $M\setminus(P \cup \partial M)$ is
an open ball and $P \cap \partial M$ is the union of the $\theta$-graphs 
marking $\partial M$. Note that if $M$ has one boundary component then
$P$ is a \emph{spine} of $M$, that is $M$ collapses onto $P$. Moreover, if $M$
is closed, then $P$ is a spine of $M$ minus a point.
It is easy to prove that every 3-manifold with marked boundary has a
skeleton.

\paragraph{Complexity}
The {\em complexity} $c(M)$ of a 3-manifold with marked boundary $M$ is
defined as the minimal number of vertices of a simple skeleton of 
$M$. When $M$ is closed, this definition coincides with Matveev's~\cite{Mat}.
Note that $c(M)$ depends on the topology of $M$ \emph{and on the marking}. In particular, if $T=\partial
M$ is one torus then every (isotopy class of a) $\theta$-graph on $T$
gives a distinct complexity for $M$.
Two properties extend from the closed case to the case with marked boundary:
complexity is still additive on connected sums, and it is finite-to-one on
orientable irreducible manifolds with marked boundary~\cite{MaPe}.

\paragraph{Examples}
Let $T$ be the torus. Consider $M=T\times I$, the boundary being
marked by $\theta_0\subset T\times\{0\}$ and $\theta_1\subset T\times\{1\}$. If $\theta_0$ and
$\theta_1$ are isotopic, the resulting 
manifold with marked boundary is called $B_0$. If $\theta_0$ 
and $\theta_1$ are related by a flip, we call the resulting manifold
with marked boundary $B_3$.
A skeleton for $B_0$ is $\theta_0 \times I$, 
while a skeleton for $B_3$ is shown in
Fig.~\ref{flip:b3:fig}-right. The skeleton of $B_0$ has no vertices, so $c(B_0)=0$. The
skeleton of $B_3$ has 1 vertex, and it can be
shown~\cite{MaPe} that there is no skeleton for $B_3$ without vertices, so $c(B_3)=1$.

\begin{figure}
\begin{center}
\mettifig{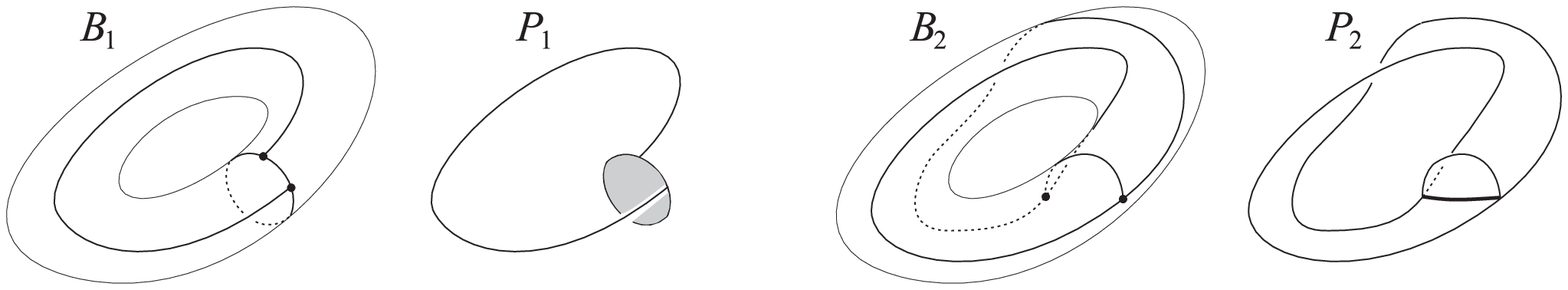, width = 13 cm} 
\nota{The solid tori with marked boundary $B_1$ and $B_2$, and skeleta for them.}
\label{B1_B2:fig}
\end{center}
\end{figure}

Two distinct marked solid tori are shown in Fig.~\ref{B1_B2:fig} and denoted
by $B_1$ and $B_2$. The same figure shows skeleta of $B_1$ and $B_2$
without vertices, so $c(B_1)=c(B_2)=0$.
Note that the skeleton of $B_1$ has a 1-dimensional portion.
All the other skeleta we will meet in this paper are fully 2-dimensional.

\begin{figure}
\begin{center}
\mettifig{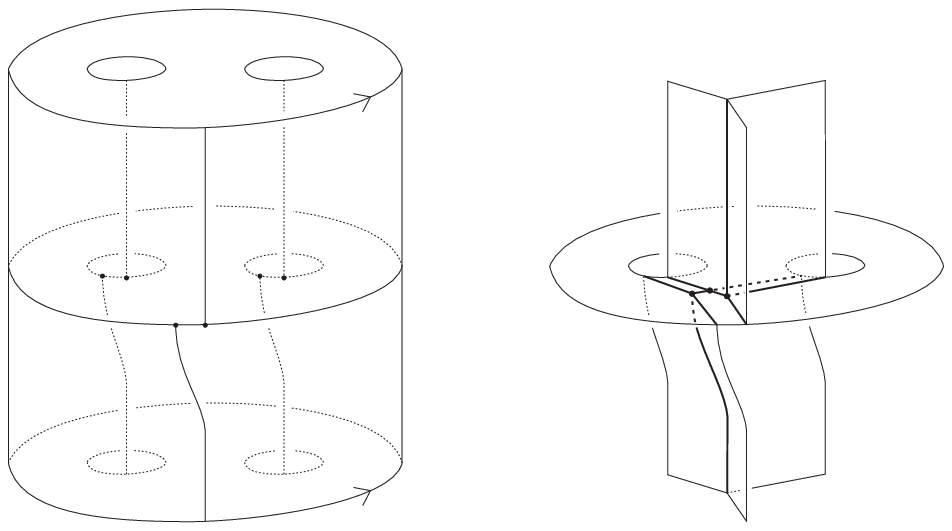}
\nota{The manifold with marked boundary $B_4$ (left) and a skeleton for it (right).} 
\label{B4:fig}
\end{center}
\end{figure}

The first irreducible orientable manifold with more than two marked
boundary components has complexity $3$ and is denoted by $B_4$. 
It is a marked $D_2\times S^1$, where $D_2$ is
the disc with two holes. To describe the markings of $B_4$, we note that on each
of the three boundary tori we have a natural figure-eight graph, obtained by
choosing basepoints in the factors of the product $(\partial
D_2)\times S^1$. There are two possible resolutions of the 4-valent
vertex of the graph, yielding $\theta$-graphs related by a flip. The
resolutions giving our $B_4$ are shown in Fig.~\ref{B4:fig}-left.
A skeleton with 3 vertices of $B_4$  is shown in Fig.~\ref{B4:fig}-right. 
It can be proved~\cite{MaPe} that $B_4$ has no skeleton with fewer vertices, so $c(B_4)=3$,
and that a different choice of marking would give higher complexity.

\paragraph{Subadditivity under assembling}
Let $M, M'$ be manifolds with marked boundary, and $P, P'$ be two
corresponding skeleta. An assembling of $M$
and $M'$ is given by a map $\psi$ that matches the $\theta$-graphs, so
$P\cup_\psi P'$ is a simple polyhedron inside $M\cup_\psi M'$. 
Moreover $P\cup_\psi P'$ is a
skeleton of the assembling $M\cup_\psi M'$, 
because gluing two balls along a boundary disc one gets a ball.
If $P,P'$ have $n, n'$ vertices, then $P\cup_\psi P'$ has $n+n'$,
therefore $c(M\cup_\psi M')\leqslant c(M)+c(M')$.

A similar construction works for self-assemblings and motivates its definition. 
Let $M'$ be
obtained by self-assembling $M$ along $\psi:T\to T'$.
By definition $\psi(\theta)$ either equals $\theta'$ or is obtained from $\theta'$
via a flip. In any case,
it is possible to isotope $\psi$ to $\psi'$ so that $\psi'(\theta)$ and $\theta'$
intersect each other transversely in 2 points.
We can now suppose $M'$ to be constructed using $\psi'$ rather than $\psi$,
and note that $P'=P\cup T$ is a skeleton of $M'$. 
Since $P'$ has at most $6$ vertices more than $P$,
we have $c(M')\leqslant c(M)+6$.

\paragraph{Bricks}
The theory culminates in a decomposition theorem. An assembling is called
\emph{sharp} if the complexity is additive and both manifolds with
marked boundary are irreducible and distinct from $B_0$, and a self-assembling is
\emph{sharp} if the complexity of the new manifold is the complexity
of the old one plus 6.
An irreducible orientable manifold with marked boundary is a {\em brick} if it is not the
result of a sharp assembling or self-assembling of other irreducible
manifolds with marked boundary. The proof of the following
result is clear: if an irreducible manifold with marked boundary is not a brick, then it can be
de-assembled. Then we repeat the analysis on each new piece. 
Since the sum of the complexities of all pieces does not increase
(and since the only possible pieces with complexity 0 are known to be
$B_1$ and $B_2$),
this iteration must stop after finite time.
\begin{teo} \label{split:teo}
Every irreducible orientable manifold with marked boundary can
be obtained from some 
bricks via a combination of sharp assemblings and sharp self-assemblings.
\end{teo}

The practical relevance of this result depends on the fact that the
set $\calB_{\leqslant n}$ of all bricks having complexity up to $n$ is
much smaller and easier to determine than the set of all manifolds up
to the same complexity.  We will soon describe $\calB_{\leqslant 9}$,
but we first note that closed bricks play a special (and marginal)
r\^ole in the theory, because they cannot be assembled at all.

\paragraph{Approximated complexity}
We can now formally introduce the functions $c_n$ to the study of
which the present paper is devoted.
Let $M$ be a manifold with marked boundary. If $M$ is prime we define
$c_n(M)$ as the minimum of $\sum_{i=1}^h c(B^{(i)})+6m$ over all
realizations of $M$ as an assembling of bricks
$B^{(1)},\ldots,B^{(h)}\in\calB_{\leqslant n}$ followed by $m$
self-assemblings. As usual, the minimum over an empty list is taken to
be $+\infty$. If $M$ is not prime we define $c_n(M)$ as the sum
of $c_n$ over the prime connected summands of $M$.
The next result is a direct consequence of Theorem~\ref{split:teo} and 
additivity of $c$ under connected sum~\cite{MaPe}:

\begin{prop}\label{approx:prop}
  For any $M$, the sequence
  $\{c_n(M)\}_{n=0}^\infty$ is monotonic non-increasing. Moreover
  $c_n(M)=c(M)$ for $n\geqslant c(M)$.
\end{prop}

\paragraph{Small bricks with boundary}
We have shown in~\cite{MaPe} that there are $11$ non-closed bricks
in $\calB_{\leqslant 9}$, namely the $B_0,\ldots,B_4$ introduced
above, and certain other $B_5,\ldots,B_{10}$.
The brick $B_5$ is the last Seifert one in our
list, being a marked $(D,(2,1),(3,1))$, where $D$ is the disc.
The other bricks $B_6,\ldots,B_{10}$ are (suitably marked)
cusped hyperbolic manifolds. Using as above the labels of~\cite{Ca-Hi-We},
we have that $B_6$,
$B_7$, and $B_8$ are $M2_2^1$, $M3_4^1$, and $M4^2_1$
respectively.  Both $B_9$ and $B_{10}$ are $M6^3_1$, with
different markings on the boundary.
We describe the markings of $B_5,\ldots,B_{10}$ in
Sections~\ref{atoroidal:Seifert:torus:section}
and~\ref{hyperbolic:complexity:section}.
Note that $M2_2^1$, $M3_4^1$, and $M4_1^2$ can
be obtained from the chain-link complement $M6^3_1$ via appropriate Dehn fillings.
Concerning complexity, we already know that $c(B_i)=0$ for $i=0,1,2$,
$c(B_3)=1$, and $c(B_4)=3$. 
In addition, $c(B_i)=8$ for
$i=5,6$ and $c(B_i)=9$ for $i=7,\ldots,10$.  

\paragraph{Minimal skeleta}
A skeleton $P$ of $M$ is \emph{minimal} if it has $c(M)$ vertices
and $P\cup\partial M$ cannot be collapsed onto a proper subpolyhedron.
The next technical result is used below.
It was also crucial for the computer program~\cite{MaPe} 
which classified $\calB_{\leqslant 9}$.

\begin{teo}\label{superstandard:teo}
  Let $B$ be a brick with $c(B)\geqslant 4$. Assume $B$ has $k$
  boundary components and $P$ is a minimal skeleton of $B$.  For 
  $i=1,\ldots, k$  let $e^i_1,e^i_2,e^i_3$ be the edges of the $\theta$-graph
  which marks the $i$-th component of $\partial M$. Then:
\begin{itemize}  
  \item for $i=1,\ldots, k$ and $j=1,2,3$ there is a $2$-cell $\sigma^i_j$ of $P$
  incident to $e^i_j$;
  \item for $i=1,\ldots, k$ there are vertices $v^i_1,v^i_2$ of $P$ such that
  $\sigma^i_j$ is incident to both of them for $j=1,2,3$;
  \item  the $2k$ vertices $\{v_\ell^i\}$ are all
  distinct from each other.
\end{itemize}  
\end{teo}

\paragraph{Small closed bricks}
Given integers $i>0$, $j>0$, $k\geqslant 0$ we employ the normalized parameters for
Seifert manifolds (Proposition~\ref{Seifert:generalities:prop}) to define
\begin{eqnarray*}
E_k & = & (S^2,(2,1), (3,1), (5+k,1), -1),\\
C_{i,j} & = & (S^2,(2,1), (1+i,1), (1+j,1), -1).
\end{eqnarray*}
In~\cite{MaPe} we have described a systematic way to construct skeleta
of $E_k$ and $C_{i,j}$, with $5+k$ and $i+j$ vertices respectively.
The following topological result is easy:

\begin{prop}\label{Mstar:other:list:prop}
The set $\calM^*$ of manifolds listed in equation~(\ref{Mstar:list:eqn})
coincides with the following one:
\begin{equation}\label{Mstar:other:list:eqn}
\begin{array}{l}
\big\{E_k:\ k\geqslant 0,\ k\ne 1\big\}\\
\qquad\bigcup\big\{C_{i,j}:\ 
1\leqslant i \leqslant j,\ (i,j)\ne (3,3),\ 
\textit{and}\ j\in\{2,3\}\ \textit{if}\ i=2\big\}.\end{array}
\end{equation}
Both lists~(\ref{Mstar:list:eqn}) and~(\ref{Mstar:other:list:eqn}) contain
no repetitions.
\end{prop}

For $M\in\calM^*$ define now $c^*(M)$ to be
$5+k$ or $i+j$ depending on the type of $M$.  In~\cite{MaPe} we have
shown that the closed bricks in $\calB_{\leqslant 9}$ are precisely the 19
elements of $\calM^*$ such that $c^*(M)\leqslant 9$. 

\paragraph{Results on $\calM^*$ and closed bricks} By 
definition if $M$ is a brick then
$c_n(M)$ is strictly bigger than $c(M)$ for all $n<c(M)$.  In
particular, the upper estimate $c_9(M)$ will never give $c(M)$ on a
closed brick $M$, except on the first 19 ones (those such
that $c(M)\leqslant 9$). A complete understanding of all closed bricks
would therefore be very desirable, and we will now state some
conjectures we have, starting from the result which motivates them.
(Notation is as above, points 1 and 2 were
proved in~\cite{MaPe}, and point 3 is proved below.)

\begin{teo}\label{Mstar:teo}
\begin{enumerate}
\item $c(M)\leqslant c^*(M)$ for all $M\in\calM^*$; 
\item If $c^*(M)\leqslant 9$ then $c(M)=c^*(M)$;
\item if $c^*(M)>9$ then $c_9(M)=c^*(M)+1$.  
\end{enumerate}
\end{teo}

This result suggests that all the elements of $\calM^*$ are bricks,
and that $c^*(M)$ is precisely the complexity of $M$ for
$M\in\calM^*$, and indeed we believe this is true.  With a slightly
smaller degree of confidence, we also believe that the 
only closed Seifert bricks are those
in $\calM^*$.  However, we do not risk a prediction on whether
geometrically toroidal or closed hyperbolic bricks exist in
complexity greater than 9.  Our experience with high complexity is too
limited to even make a conjecture. We are planning to investigate
manifolds and bricks of complexity 10 in the close future, and we hope
that the experimental data will provide some insight to these matters.

\section{Slopes and $\theta$-graphs on the torus}\label{slope:triod:section}
Our computation of $c_n$ will be based, among other things, on a
convenient and unified geometric encoding for simple closed curves and
$\theta$-graphs on the torus.  Our encoding uses the so-called Farey
tessellation and its dual graph, and it was partially inspired
by~\cite{Agol}. This approach leads to a very direct proof of the
fact, previously established in~\cite{Anisov}, that any two
$\theta$-graphs on the torus can be transformed into each other by a unique
minimal sequence of flips.  Quite surprisingly, it also leads to
geometric proofs of Propositions~\ref{good:|p,q|:for:lens:prop}
and~\ref{matrix:decomposition:prop} 
(algebraic proofs are also possible but quite awkward, so we have opted
for the geometric approach).

\paragraph{Slopes as points at infinity in hyperbolic plane}
Let us fix for the rest of the section a torus $T$, and let us denote
by $\calS(T)$ the set of \emph{slopes} on $T$, \emph{i.e.}~the isotopy
classes of non-contractible simple closed curves (without orientation)
on $T$. It is well-known that, given a basis $\basis=(a,b)$ of
$H_1(T)$, we have a bijection
$\Phi_\basis:\calS(T)\to\matQ\cup\{\infty\}$, where
$\Phi_\basis(\gamma)=p/q$ if the homology class of $\gamma$ is
$\pm(p\cdot a+q\cdot b)$. We note now that $\matQ\cup\{\infty\}$ is a
subset of the boundary $\matR\cup\{\infty\}$ of the upper half-plane
model of hyperbolic plane $\matH^2$, and for this reason we will
denote it by $\partial_\matQ\matH^2$. In the sequel it will be
convenient to freely switch from the half-plane to the Poincar\'e disc
model of $\matH^2$ without changing notation.  We do this by fixing
once and forever an isometry from the half-plane to the disc, namely
the orientation-preserving isometry that maps $\infty$ to $i$ and
$\{\infty,0,1\}$ to the vertices of a (Euclidean) equilateral
triangle.

It is now appropriate to recall that $\GL_2(\matZ)$ acts on $\matH^2$,
using the half-plane model, by fractional linear or anti-linear
transformations (depending on the sign of the determinant).  This
action extends to $\partial\matH^2$ and leaves $\partial_\matQ\matH^2$
invariant. On the other hand, given a basis $\basis$ of $H_1(T)$, we
also have an action of $\GL_2(\matZ)$ on $\calS(T)$, because every
$A\in\GL_2(\matZ)$ determines up to isotopy precisely one automorphism
$f$ of $T$ such that $f_*:H_1(T)\to H_1(T)$ has matrix $A$ with
respect to $\basis$.  It is now a routine matter to show that these
actions of $\GL_2(\matZ)$ are equivariant with respect to the
bijection $\Phi_\basis:\calS(T)\to\partial_\matQ\matH^2$.

\paragraph{The Farey tessellation}
We denote now by $\calL$ the set of lines in $\matH^2$ joining two
points
$\Phi_\basis(\gamma),\Phi_\basis(\gamma')\in\partial_\matQ\matH^2$
whenever $\gamma$ and $\gamma'$ intersect each other transversely in
only one point (up to isotopy). Using the $\matQ\cup\{\infty\}$ model
of $\partial_\matQ\matH^2$ we see that there is a line from $p/q$ to
$s/t$ precisely if $|p\cdot t-s\cdot q|=1$, so indeed $\calL$ is
independent of $\basis$. Moreover $\calL$ is invariant under the
action of $\GL_2(\matZ)$ on $\matH^2$.

\begin{lemma}\label{Farey:lem}
\begin{enumerate}
\item The lines in $\calL$ are pairwise disjoint (except at their
  ends);
\item For all $\ell\in\calL$ there exist precisely two points of
  $\partial_\matQ\matH^2$ which are joined by lines of $\calL$ to both
  the ends of $\ell$;
\item $\calL$ is the $1$-skeleton of an ideal triangulation of
  $\matH^2$.
\end{enumerate}
\end{lemma}

\begin{figure}\begin{center}
\mettifig{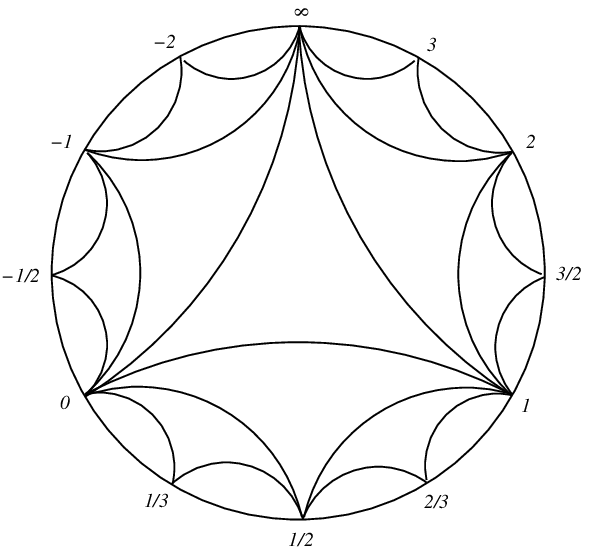} 
\nota{The Farey tessellation: an ideal
      triangulation of the Poincar\'e disc.}\label{tesselation:fig}
\end{center}\end{figure}

The ideal triangulation whose $1$-skeleton is $\calL$ is usually
called the \emph{Farey tessellation}~\cite{Farey1,Farey2}, and it is
shown in Fig.~\ref{tesselation:fig}.  For
the sake of simplicity
our pictures of $\calL$ will always be only combinatorially correct,
because metrically correct pictures would be harder to understand.

\dimo{Farey:lem} For point 1, let $\ell,\ell'\in\calL$. If $\ell$ and
$\ell'$ have a common end then they do not meet in $\matH^2$.  So we
can assume they have distinct ends. Using the fact that $\GL_2(\matZ)$
acts in a triply transitive way on $\partial_\matQ\matH^2$, we can now
assume in the half-space model that $\ell$ joins $0$ to $1$ and
$\ell'$ has one end at $\infty$.  The other end of $\ell'$ is then an
integer, and point 1 follows.  Point 2 is particularly easy to show
using the bijection between $\partial_\matQ\matH^2$ and $\calS(T)$:
the ends of $\ell$ represent two slopes $a$ and $b$ having geometric
intersection 1, and we must show that there are precisely two other
slopes having intersection 1 with both $a$ and $b$. If we choose
orientations on $a$ and $b$ we get a basis of $H_1(T)$, and the two
other slopes we are looking for are $\pm(a+b)$ and $\pm(a-b)$.  Point
3 now follows from points 1 and 2, and the density of
$\partial_\matQ\matH^2$ in $\partial\matH^2$.  \finedimo

\begin{figure}\begin{center}
\mettifig{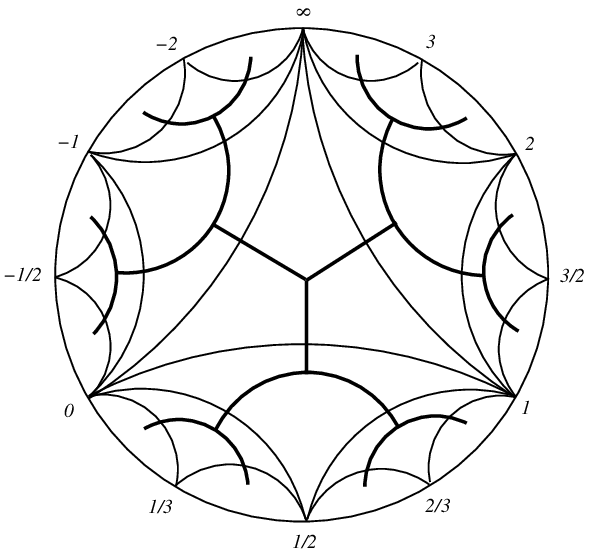}
    \nota{The
      trivalent tree dual to the Farey triangulation of $\matH^2$.}
\label{tree:fig}
\end{center}\end{figure}

\paragraph{$\theta$-graphs as vertices of the dual tree}
Let us denote now by $\Theta(T)$ the set of isotopy classes of $\theta$-graphs
in $T$ having a disc as a complement. It is not hard to see that the
elements of $\Theta(T)$ correspond bijectively to the unordered
triples of elements of $\calS(T)$ having pairwise geometric
intersection 1, the three slopes corresponding to a given
$\theta\in\Theta(T)$ being those contained in $\theta$ (up to
isotopy). So a triple $\{\gamma,\gamma',\gamma''\}$ defines a $\theta$-graph if
and only if
$\Phi_\basis(\gamma),\Phi_\basis(\gamma'),\Phi_\basis(\gamma'')$ are
the vertices of a triangle of the Farey triangulation of $\matH^2$.
Therefore, if we denote by $\calT$ the graph dual to the Farey
triangulation, we get a bijection $\Psi_\basis$ from $\Theta(T)$ to
the set $\calT^{(0)}$ of vertices of $\calT$.  Note that $\calT$ is a
tree embedded in $\matH^2$, see Fig.~\ref{tree:fig}.

The correspondence between $\Theta(T)$ and $\calT^{(0)}$ has a
geometric consequence. Note first that $\calT$ is trivalent,
\emph{i.e.}~there are precisely three $\theta$-graphs joined to a given
$\theta\in\Theta(T)$ by a single edge of $\calT$. Now one sees quite
easily that these three $\theta$-graphs are precisely those obtained from
$\theta$ by the elementary transformation called \emph{flip} and
described above in Fig.~\ref{flip:b3:fig}-left.  So we deduce:

\begin{prop}\label{unique:path:cor}
  Given distinct $\theta,\theta'\in\Theta(T)$ there exists a
  unique sequence $\{\theta_i\}_{i=0}^n\subset\Theta(T)$ of pairwise
  distinct $\theta$-graphs such that $\theta_0=\theta$, $\theta_n=\theta'$, and
  $\theta_i$ is obtained from $\theta_{i-1}$ by a flip, for
  $i=1,\ldots,n$.
\end{prop}

\paragraph{Distance between slopes and $\theta$-graphs}
Let us endow $\calT^{(0)}$ with the distance defined as usual for the
vertices of a graph. Namely, given $v,v'\in\calT^{(0)}$, their
distance $d(v,v')$ is the number of edges of the only simple path in
$\calT$ that joins $v$ to $v'$ (this path is unique because $\calT$ is
a tree). Using $\Psi_\basis$ we deduce a distance on $\Theta(T)$,
again denoted by $d$. This distance is independent of $\basis$,
because $d(\theta,\theta')$ is intrinsically interpreted as the number
of flips required to transform $\theta$ into $\theta'$.

We will need in the sequel a notion of ``distance'' between a slope
$\gamma\in\calS(T)$ and a $\theta$-graph $\theta\in\Theta(T)$. We define
$d(\gamma,\theta)$ to be $n-1$, where $n$ is the number of lines in
$\calL$ intersected by the half-line in $\matH^2$ that joins
$\Psi_\basis(\theta)$ to $\Phi_\basis(\gamma)$. Independence of
$\basis$ is easy and left to the reader. Of course our function $d$ is
not quite a distance, because $d(\gamma,\theta)=-1$ if $\gamma$ is
contained in $\theta$ up to isotopy.

\paragraph{Action of $\GL_2(\matZ)$ on $\theta$-graphs}
As above for slopes, we have now two different ways to define an
action of $\GL_2(\matZ)$ on $\theta$-graphs. On one hand, we can use the basis
$\basis$ to associate an automorphism of $T$ to a given matrix, and
then we can let this automorphism act on $\Theta(T)$. On the other
hand, the matrix acts by fractional linear or anti-linear
transformation on $\matH^2$ leaving $\calT^{(0)}$ invariant, so we
have an action defined through the bijection
$\Psi_\basis:\Theta(T)\to\calT^{(0)}$.  The reader is invited to check
that these actions of $\GL_2(\matZ)$ actually coincide.  

Our next results concern a relation between 
the distance function on the Farey tessellation and the complexity
$|p,q|$ for a coprime pair $(p,q)$ defined in
Section~\ref{statement:section}. 
To state the relation we introduce a notation repeatedly used below, namely
for $i\in\matZ$ we denote by $\theta^{(i)}\in\calT^{(0)}$ the 
centre of the triangle with vertices at $i$, $i+1$, and $\infty$. 
Note that $\theta^{(0)}$ is the
centre of the disc model of $\matH^2$.

\begin{prop}\label{lines:and:|p,q|:prop}
  Let $(p,q)$ be a pair of
  coprime integers.  Consider the point corresponding to
  $p/q$ in $\partial_\matQ\matH^2$.  Then
  $|p,q|$ equals the number of lines in $\calL$ met by the half-line
  that starts at $\theta^{(0)}$ and ends at $p/q$. Any
  other path from $\theta^{(0)}$ to $p/q$ meets at least $|p,q|$ lines of
  $\calL$.
\end{prop}

\begin{proof}
For $\alpha,\alpha'\in\calT^{(0)}\cup\partial_\matQ\matH^2$ 
let $\ell(\alpha,\alpha')$ be the geodesic segment, line, or half-line 
with ends at $\alpha$ and $\alpha'$, and let $n(\alpha,\alpha')$ be the 
number of intersections of $\ell(\alpha,\alpha')$ with the lines of $\calL$.

We must show that $|p,q|=n(p/q,\theta^{(0)})$.
We proceed by induction on $|p,q|$ and note that the equality is obvious when
$|p,q|=0$, namely for $p/q \in\{0,1,\infty\}$. Let us then consider
the case $|p,q|>0$. Suppose first that $p>q>0$. Let $f$ be the automorphism
of $\matH^2$ associated to $\tiny{\matr 1{-1}01}$. Then $f(\calL)=\calL$, 
$f(p/q)=(p-q)/q$, and $f(\theta^{(0)})=\theta^{(-1)}$, so
$n(p/q,\theta^{(0)})=n((p-q)/q,\theta^{(-1)})$. 
Now we note that $\ell(\theta^{(0)},\theta^{(-1)})$ only meets $\ell(0,\infty)$,
and $\ell(0,\infty)$ cuts the closed disc in two parts, leaving
$(p-q)/q$ and $\theta^{(0)}$ on one side and 
$\theta^{(-1)}$ on the other side. So 
$$n((p-q)/q,\theta^{(-1)})=n((p-q)/q,\theta^{(0)})+1.$$
Since $|p-q,q|=|p,q|-1<|p,q|$, the induction assumption now implies that
$n((p-q)/q,\theta^{(0)})=|p,q|-1$, whence the conclusion.

If $q>p>0$ we proceed similarly, using the automorphism
$\tiny{\matr 10{-1}1}$.  If $p>0>q$ we use the automorphism
$\tiny{\matr 1{-1}0{-1}}$, which fixes $\theta^{(0)}$ and maps $p/q$ to
$(p-q)/(-q)$.  So $n(p/q,\theta^{(0)})=n((p-q)/(-q),\theta^{(0)})$, 
and $n((p-q)/(-q),\theta^{(0)})=|p-q,-q|$ by the cases already discussed.
By definition of $|\,\cdot,\cdot\,|$ we now have $|p-q,-q|=|p,-q|+1=|p,q|$, whence
the conclusion also in this case. For $p<0$ we simply note that
$p/q=(-p)/(-q)$ and $|p,q|=|-p,-q|$ by definition.  The first assertion
of the statement is now proved, and the second one is easy, because
if a line $\ell$ of $\calL$ meets $\ell(p/q,\theta^{(0)})$ then it leaves 
$\theta^{(0)}$ and $p/q$ on opposite sides,
so any path joining
$\theta^{(0)}$ to $p/q$ must meet $\ell$.
\end{proof}

\begin{cor} \label{distance:cor}
If $p$ and $q$ are coprime then $d(p/q,\theta^{(i)})=|p-iq,q|-1$.
\end{cor}

\begin{proof}
The case $i=0$ follows from Proposition~\ref{lines:and:|p,q|:prop}.
The other cases follow using the automorphism $\tiny\matr 1{-i}01$
of $\matH^2$ that sends $\theta^{(i)}$ to $\theta^{(0)}$ and $p/q$ to
$(p-iq)/q$.
\end{proof}

With the geometric tools we have developed we can now establish 
the two technical results stated without proof in Section~\ref{statement:section}.

\dimo{good:|p,q|:for:lens:prop} If $p,q\geqslant 1$ then $|p,q|$ is
the number of operations required to reach $(p,q)$ from $(1,1)$ by
either adding the left coordinate to the right one, or the right
coordinate to the left one. Since in $(1,1)$ left and right are the
same, point 1 follows for $p,q\geqslant 1$.  The same assertion for
$p\leqslant 0$ or $q\leqslant 0$ is now easy and left to the reader.  

Using point 1 we can now very easily establish 
point 2: $|p,q|=|p-q,q|+1=|q,p-q|+1=|p,p-q|$.  

Let us prove point 3. By point 2 we can suppose $qq' \equiv 1
\ ({\rm mod}\ p)$. By Proposition~\ref{lines:and:|p,q|:prop} and point 1,
to prove that $|p,q|=|p,q'|$ it is now sufficient to show that 
\begin{equation}\label{pq:qp:eqn}
d(\theta^{(0)},p/q)=d(\theta^{(0)},q'/p). 
\end{equation}
We also note that the desired equality is obvious if $q=1$, so we proceed assuming
$p/q$ is not an integer.

For $r\in\partial_\matQ\matH^2$ we denote now by $\theta_r$ the
$\theta$-graph closest to $\theta^{(0)}$ among those that contain $r$. So
$d(r,\theta^{(0)})=d(\theta_r,\theta^{(0)})-1$.
Moreover, if $r$ is neither an integer nor $\infty$, one sees
that $\theta_r$ can be characterized as the only $\theta$-graph which contains
$r$ and two other slopes $r'$ and $r''$ such that $r'<r<r''$ with respect to the 
linear order of $\matQ$.

Now let $s>0$ be such that $qq'=sp+1$ and consider the following triples
of slopes:
$$\left(\frac {p-q'}{q-s}, \frac pq, \frac {q'}s\right), \qquad
\left(\frac sq,\frac {q'}p,\frac {q'-s}{p-q}\right).$$
Using the relation $qq'=sp+1$ one readily sees that in both triples the
slopes intersect pairwise in a single point, so the triples
define $\theta$-graphs.  The same relation and its
immediate consequence that $q>s$ also easily imply that in both triples the
slopes are arranged in increasing order.  This proves 
that the triples define $\theta_{p/q}$ and $\theta_{q'/p}$ respectively.
Consider now the automorphism $f$ corresponding to $\tiny{\matr p{-q'}q{-s}}$.
Then $f(\theta^{(0)})=\theta_{p/q}$ and $f(\theta_{q'/p})=\theta^{(0)}$. Since $f$ 
acts isometrically on $\calT^{(0)}$, relation~(\ref{pq:qp:eqn}) follows immediately, and
the proof is complete.
\finedimo

\dimo{matrix:decomposition:prop} 
Our argument is based on the action of $\GL_2(\matZ)$ on $\calT^{(0)}$,
that we analyze having in mind the correspondence between $\calT^{(0)}$
and the set of $\theta$-graphs on a torus with fixed homology basis.
We begin by providing an alternative description of $\calT^{(0)}$
that allows to understand the action of $\GL_2(\matZ)$ better.
We do this by noting that a $\theta$-graph can be viewed as a triple
of slopes intersecting pairwise once. So $\calT^{(0)}$
can be identified to the quotient of the set of $2\times
3$ integer matrices such that all three $2\times 2$ submatrices are
invertible over the integers, where two matrices
are identified if they are
obtained from each other by reordering and sign-switching of columns. With this
description of $\calT^{(0)}$, a matrix 
$A\in\GL_2(\matZ)$ acts on $[B]\in\calT^{(0)}$ giving just
$[A\cdot B]$. In addition, $[B]$ and $[B']$ are adjacent in $\calT$
(\emph{i.e.}~they are joined by an edge) if and only if they share two
columns (up to sign). We prove now a series of claims which eventually
will lead to the conclusion of the proof.

\emph{Claim 1: The matrices $S_1,S_2,S_3$ generate a subgroup of
  $\GL_2(\matZ)$ isomorphic to $\permu_3$, in which they play the
  r\^ole of the generating transpositions.}  The relations $S_i^2=I$
for $i=1,2,3$ and $S_{i_1}S_{i_2}S_{i_1}=S_{i_3}$ for
$\{i_1,i_2,i_3\}=\{1,2,3\}$ are proved by direct inspection and imply
the conclusion.

\emph{Claim 2: The stabilizer of $\theta^{(0)}$ in $\GL_2(\matZ)$ is
  generated by $-I,S_1,S_2,S_3$ and isomorphic to
  $\big(\matZ/_2\big)\times\permu_3$.}  This is established very easily
noting that
$$\theta^{(0)}=\left[ \begin{array}{ccc} 1 & 0 & 1 \\ 0 & 1 & 1 \\
  \end{array} \right].$$

\emph{Claim 3: If $A\in\GL_2(\matZ)$ and $A(\theta^{(0)})=\theta^{(0)}$ 
  then $A$ can
  be written in a unique way as $A=\varepsilon\cdot S_i\cdot S_1^m$
  with $\varepsilon\in\{\pm1\}$, $m\in\{0,1\}$, and $i\in\{1,2,3\}$.}  This is a
direct consequence of the previous claims.

\emph{Claim 4: Given $A\in\GL_2(\matZ)$, the three $\theta$-graphs in
  $\calT^{(0)}$ adjacent to $A(\theta^{(0)})$ are the following ones:}
$$(A \cdot S_1\cdot J)(\theta^{(0)}),\qquad (A\cdot S_2\cdot
J)(\theta^{(0)}),\qquad (A\cdot S_3\cdot J)(\theta^{(0)})=(A\cdot
J)(\theta^{(0)}).$$
The action of $A$ on $\calT^{(0)}$, being the
restriction of one on $\calT$, preserves adjacency of $\theta$-graphs. So it is
sufficient to check the claim with $A=I$. A straight-forward
computation shows that
$$(S_1\cdot J)(\theta^{(0)})= \left[ \begin{array}{ccc} 1 & 1 & 2 \\ 0 & 1
    & 1 \\ \end{array} \right],\qquad (S_2\cdot J)(\theta^{(0)})= \left[
  \begin{array}{ccc} 0 & 1 & 1 \\ 1 & 1 & 2 \\ \end{array} \right],$$
$$(S_3\cdot J)(\theta^{(0)})=J(\theta^{(0)})= \left[ \begin{array}{ccc} 1 & 0
    & 1 \\ 0 & 1 & -1 \\ \end{array} \right]$$
whence the conclusion.

\emph{Claim 5: Given $n\geqslant 1$, $i_0\in\{1,2,3\}$, and
  $i_1,\ldots,i_{n-1}\in\{1,2\}$, the sequence of $\theta$-graphs
\begin{eqnarray*}
\theta_0 & = & \theta^{(0)} \\
\theta_1 & = & (S_{i_0}\cdot J)(\theta^{(0)}) \\
\theta_2 & = & ((S_{i_0}\cdot J)\cdot (S_{i_1}\cdot J))(\theta^{(0)}) \\
& \cdots & \\
\theta_n & = & ((S_{i_0}\cdot J)\cdot (S_{i_1}\cdot J)\cdots(S_{i_{n-1}}\cdot J))(\theta^{(0)})
\end{eqnarray*}
is the sequence of vertices of a simple simplicial path in $\calT$.}
Using Claim 4 with $A=I$, this is obvious for $n=1$. Using Claim 4
again we only have to show for $n>1$ that $\theta_n\ne\theta_{n-2}$.
Define $A$ so that $\theta_{n-2}=A(\theta^{(0)})$, whence
$\theta_{n-1}=(A\cdot (S_{i_{n-2}}\cdot J))(\theta^{(0)})$ and
$\theta_n=(A\cdot(S_{i_{n-2}}\cdot J)\cdot(S_{i_{n-1}}\cdot
J))(\theta^{(0)})$.  We know that there will exist precisely one
$i\in\{1,2,3\}$ such that $(A\cdot (S_{i_{n-2}}\cdot J)\cdot (S_i\cdot
J))(\theta^{(0)})=\theta_{n-2}$. If we show that this is the case for
$i=3$, knowing that $i_{n-1}\in\{1,2\}$, we get the desired
conclusion. And, indeed:
\begin{eqnarray*}
(A\cdot (S_{i_{n-2}}\cdot J)\cdot (S_3\cdot J))(\theta^{(0)})& = &
(A\cdot(S_{i_{n-2}}\cdot J)\cdot J)(\theta^{(0)}) \\
& = & (A\cdot S_{i_{n-2}})(\theta^{(0)})=A(\theta^{(0)})=\theta_{n-2}.
\end{eqnarray*}

\emph{Claim 6: Given a length-$n$ simple simplicial path in $\calT$
  starting from $\theta^{(0)}$, the sequence of its vertices arises as in
  Claim 5 from a unique choice of $i_0\in\{1,2,3\}$ and
  $i_1,\ldots,i_{n-1}\in\{1,2\}$.}  This is a by-product of the
previous argument.

\emph{Conclusion.} Now let $A\in\GL_2(\matZ)$ and choose
$i_0\in\{1,2,3\}$ and $i_1,\ldots,i_{n-1}\in\{1,2\}$ to describe the
only simple path in $\calT$ from $\theta^{(0)}$ to $A(\theta^{(0)})$.  So
$A=S_{i_0}\cdot J\cdot S_{i_1} \cdot J\cdots J\cdot S_{i_{n-1}}\cdot
A'$ for some $A'$ in the stabilizer of $\theta^{(0)}$.  Claims 3 and 6 now
readily imply the conclusion.  \finedimo

Using the notation of the proof just completed, we note that the same
argument combined with the transitivity of the action of
$\GL_2(\matZ)$ on $\calT^{(0)}$ shows the following:

\begin{cor} \label{||:cor} 
  Given $A\in\GL_2(\matZ)$, if $|A|$ and $||A||$ are the functions
  defined in Section~\ref{statement:section} and $d$ is the graph
  distance on $\calT^{(0)}$, we have
  $$|A|=d(\theta^{(0)},A(\theta^{(0)})),\qquad
  ||A||=\min_{\theta\in\calT^{(0)}}d(\theta,A(\theta)).$$
\end{cor}

\section[Complexity of atoroidal manifolds: generalities 
and Seifert case]{Complexity of atoroidal manifolds:\\
  generalities and the Seifert
  case}\label{atoroidal:Seifert:torus:section} The computation of
$c_n$ for atoroidal manifolds with marked boundary is
considerably easier than for general ones. One reason is that we do
not have to consider self-assemblings.  More reasons will be clear
soon.

\paragraph{Effect of small bricks}
We begin by describing the effect
of the assembling with one of the smallest bricks $B_0,\ldots,B_3$.
The proof is immediate.

\begin{prop}\label{easy:assemblings:prop}
Let $M$ be a manifolds with marked boundary. Let $T$ be a component of
    $\partial M$, marked by some $\theta$.
    For $0\leqslant i\leqslant 3$, the effect of assembling $B_i$ to $M$ along $T$  is as
    follows: 
\begin{itemize}   
\item if $i=0$ then $M$ remains unaffected;
\item if $i=1$, a
    Dehn filling is performed on $M$ that kills one of the three
    slopes contained in $\theta$, \emph{i.e.}~a slope having distance
    $-1$ from $\theta$; all such fillings can be realized by an 
    assembling of $B_1$ with an appropriate map;
\item if $i=2$, a Dehn filling is performed on $M$
    that kills one of the three slopes having distance $0$ from
    $\theta$; all such fillings can be realized by an 
    assembling of $B_2$ with an appropriate map;
\item if $i=3$, the graph $\theta$ gets replaced by one of the three
    $\theta$-graphs having distance $1$ from $\theta$;
    all such replacements can be realized by an 
    assembling of $B_3$ with an appropriate map.
\end{itemize}
\end{prop}

We note here that the result just stated, together with the
knowledge of $\calB_{\leqslant 9}$, readily implies Theorem~\ref{finiteness:teo}.

\paragraph{Assemblings with a centre}
If $B$ is a brick, let us call \emph{assembling with centre} $B$ 
an assembling of $B$ (the centre) with some copies of $B_2$ and $B_3$.
We will prove in this paragraph that $c_n(M)$ is realized by an
assembling with centre whenever $M$ is atoroidal and $n\leqslant 9$.
This result will be crucial in the sequel. Its proof
relies on the following
experimental evidence, that we have checked for 
$n\leqslant 9$: \emph{the value of $c_n$ equals that of $c_1$ for marked
$\Tsolid$'s or $T\times I$'s}.  We formalize this evidence as
follows:

\begin{defn}
  \emph{An integer $n\geqslant 1$ is said to be \emph{well-behaved on
      tori} if $c_n(M)$ equals $c_1(M)$ whenever $M$ is a marked $\Tsolid$
    or $T\times I$.}
\end{defn}

\begin{conj} \label{torus:conj}
  All $n\geqslant 1$ are well-behaved on tori. That is, $c_1(M)=c(M)$ 
whenever $M$ is a marked $\Tsolid$ or $T\times I$.
\end{conj}

\begin{rem}
  \emph{If $n$ is well-behaved on tori, then all integers smaller than
    $n$ also are.}
\end{rem}

We will prove below that $n=9$ is well-behaved on tori.  To show the
main result of this paragraph we start with the following:

\begin{lemma}\label{no:B_1:lem}
  Let $B$ be a brick of complexity at most $n$ and let $M$ be the
	result of an assembling of $B$ with $h$
  copies of $B_1$.  Then
$$c_n(M) \leqslant \max\{c(B)-2h,0\}.$$
If $B\not\in\{B_0,\ldots,B_4\}$, then $c(B)-2h\geqslant 0$. If
$B\in\{B_1,\ldots,B_4\}$, then either $c_n(M)\leqslant c(B)-h$ or
$M$ is closed and $c(M)=0$.
\end{lemma}

\begin{proof}
Of course we may well assume that $c(B)=n$. The
assumptions imply that $c(M)\leqslant n$, so $c_n(M)=c(M)$. 
Therefore it is sufficient to 
prove the inequalities for $c(M)$ instead of $c_n(M)$.

Let us assume first
that $B$ is distinct from $B_0,\ldots, B_4$. 
Let $P$ be a minimal skeleton of $B$. Recall that a minimal
skeleton for $B_1$ consists of a meridinal disc and an arc on the boundary,
as shown in Fig.~\ref{B1_B2:fig}.
A skeleton $Q$ for $M$ is then obtained as follows:
for each component $T_i$ of $\partial M$ to which a $B_1$ is
assembled, a disc is glued to $P$ along a loop $\gamma$ contained in the
corresponding marking $\theta_i\subset T_i$. Therefore
$\theta_i\setminus\gamma$ is an edge of $\theta_i$, and 
Theorem~\ref{superstandard:teo}
implies that there is a 2-cell of $P$ incident
to it. In $Q$ this cell has a free edge, so it can
be collapsed.  Since we are assuming that 
$B\not\in\{B_0,\ldots,B_4\}$,
Theorem~\ref{superstandard:teo} again implies that collapsing this face one
eliminates two vertices from $Q$, and that the $2h$ vertices thus
eliminated are all distinct, whence $c(M)\leqslant c(B)-2h$, as required.

Let us prove now the first assertion for $B\in\{B_0,\ldots, B_4\}$.
For $B\in\{B_0,B_1,B_2\}$ the conclusion is obvious because
$c(M)=0$.  For $B=B_3$ and $B=B_4$ one checks by hand, looking at the
singular set of the skeleton of $B$, that
at least $\min\{2h,c(B)\}$ distinct vertices get eliminated when
collapsing $Q$.  

We are left to prove the last assertion for 
$B\in\{B_1,\ldots,B_4\}$.  Note that $B$ has at most $c(B)+1$ boundary components.
Assuming $c(M)>c(B)-h$ we deduce that $c(M)=0$ by the first assertion, so
$h>c(B)$. Then $h=c(B)+1$ and $M$ is closed.
\end{proof}

\begin{prop} \label{B_2:or:B_3:except:one:prop}
Let $n\geqslant 1$ be well-behaved on tori and $M$ be a manifold
with marked boundary. Assume $M$ is irreducible and atoroidal,
$c(M)>0$, and $c_n(M)<\infty$.
Then there is an assembling with centre realizing $c_n(M)$.
\end{prop}

\begin{proof}
Let us fix an assembling realizing $c_n(M)$ with a
minimal number $p\geqslant 0$ of bricks different from $B_2$ and
$B_3$. We must show that $p\leqslant 1$, so we assume by contradiction
that two bricks $B^{(1)}$ and $B^{(2)}$ are different from $B_2$ and $B_3$.
Let us concentrate on the boundary component $T$ of $B^{(1)}$ that leaves
$B^{(1)}$ and $B^{(2)}$ on opposite sides.
Since $M$ is
  irreducible and atoroidal, one of the following holds: (i) $T$
  bounds a solid torus in $M$; (ii) $T$ is parallel to a component of
  $\partial M$; (iii) $T$ is contained in a ball. Using the fact that
  each brick in the assembling is irreducible one easily sees that if (iii) holds
  then (i) also does. This implies that there is a sub-assembling
involving either $B^{(1)}$ or $B^{(2)}$ 
and giving a marked $\Tsolid$ or $T\times I$, which we denote by $N$.
Since $n$ is well-behaved on tori, we can replace 
the sub-assembling with one which still realizes $c_n(N)$ and involves
$B_0$, $B_1$, $B_2$, and $B_3$ only.
We can then dismiss the $B_0$'s. Moreover,
Proposition~\ref{no:B_1:lem} and relation $c_n(M)\geqslant c(M)>0$ imply that no $B_1$
is present, whence the conclusion.
\end{proof}

\paragraph{Complexity of Dehn fillings}

Let $N$ be a manifold with $\partial N=T_1\sqcup\ldots\sqcup T_k$.
Consider a $k$-tuple $\alpha_1,\ldots,\alpha_k$ where
$\alpha_i\in\calS(T_i)\cup\Theta(T_i)$, namely $\alpha_i$ is either a
slope or a $\theta$-graph in $T_i$.  We define a manifold with marked boundary
$N_{\alpha_1,\ldots,\alpha_k}$ obtained as follows from $N$: when
$\alpha_i\in\Theta(T_i)$, we mark $T_i$ with the $\theta$-graph $\alpha_i$;
when $\alpha_i\in\calS(T_i)$, we perform on $T_i$ the Dehn filling
that kills $\alpha_i$. 
We will use the notation $N_{\theta}$ with the same meaning, when 
$\theta=\{\theta_1,\ldots,\theta_k\}$ contains $\theta$-graphs only.

Recall now that we have defined in Section~\ref{slope:triod:section} a
distance $d(\alpha,\alpha')$ when $\alpha$ is either a slope or a
$\theta$-graph and $\alpha'$ is a $\theta$-graph on a given torus.  Given tori
$T_1,\ldots,T_k$ and two $k$-tuples
$\alpha=\{\alpha_1,\ldots,\alpha_k\}$ and
$\alpha'=\{\alpha_1',\ldots,\alpha_k'\}$, where $\alpha_i$ is either a
slope or a $\theta$-graph and $\alpha'$ is a $\theta$-graph on $T_i$, we now set
$d(\alpha,\alpha')=\sum_{i=1}^k d(\alpha_i,\alpha_i')$.

\begin{prop} \label{atoroidal:prop}
  Let $n\geqslant 1$ be well-behaved on tori and $M$ be 
  a manifold with marked boundary. Assume $M$ is irreducible and atoroidal, and
  $c(M)>0$. Then $c_n(M)$ equals the minimum of $c(N_\theta)+d(\alpha,\theta)$ over all
  realizations of $M$ as $N_{\alpha}$ and all
  $N_\theta\in\calB_{\leqslant n}\setminus\{B_1,B_3\}$, with $N_\theta\neq
  B_0$ if $M$ is not a marked $T\times I$.
\end{prop}

\begin{proof}
  We first prove inequality $\geqslant$, so we suppose
  $c_n(M)<\infty$.  Proposition~\ref{B_2:or:B_3:except:one:prop}
  implies that there exists an assembling that realizes $c_n(M)$
  involving one $B\in\calB_{\leqslant n}$ together with some $B_2$'s
  and $B_3$'s. We have $B\neq B_1$ by Lemma~\ref{no:B_1:lem}.  If
  $B=B_3$ and $M\neq T\times I$, at least one $B_2$ must occur in
  the assembling, and we redefine $B$ as one of these $B_2$'s.  If
  $B=B_3$ and $M=T\times I$, we can insert one $B_0$ and redefine
  $B$ to be this $B_0$. So in the assembling realizing $c_n(M)$ with
  centre $B$ we can actually assume that $B\in\calB_{\leqslant
  n}\setminus\{B_1,B_3\}$ and $B=B_0$ only if $M=T\times I$.

  Now let $h$ be the number of copies of $B_3$ appearing in our
  assembling, so that $c_n(M)=c(B)+h$. Let $B=N_\theta$, with
  $\theta=\{\theta_1,\ldots,\theta_k\}$. 
  Proposition~\ref{easy:assemblings:prop} shows that a successive assembling along
  a component $T_i$ of $\partial N$ of some $j$ copies of $B_3$
  corresponds to replacing $\theta_i$ by some $\theta'_i$ with
  $d(\theta_i,\theta'_i)\leqslant j$. The same proposition implies that if we
  first assemble $j$ copies of $B_3$ and then a copy of $B_2$ then $M$ gets
  filled along a slope $\alpha_i$ with $d(\alpha_i,\theta_i)\leqslant
  j$.  This shows that $M=N_\alpha$ with $d(\alpha,\theta)\leqslant h$,
  whence $c_n(M)\geqslant c(B)+d(\alpha,\theta)$, and inequality
  $\geqslant$ follows.

  Turning to inequality $\leqslant$, suppose the right-hand side
  equals $c(N_\theta)+d(\alpha,\theta)$ for some brick $N_\theta\in\calB_{\leqslant
    n}\setminus\{B_1,B_3\}$ and some realization $M=N_\alpha$.
  Suppose $\partial N=T_1\sqcup \ldots\sqcup T_k$,
  $\theta=\{\theta_1,\ldots,\theta_k\}$,
  $\alpha=\{\alpha_1,\ldots,\alpha_k\}$.  We can now associate to the
  realization $N_\alpha$ of $M$ an assembling yielding $M$.
  When $\alpha_i\in\Theta(T_i)$  we assemble $d(\alpha_i,\theta_i)$
  copies of $B_3$, following the unique simple path in $\Theta(T_i)$
  from $\theta_i$ to $\alpha_i$. When $\alpha_i\in\calS(T_i)$ and
  $d(\alpha_i,\theta_i)\geqslant0$ we assemble $d(\alpha_i,\theta_i)$
  copies of $B_3$ and one of $B_2$.  When $\alpha_i\in\calS(T_i)$ and
  $d(\alpha_i,\theta_i)=-1$ we assemble one copy of $B_1$.
	Now we deduce from the last assertion of Lemma~\ref{no:B_1:lem} 
	that the partial assembling of $N$ with the $B_1$'s only has complexity $c_n$
	at most $c(N)-h$, where $h$ is the number of $B_1$'s. Adding also the $B_2$'s and
	$B_3$'s we then get $c_n(M)\leqslant c(N_\theta)+ d(\alpha,\theta)$, whence the conclusion.
\end{proof}

Let $N_\theta$ be a brick of complexity $n\geqslant 1$.  From
Proposition~\ref{atoroidal:prop} we know that $c_n(N_{\alpha})\leqslant
n+d(\alpha,\theta)$ for any $k$-tuple $\alpha$. We then say that $\alpha$
is \emph{degenerate} if $c_n(N_{\alpha})< n+d(\alpha,\theta)$ or
$c(N_{\alpha})=0$.  If $S$ is a set of slopes lying on pairwise
distinct components of $\partial N$, we add to $S$ the $(\#\theta)-(\#S)$ $\theta$-graphs
of $\theta$ on the tori that $S$ does not intersect, and we call $S$ degenerate
if the resulting $k$-tuple is.

\begin{prop} \label{degenerate:prop}
Let $N_\theta$ be a brick of complexity $n$, with $\partial N=T_1\sqcup\ldots\sqcup T_k$.
\begin{enumerate}
\item\label{degenerate:contained:point} If $n\geqslant 1$ then every slope
  $s_i\in\calS(T_i)$ with $d(s_i,\theta_i)=-1$ is degenerate; 
\item\label{degenerate:superset:point} If $n\geqslant 2$
  then any $k$-tuple containing a degenerate set of slopes is
  degenerate.
\end{enumerate}
\end{prop}
\begin{proof}
  The first assertion follows from Lemma~\ref{no:B_1:lem} and
  Proposition~\ref{easy:assemblings:prop}. To prove the second one, let
  $\alpha_1,\ldots,\alpha_k$ be a $k$-tuple where
  $\{\alpha_1,\ldots,\alpha_j\}$ is a degenerate set of slopes for
  some $1\leqslant j\leqslant k$.  Namely, one of the following holds:
\begin{eqnarray*}
c_n(N_{\alpha_1,\ldots,\alpha_j,\theta_{j+1},\ldots,\theta_k})&=& 0,\qquad{\rm or} \\
c_n(N_{\alpha_1,\ldots,\alpha_j,\theta_{j+1},\ldots,\theta_k})&<&
n+\sum_{i=1}^j d(\alpha_i,\theta_i).
\end{eqnarray*}
If the second relation does not hold, so the first does, we deduce
that $n\leqslant k$, whence $n\leqslant 3$ by
Theorem~\ref{superstandard:teo}. But $n\geqslant 2$, so $N_\theta=B_4$
and $k=j=3$, and we have nothing to prove. So we can proceed assuming
the second relation holds. Now Lemma~\ref{no:B_1:lem} and
Proposition~\ref{easy:assemblings:prop} again imply that either
$c_n(N_{\alpha_1,\ldots,\alpha_k})=0$ or
$$c_n(N_{\alpha_1,\ldots,\alpha_k})\leqslant
c_n(N_{\alpha_1,\ldots,\alpha_j,\theta_{j+1},\ldots,\theta_k})+
\sum_{i=j+1}^k d(\alpha_i,\theta_i),$$
so
$c_n(N_{\alpha_1,\ldots,\alpha_k})<n+d(\alpha,\theta)$.
\end{proof}

The following result will be used later in this section and 
in Section~\ref{hyperbolic:complexity:section} 
to show that $n$ is well-behaved on tori for $n\leqslant 9$.

\begin{prop}\label{well:sufficient:prop}
Let $n>1$ and assume the following holds for every brick $B$ with 
$1<c(B)\leqslant n$:
\begin{itemize}
\item $B$ is not a marked $\Tsolid$ or $T\times I$;
\item If $\alpha$ is a non-empty set of \emph{slopes} on $\partial B$
and $B_\alpha$ is a marked $\Tsolid$ or $T\times I$
then $\alpha$ is degenerate.
\end{itemize}
Then $n$ is well-behaved on tori.
\end{prop}

\begin{proof}
Assume $n$ is not well-behaved on tori.
An argument similar to the proof of Proposition~\ref{B_2:or:B_3:except:one:prop} 
shows that there is a marked $\Tsolid$ or $T\times I$, which we denote by $M$, and
an assembling with centre $B$ realizing $c_n(M)<c_1(M)$. Of course $1<c(B)\leqslant n$.
Then, as in the proof of
Proposition~\ref{atoroidal:prop}, one shows that there is 
a $k$-tuple $\beta$ such that $B_\beta = M$
and $\beta$ is non-degenerate.
Let $\alpha$ be obtained from $\beta$ by dismissing the $\theta$-graphs and
keeping the slopes.  Then $B_\alpha$ is also 
a marked $\Tsolid$ or $T\times I$, so $\alpha$ is non-empty and
degenerate by assumption. But $\beta$ contains $\alpha$, which contradicts
Proposition~\ref{degenerate:prop}-(\ref{degenerate:superset:point}).
\end{proof}

\begin{rem}\label{few:cases:rem}
\emph{We emphasize that the set $\alpha$ in the previous statement consists of
slopes only, not of $\theta$-graphs too. We also recall that 
that none of the bricks of complexity between 2 and 9 is a marked $\Tsolid$ or
$T\times I$. Moreover, of the non-closed bricks
$B_4,\ldots,B_{10}$, we know that $B_5,B_6,B_7$ have only one boundary component.
Therefore, to conclude that $n=9$ is well-behaved on tori, 
we will only need to show that all the fillings of $B_4,B_8,B_9,B_{10}$
giving $T\times I$ or $\Tsolid$ are degenerate. Moreover $B_4$ alone
already suffices for $n$ up to $8$.}
\end{rem}

\paragraph{Computation of $c_1$}
We start by computing $c_n(M)$ for $n=1$ and any irreducible and
atoroidal marked $M$.  For $n\geqslant 3$ we will confine ourselves to the
closed case.
The bricks of complexity $0$ and $1$ are marked $T\times I$'s and $\Tsolid$'s,
so $c_1(M)<\infty$
precisely when $M$ is $T\times I$, $\Tsolid$, a lens space, $S^3$,
or $\matRP^3$. And of course $c_1(S^3)=c_1(\matRP^3)=0$. 

We now establish two results that will be used to
compute $c_1$ on torus bundles and lens spaces respectively.
Here and in the sequel we 
will denote by $f:\Theta(T\times\{0\})\to \Theta(T\times\{1\})$ 
the obvious map given by the product structure of $T\times I$.

\begin{prop} \label{T:times:I:prop}
If $\theta_i$ marks $T\times\{i\}$ then
$c_1\big((T\times I)_{\theta_0,\theta_1}\big)= d(f(\theta_0),\theta_1)$.
\end{prop}  
\begin{proof}
Set $M=(T\times I)_{\theta_0,\theta_1}$. Proposition~\ref{atoroidal:prop} implies that 
$c_1(M)$ equals the minimum of
$d(\theta'_0,\theta)+d(\theta'_1,f(\theta))$ over all realizations of $M$
as $(T\times I)_{\theta'_0,\theta'_1}$, where $\theta'_0$ and $\theta'_1$ vary and
$\theta$ is fixed. Now
$d(\theta'_0,\theta)+d(\theta'_1,f(\theta))\geqslant
d(f(\theta'_0),\theta'_1)$.  Moreover $d(f(\theta'_0),\theta'_1)=d(f(\theta_0),\theta_1)$, 
and it is possible to realize $M$ with
$\theta'_0=\theta$, whence the equality.
\end{proof}

\begin{figure}\begin{center}
\mettifig{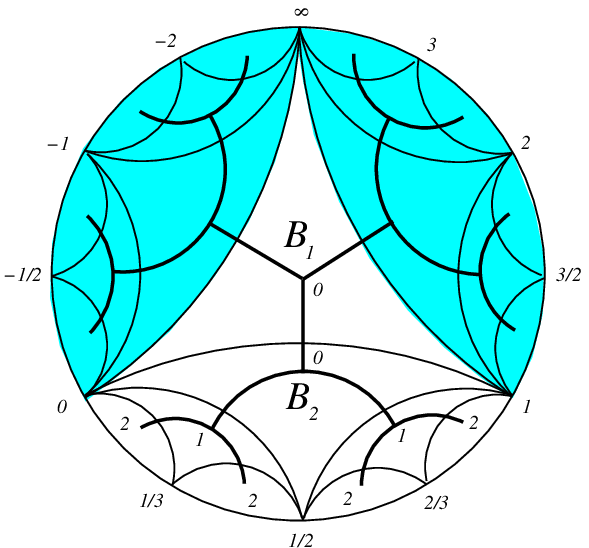}
\nota{Marked solid tori and their complexity.}
\label{solid:fig}
\end{center}\end{figure}

\begin{prop} \label{Tsolid:prop}
There is a natural correspondence between the set of \emph{oriented} marked solid tori and the
vertices in the non-shadowed region of Fig.~\ref{solid:fig}. 
If $d$ denotes the distance induced by this correspondence then
$c_1(\Tsolid_\theta) = d(\Tsolid_\theta, B_2)$ for every marked $\Tsolid_\theta$ distinct from $B_1$.
\end{prop}

\begin{proof}
Choosing a meridian-longitude basis $\basis$
for $H_1(\partial\Tsolid)$, we can identify the space of marked solid tori
up to orientation-preserving homeomorphism to the quotient of
$\Theta(\partial\Tsolid)$ under the action of the subgroup of
$\SL_2(\matZ)$ consisting of upper triangular matrices.  Via
the map $\Psi_\basis$ of Section~\ref{slope:triod:section} we can also
identify $\Theta(\partial\Tsolid)$ to the set of vertices of the graph
$\calT$ dual to the Farey tessellation of $\matH^2$.  Upper triangular
matrices in $\SL_2(\matZ)$ act on $\matH^2$
as parabolic isometries fixing point $\infty\in\partial_\matQ\matH^2$, and a
fundamental domain is delimited by two consecutive lines of $\calL$
ending at $\infty$.
Finally,
Proposition~\ref{atoroidal:prop} implies that if $c(\Tsolid)>0$
then $c_1(\Tsolid)$ equals $d(\Tsolid, B_2)$.
\end{proof}

Figure~\ref{solid:fig} shows the positions of
$B_1$ and $B_2$ and some values of $c_1$.
Orientation-reversal corresponds in the figure
to the reflection in the vertical axis,
which explains why the values of $c_1$ are symmetric with respect to it.

Turning to lens spaces, we have:

\begin{prop} \label{c1:prop} Let $L_{p,q}$ be a lens space. Then
$c_1(L_{p,q})=|p,q|-2$.
\end{prop}
\begin{proof}
We have $L_{p,q} = \Tsolid_\alpha$ with a slope $\alpha$ having coordinates $(q,p)$ with
respect to a meridian-longituse basis.  Figure~\ref{solid:fig} and
Proposition~\ref{lines:and:|p,q|:prop} now imply that the distance
between $\alpha$ and the $\theta$-graph of $B_2$ is $|q,p|-2$. Our result follows
from Proposition~\ref{atoroidal:prop},
since $|q,p|=|p,q|$ by Proposition~\ref{good:|p,q|:for:lens:prop}-(2).
\end{proof}

Note that by convention we view $L_{p,q}$ as a honest lens space
precisely when $|p,q|\geqslant 2$.

\paragraph{A preview on torus bundles}
Even if this section is devoted to atoroidal manifolds,
to proceed we need to state the (easy) computation of $c_1$ on torus
bundles. Proposition~\ref{cn:torus:bundle:prop} will show that 
$c_n$ equals $c_1$ on these manifolds for
$n=2,\ldots,9$.

\begin{prop} \label{c1:torus:bundle:prop} 
If $A\in \SL_2(\matZ)$ and $T_A$ is the torus bundle
over the circle with
monodromy $A$ then:
  \begin{itemize}
  \item $c_0(T_A)$ equals $6$ if $||A||\leqslant 1$ and
    $\infty$ otherwise;
  \item $c_1(T_A) = \max\{||A||+5,6\}$.
\end{itemize}
\end{prop}

\begin{proof}
  The only way to obtain a torus bundle $T_A$ using bricks having complexity
  0 is to do a self-assembling of
  $B_0=(T\times I)_{\theta, f(\theta)}$ along a map $\psi$
  acting as (a matrix conjugate to) $A$ on homology.
  Recall that $\psi(\theta)$ has distance at most 1 from $f(\theta)$.
  It follows from Corollary~\ref{||:cor} that we can realize precisely
  those bundles $T_A$ such that $\|A\|\leqslant 1$, and $c_0=6$ on them. This shows the first
  assertion.

  If we also allow bricks of complexity 1, we can realize $T_A$ only
  as a self-assembling along a $\psi$ as above
  of some $(T\times I)_{\theta_0, \theta_1}$,
  which is itself either $B_0$ or an assembling of copies of $B_3$.
  Moreover
  $c_1((T\times I)_{\theta_0, \theta_1})=d(f(\theta_0),\theta_1)$
  by Proposition~\ref{T:times:I:prop}.
  So, viewing $A$ as a fixed map from $T\times\{0\}$ to $T\times\{1\}$, we have
  $$c_1(T_A)  =  6+\min\big\{
  d(f(\theta_0),\theta_1): \theta_i\in\Theta(T\times i),\
  d(A\theta_0,\theta_1)\leqslant 1\big\}.$$
  Corollary~\ref{||:cor} then shows that $c_1(T_A)=6+\max\{0,\|A\|-1\}$, and 
  the second assertion is proved.
\end{proof}

\paragraph{Computation of $c_3$}
The brick $B_4$ is a marked $D_2\times S^1$, where $D_2$ is the disc
with two holes. Let us now
fix a canonical homology basis $(a_i,b_i)$ on the $i$-th
component of $\partial (D_2\times S^1)$, where 
$a_i$ is a component of $\partial D_2$ with induced orientation,
$b_i$ is a copy of the fibre $S^1$, and the basis is positive.
Using these bases we will henceforth use rational numbers 
or $\infty$ to describe slopes.
Recalling also that we are denoting 
by $\theta^{(i)}$ the $\theta$-graph containing the slopes
$i$, $i+1$, and $\infty$, it is now easy to see from Fig.~\ref{B4:fig} that
$$B_4=(D_2\times S^1)_{\theta^{(0)},\theta^{(0)},
\theta^{(-1)}}.$$

Using Remark~\ref{few:cases:rem},
the next result shows that $n=8$ is well-behaved on tori.

\begin{prop}\label{c3:atoroidal:prop} 
\begin{enumerate}
\item\label{c3:atoroidal:tori:point} If $\alpha$ consists of
some slopes on $\partial (D_2\times S^1)$ and $(D_2\times S^1)_\alpha$ is 
$\Tsolid$ or $T\times I$ then $\alpha$ contains
a slope $0$ or $1/q$ for $q\in\matZ$;
\item\label{c3:atoroidal:degenerate:point} All slopes $0$ and $1/q$ for $q\in\matZ$ on any 
component of $\partial B_4$ are degenerate;
\item\label{c3:atoroidal:non:degenerate:point} A filling of $B_4$ along a
non-degenerate triple of slopes is either a genuine closed
atoroidal Seifert manifold or $T_{\tiny{\matr 01{-1}{-1}}}$;
\item\label{c3:atoroidal:three-fibres:point} If $0<q_i<p_i$ for $i=1,2,3$ and 
$t\geqslant -1$ then
\begin{equation}\label{c3:formula}
\hspace{-.2cm}c_3\big(S^2,(p_1,q_1),(p_2,q_2),(p_3,q_3),t\big)=
  |p_1,q_1|+|p_2,q_2|+|p_3,q_3|+t+1
\end{equation}
  except when $t=-1$, $q_1=q_2=q_3=1$ and $\{p_1,p_2,p_3\}$ equals either $\{2,2,2\}$,
  or $\{2,2,3\}$, or $\{2,3,6\}$, or $\{2,4,4\}$. For these cases the value of
  $c_3$ is respectively $2$, $3$, $6$, and $6$.
\end{enumerate}
\end{prop}
\begin{rem} 
{\em Let us explain why there are four special manifolds
excluded from the range of applicability of~(\ref{c3:formula}), but
$\big(S^2,(3,1),(3,1),(3,1),-1\big)=T_{\tiny{\matr 01{-1}{-1}}}$,
which is also quite special, is not excluded.
The first two manifolds excluded are the two closed bricks having complexity at most 3
(denoted by $C_{1,1}$ and $C_{1,2}$ in Proposition~\ref{Mstar:other:list:prop}).
Their complexity is 2 and 3 respectively,
and~(\ref{c3:formula}) is incorrect for them, because it would give 
3 and 4 respectively.
Turning to the torus bundles
$$T_{\tiny{\matr 01{-1}1}},\quad T_{\tiny{\matr 01{-1}0}},\quad
T_{\tiny{\matr 01{-1}{-1}}}$$
(see Proposition~\ref{Seifert:generalities:prop}), we know~\cite{MaPe,
Mat}
that their complexity is $6$, and we have seen in
Proposition~\ref{c1:torus:bundle:prop} that 
$c_1$ already equals $6$ on them, while 
formula~(\ref{c3:formula}) gives $8$, $7$, and $6$ respectively.
This means that the formula is incorrect for the first two of them, but it is
correct for the last one. It also means that $T_{\tiny{\matr 01{-1}{-1}}}$
    can be expressed both as a sharp self-assembling of one $B_0$, and as 
    a sharp assembling of one $B_4$ and three $B_2$'s. 
We have here one of the (not many) instances we know of 
non-uniqueness in Theorem~\ref{split:teo}.}
\end{rem}

\dimo{c3:atoroidal:prop}
If all slopes in $\alpha$ are different from $0$ and $1/q$ then
$\big(D_2\times S^1\big)_\alpha$ is Seifert fibred with $\#\alpha$
exceptional fibres, and point 1 easily follows.
Let us turn to point 2.
The slope $0$ is degenerate on each boundary component 
by Proposition~\ref{degenerate:prop}-(\ref{degenerate:contained:point}), 
because it is contained in
both $\theta^{(0)}$ and $\theta^{(-1)}$.  Let us prove that a slope
$1/q$ on the $j$-th boundary component is degenerate.
If we perform a $1/q$-Dehn filling on the $j$-th component we get some
$(T\times I)_{\theta_0,\theta_1}$.  
It is now easy to prove that:
\begin{eqnarray*}
  d(f(\theta_0),\theta_1) & = & \left\{\begin{array}{ll}
      |q| & \quad{\rm if}\ j=1,2\\
      |q-1|& \quad{\rm if}\ j=3\end{array}\right. \\
  d(1/q,\theta^{(i)}) & = & \left\{\begin{array}{ll}
      \max\{q-2,-q-1\} & \quad{\rm if}\ i=0\\
      \max\{q-1,-q-2\} & \quad{\rm if}\ i=-1\end{array}\right.
\end{eqnarray*}
(see also Corollary~\ref{distance:cor}).
In all cases we get 
$$c_1(T\times I,\{\theta_0,\theta_1\})=
d(f(\theta_0),\theta_1) <3 + d(1/q,\theta^{(i)})$$ 
by Proposition~\ref{T:times:I:prop}, and point 2 is proved.

To prove point 3, let us take a triple $\alpha$ of slopes not
containing $0$'s and $1/q$'s.
Then $(B_4)_\alpha$ is Seifert fibred over $S^2$ with three exceptional fibres, so it is
either genuine and atoroidal or one of the three exceptional torus bundles
of Proposition~\ref{Seifert:generalities:prop}-(1). We will show together with point 4 
that $\alpha$ is degenerate when $(B_4)_\alpha$ is
$T_{\tiny{\matr 01{-1}1}}$ or $T_{\tiny{\matr 01{-1}0}}$, and it is not 
when $(B_4)_\alpha$ is $T_{\tiny{\matr 01{-1}{-1}}}$.

So we turn to point 4, treating first the case of a genuine atoroidal manifold
$$M=\big(S^2,(p_1,q_1),(p_2,q_2),(p_3,q_3),t\big).$$
The shape of $B_0,\ldots,B_3$ and the knowledge of the closed bricks 
of small complexity
show that only for $M=C_{1,1}$ and $M=C_{1,2}$ it is 
possible to realize $M$ as an assembling of bricks with complexity at most
$3$ but without $B_4$'s.
Of course $c_3(M)=c(M)$ for $M=C_{1,1}$ and $M=C_{1,2}$, so we exclude these
two cases. By Proposition~\ref{atoroidal:prop} the conclusion will now
follow as soon as we prove that the minimum of
$$3+d(\alpha_1,\theta^{(0)})+d(\alpha_2,\theta^{(0)})+d(\alpha_3,\theta^{(-1)})$$
over all realizations of
$M$ as $(D_2\times S^1)_{\alpha_1,\alpha_2,\alpha_3}$ is
$|p_1,q_1|+|p_2,q_2|+|p_3,q_3|+t+1$.  
It is well-known that, up to permutation, the choices of
$\alpha_1,\alpha_2,\alpha_3$ leading to one such realization are
precisely those of the form
$$(\alpha_1,\alpha_2,\alpha_3)= \big(sp_1/(q_1+k_1\cdot p_1),
sp_2/(q_2+k_2\cdot p_2), sp_3/(q_3+k_3\cdot p_3)\big)$$
for $s\in\{1,-1\}$ and integers $k_1,k_2,k_3$ such that $k_1+k_2+k_3=t$.

Before computing the minimum we recall that 
$d(p/q,\theta^{(0)})=|p,q|-1$ and $d(p/q,\theta^{(-1)})=|p,p+q|-1$
by Corollary~\ref{distance:cor}. Moreover if $0<q<p$ we have
$|p,q-p|=|p,q|$ by Proposition~\ref{good:|p,q|:for:lens:prop}-(2).
So, using the definition of $|\cdot\,,\,\cdot|$ we have for $h>0$ that
\begin{eqnarray*}
|p,q+hp| & = & |p,q|+h,\\
|p,q-hp| & = & 1+|p,hp-q|=1+(h-1)+|p,p-q|=|p,q|+h,\\
&\Rightarrow & |p,q+k\cdot p|=|p,q|+|k|\qquad\forall k\in\matZ.
\end{eqnarray*}
Similarly:
\begin{eqnarray*}
|p,-q-hp| & = & 1+|p,q+hp|=1+h+|p,q|,\\
|p,-q+hp| & = & h-1+|p,p-q|=-1+h+|p,q|,\\
&\Rightarrow & |p,-q-k\cdot p|=|p,q|+|k+1|\qquad\forall k\in\matZ.
\end{eqnarray*}

Back to our minimum, let us first compute it over realizations with
$s=1$.  From the previous discussion we have
\begin{eqnarray*}
  3&+&d(\alpha_1,\theta^{(0)})+d(\alpha_2,\theta^{(0)})+d(\alpha_3,\theta^{(-1)}) \\
  & = & \big|p_1,q_1+k_1\cdot p_1\big|+ \big|p_2,q_2+k_2\cdot p_2\big|
  +\big|p_3,q_3+(1+k_3)\cdot p_3\big| \\
  & = & |k_1|+|k_2|+|1+k_3|+\sum_{i=1}^3|p_i,q_i|
\end{eqnarray*}
therefore the minimum is attained when $k_1=k_2=0$ and $k_3=t$, and
its value is $\sum_i|p_i,q_i|+t+1$.  Note that this value is unchanged
under reordering of the $p_i/q_i$'s, and that the assumption $t\geqslant -1$ is
used here.

Let us now consider the realizations with $s=-1$. From the properties
of $|\cdot\,,\cdot\,|$ shown above we deduce that
\begin{eqnarray*}
  3&+&d(\alpha_1,\theta^{(0)})+d(\alpha_2,\theta^{(0)})+d(\alpha_3,\theta^{(-1)}) \\
  & = & \big|p_1,-q_1-k_1\cdot p_1\big|+
  \big|p_2,-q_2-k_2\cdot p_2\big|+\big|p_3,-q_3-(k_3-1)\cdot p_3\big|\\
  & = & |1+k_1|+|1+k_2|+|k_3|+\sum_{i=1}^3|p_i,q_i|
\end{eqnarray*}
so the minimum is attained when $k_1=k_2=-1$ and $k_3=t+2$, and its
value is $\sum_i|p_i,q_i|+t+2$. Again this value is invariant under
reordering, and the conclusion of point 3 easily follows.

We are left to consider the three Seifert manifolds which are not genuine and atoroidal
because they can also be described as torus bundles.
For these manifolds $c_1$ equals $6$ by Proposition~\ref{c1:torus:bundle:prop}.
Now suppose there is an assembling which realizes $c_3(T_A)$ and 
includes some copy of $B_4$.
Then $c_3(T_A)$ is computed just as above in the atoroidal case, taking
a minimum over all fillings of $B_4$ that give $T_A$. 
This minimum is respectively $8$, $7$, and $6$ when
$A$ is ${\tiny{\matr 1{-1}10}}$, ${\tiny{\matr 01{-1}0}}$, and
${\tiny{\matr 01{-1}{-1}}}$. In the first two cases, we deduce that the 
triple giving $T_A$ is degenerate.  
In the last case,
the triple is non-degenerate.
\finedimo

\paragraph{Computation of $c_8$ and $c_9$}
Recall that there are two bricks $B_5$ and $B_6$ having complexity 8, and four
$B_7,\ldots,B_{10}$ having complexity 9.  The bricks $B_6,\ldots,
B_{10}$ are hyperbolic, and it will follow from
Propositions~\ref{non:neg-curved:degenerate:prop}
and~\ref{non:neg-curved:degenerate:prop2} that they give no
contribution to the values of $c_8$ and $c_9$ for non-hyperbolic manifolds
(see Proposition~\ref{no:Bi:for:Seifert:prop}). 
Therefore
we have $c_8=c_9$ on Seifert manifolds,
and we conclude this section by studying the contribution of the brick
$B_5$ to $c_8=c_9$ on the atoroidal ones. Let us first describe $B_5$.
We take $D_2\times S^1$ with the basis on the boundary defined above
for $B_4$, and again we use the rationals to denote slopes.  Then $B_5
= (D_2\times S^1)_{2,3,\theta}$ where $\theta$ is the $\theta$-graph containing
the slopes $-1,-6/5,-5/4$.

\begin{prop}\label{c8:atoroidal:prop}
\begin{enumerate}
\item A slope $\alpha$ on $B_5$ is non-degenerate if and only if
$-5/4<\alpha< -1$ and $\alpha\not\in\{-6/5,-7/6,-8/7\}$;
\item The fillings of $B_5$ along non-degenerate slopes are precisely
  the Seifert manifolds of the form
  $\big(S^2,(2,1),(3,1),(p,q),-1\big)$ with $p/q>5,\ p/q\not\in\{6,7,8\}$.
\item For $p/q>5, p/q\not\in\{6,7,8\}$ we have $c_8\big(S^2,(2,1),(3,1),(p,q),-1\big)=
  |p,q|+2.$
\end{enumerate}
\end{prop}

\begin{rem}\label{after:c8:atoroidal:rem}
\begin{itemize}
\item \emph{The formula for $c_8$ in point 3 indeed does
not apply for $p/q\in\{6,7,8\}$, because it would give 
$p/q+1$, but we can actually show that $c_8$ equals $p/q$ in these three cases:
for $p/q=6$ the corresponding manifold is 
$T_{\tiny{\matr 01{-1}1}}$, whose complexity is $6$, so
$c_8$ and $c_9$ also are $6$ (and we have seen in 
Proposition~\ref{c1:torus:bundle:prop} that $c_1$ already is $6$).    For
    $p/q=7,8$ we get the closed bricks $E_2, E_3\in\calM^*$, having 
complexity $7$ and $8$ respectively.}
  \item \emph{Comparison with Proposition~\ref{c3:atoroidal:prop} shows
    that $c_8=c_3-1$ on manifolds of the form
    $\big(S^2,(2,1),(3,1),(p,q),-1\big)$ for $p/q>5,\ p/q\not\in\{6,7,8\}$.
    }
\end{itemize}
\end{rem}

\dimo{c8:atoroidal:prop} 
We first show that the slopes outside $(-5/4,-1)$ are degenerate. Later we will 
prove that the slopes inside $(-5/4,-1)$, except $-6/5,-7/6,-8/7$, are
non-degenerate, and at the same we will describe the fillings they determine and 
the complexity $c_8$ of these fillings.

\begin{figure}\begin{center}
\mettifig{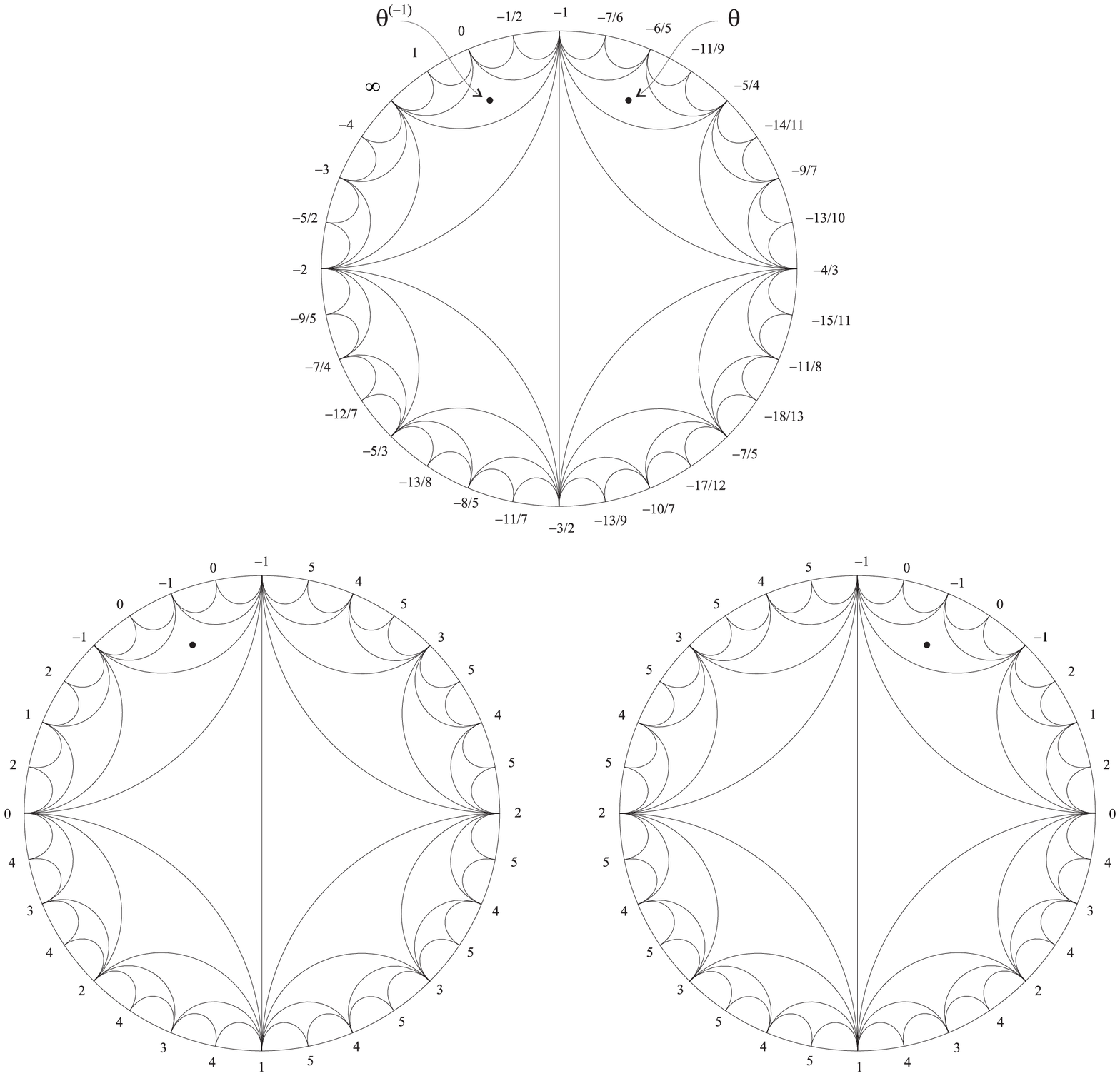, width=14cm}
\nota{Position of $\theta^{(-1)}$ and $\theta$, and computation of
      $d(\alpha,\theta^{(-1)})$ and $d(\alpha,\theta)$.}
\label{c8:fig}
\end{center}\end{figure}

Let us then fix $\alpha\in\partial_\matQ\matH^2$ and define $M=(B_5)_\alpha$,
\emph{i.e.}~$M=(D_2\times S^1)_{2,3,\alpha}$. We recall that 
$B_4=(D_2\times S^1)_{\theta^{(0)},\theta^{(0)},\theta^{(-1)}}$
and $c(B_4)=3$. Hence, by Proposition~\ref{atoroidal:prop},
$$c_8(M)\leqslant 3+0+1+d(\alpha,\theta^{(-1)}).$$
A sufficient condition for $\alpha$ to be degenerate is then that
\begin{equation}\label{suff:c8:deg:eqn}
4+d(\alpha,\theta^{(-1)})<8+d(\alpha,\theta)
\end{equation}
(where $\theta$ is the $\theta$-graph on $\partial B_5$).
Figure~\ref{c8:fig}, which uses an alternative embedding of the Farey tessellation,
easily proves that~(\ref{suff:c8:deg:eqn}) holds if $\alpha\not\in[-5/4,-1)$.
So all the slopes outside $[-5/4,-1)$ are degenerate, but 
$-5/4$ also is, because it is contained in $\theta$
(see Proposition~\ref{degenerate:prop}), and our first task is accomplished.

We can now concentrate on the slopes inside $(-5/4,-1)$, which are precisely
those of the form $p/(q-p)$ for $0<q<p$ and $p/q>5$. Note that $M$ is then
given by
$$M=\big(S^2,(2,1),(3,1),(p,q-p)\big)=
\big(S^2,(2,1),(3,1),(p,q),-1\big)$$
and the latter expression is the unique normalized one.
Now $\alpha=-6/5$ is degenerate because it is contained in $\theta$.
For $\alpha=-7/6$ we get $M=E_2$ and $c(M)=7<8=8+d(\alpha,\theta)$, so
$\alpha$ is again degenerate.  For $\alpha=-8/7$ we get $M=E_3$
and $c(M)=8$, and $\alpha$ is degenerate for the same reason.

Now take $\alpha$ in $(-5/4,-1)\setminus\{-6/5,-7/6,-8/7\}$. Then 
$M$ is not a closed brick of complexity up to $8$. Hence, by
Propositions~\ref{atoroidal:prop} and~\ref{c3:atoroidal:prop}
we have
$$c_8(M)=\min\{c_3(M), c(B_5)+d(\alpha,\theta)\} = 
\min\{3+|p,q|, 8+d(\alpha,\theta)\}$$
where $\alpha = p/(q-p)$. Corollary~\ref{distance:cor} gives
$3+|p,q|=4+d(\alpha,\theta^{(-1)})$ and Fig.~\ref{c8:fig} gives
$d(\alpha,\theta) = d(\alpha,\theta^{(-1)})-5$, hence
$c_8(M)=c_3(M)-1=|p,q|+2$, and our argument is complete.
\finedimo

\section{Complexity of hyperbolic manifolds}\label{hyperbolic:complexity:section}
As we know from Theorem~\ref{finiteness:teo}, for a closed hyperbolic $M$ 
we have $c_9(M)<\infty$ if and only if $M$ is a Dehn filling of $M6_1^3$.
We denote in this section this manifold by $N$, and we carefully analyze the complexity
of its fillings. In particular, 
we compute $c_8$ and $c_9$ on the hyperbolic fillings 
(proving Theorem~\ref{hyp:Dehn:fill:teo}) and we show that all the other fillings
are degenerate. The latter result will be used in Section~\ref{Seifert:section}
to conclude the computation of $c_n$ on Seifert manifolds. 

Recall that $N$ is a link complement in $S^3$ and that the slopes on $\partial N$
are described by rational numbers or $\infty$ using the meridian-longitude
homology bases.

\paragraph{Computation of $c_8$} 
The non-closed bricks of complexity 8 are the Seifert brick $B_5$
already discussed above and $B_6 = N_{1,-4,\theta^{(-1)}}$, a marking of
the ``figure-8 knot sister'' $N_{1,-4}=M2_2^1$. 
We start with the following easy result, which implies the first half
of Theorem~\ref{hyp:Dehn:fill:teo}:

\begin{prop}\label{hyperbolic:c8:prop}
Let $M$ be closed hyperbolic.
Then $c_8(M)$ equals the minimum of
  $7+|p,-q|$ over all realizations of $M$ as $N_{1,-4,p/q}$.
\end{prop}

\begin{proof}
By Proposition~\ref{atoroidal:prop}, $c_8(M)$ is the minimum of
$c(M'_\theta)+d(\alpha,\theta)$ with $M'_\theta\in\calB_{\leqslant 8}$ and $M=M'_\alpha$.
Since no filling of a Seifert manifold is hyperbolic,
we deduce that $M'_\theta$ must be $B_6$. 
So $c_8(M)$ is the minimum of $8+d(p/q,\theta^{(-1)})$
over all realizations of $M$ as $N_{1,-4,p/q}$. Now 
$d(p/q,\theta^{(-1)})=|p+q,q|-1$ by Corollary~\ref{distance:cor},
and it is easy to see that $|p+q,q|=|p,-q|$, using 
Proposition~\ref{lines:and:|p,q|:prop} and the matrix $S_1$ 
of Section~\ref{statement:section},
which fixes $\theta^{(0)}$ and maps $(p+q)/q$ to $-p/q$.
\end{proof}

\paragraph{Special fillings}
There is one more crucial ingredient that we must introduce
before computing $c_9$ and characterizing the degenerate fillings of $B_6,\ldots,B_{10}$.
This ingredient is the 
complete classification, established in~\cite{MaPe:chain} and highlighted in
Theorem~\ref{short:chain:fill:teo} above,
of the closed non-hyperbolic fillings of $N$. We need here only a simplified version
of~\cite[Theorem~1.3]{MaPe:chain}, which is 
Proposition~\ref{recognized:fillings:prop} below.
Later we will need to refer to~\cite[Theorem~1.3]{MaPe:chain}
only for the 14 exceptional filling triples listed in the last point
of Theorem~\ref{short:chain:fill:teo}. 

Before proceeding we must introduce
some notation. As we did before the statement of Proposition~\ref{c3:atoroidal:prop} for $i=2$,
we denote by $D_i$ the disc with $i$ holes, and we take product positive homology
bases on the components of $\partial (D_i\times S^1)$.
We use these bases to describe the filling
coefficients on $D_i\times S^1$
as rational numbers or $\infty$, and to describe gluing maps of boundary tori 
as matrices in $\GL_2(\matZ)$.

\begin{prop}[\cite{MaPe:chain}] \label{recognized:fillings:prop}
The following homeomorphisms hold:
\begin{eqnarray}
  N_{\infty,\frac rs,\frac tu} & = & \big(D_1\times S^1\big)_{\frac rs,-\frac ut}\label{521:eqn}
  \\
  N_{-3,\frac rs,\frac tu} & = & \big(D_2\times S^1\big)_{2,\frac{r+s}{r+2s}}
    \bigcup\nolimits_{\tiny{\matr110{-1}}}
  \big(D_2\times S^1\big)_{2,\frac{t+u}{t+2u}}\label{522:eqn} \\
  N_{-2,\frac rs,\frac tu} & = & \big(D_2\times S^1\big)_{\frac 32,-2-\frac rs,-2-\frac
  tu}\label{523:eqn} \\
  N_{-1,\frac rs,\frac tu} & = & \big(D_2\times S^1\big)_{2,-3-\frac rs,-3-\frac tu}
  \label{524:eqn} \\
  N_{0,\frac rs,\frac tu} & = & \big(D_2 \times S^1\big)_{\frac s{r+2s},\frac u{t+2u}}
    \bigcup\nolimits_{\tiny{\matr 0{-1}11}}
  \big(D_2\times S^1\big)_{2,3} \label{525:eqn} \\
  N_{1,1,\frac tu} & = & \big(D_2\times
  S^1\big)_{-\frac u{t+2u}}\Big/_{\tiny{\matr1{-1}{-1}0}}. \label{526:eqn} 
\end{eqnarray}
\end{prop}

\paragraph{Complexity estimates}
We prove now certain complexity
estimates needed below. Our statement involves the norm $|A|$ defined
in Section~\ref{statement:section} for $A\in\GL_2(\matZ)$.

\begin{lemma} \label{c3:graph:lem}
\begin{enumerate}
\item If $\alpha_1$ and $\alpha_2$ are slopes or $\theta$-graphs on distinct boundary
components of $D_2\times S^1$ then
$$c_3\big((D_2\times S^1)_{\alpha_1,\alpha_2,\theta^{(-1)}}\big)\leqslant 3+
d(\theta^{(0)},\alpha_1)+d(\theta^{(0)},\alpha_2).$$
\item Let $\tiny\matr ijhk\in\GL_2(\matZ)$. Assume $(\alpha_1,\alpha_2)$ and
$(\alpha_3,\alpha_4)$ are pairs as in the previous point. Then:
$$\hspace{-12pt}c_3\left((D_2\times S^1)_{\alpha_1, \alpha_2}
\bigcup\nolimits_{\tiny{\matr ijhk}}(D_2\times S^1)_{\alpha_3,\alpha_4}\right)
\leqslant 6+\sum_{l=1}^4 d(\alpha_l,\theta^{(0)})+ \left|{\tiny\matr i{-j}{-h}k}\right|.$$
\item If $\alpha$ is a slope or a $\theta$-graph on $\partial (D_2\times S^1)$ then
$$c_3\left( (D_2\times S^1)_{\alpha}\Big/_{\tiny{\matr1{-1}{-1}0}}\right)
\leqslant 9+d(\alpha,\theta^{(-1)}).$$
\end{enumerate}
\end{lemma}

\begin{proof}
Point 1 readily follows from Proposition~\ref{easy:assemblings:prop},
using $B_4=(D_2\times S^1)_{\theta^{(0)}, \theta^{(0)}, \theta^{(-1)}}$, because $c(B_4)=3$.

For point 2 we apply the inequality of point 1 to both the blocks involved, and
we change the homology bases on the boundary of both blocks, performing
a change of type $(a,b)\mapsto(-a,b)$. In the new basis the old 
$\theta$-graph $\theta^{(-1)}$ is transformed into $\theta^{(0)}$, and the gluing matrix becomes
$A=\tiny\matr i{-j}{-h}k$. Now Proposition~\ref{easy:assemblings:prop} implies that
the gluing matching $\theta^{(0)}$ to $A\theta^{(0)}$ can be realized by insertion of 
$d(\theta^{(0)},A\theta^{(0)})$ copies of $B_3$ and assembling. But 
$d(\theta^{(0)},A\theta^{(0)})=|A|$ by Corollary~\ref{||:cor}
and the conclusion follows.

Turning to point 3, we note that 
$c_3(D_2\times S^1)_{\alpha, \theta^{(0)}, \theta^{(0)}}\leqslant 3+d(\alpha,\theta^{(-1)})$,
again from Proposition~\ref{easy:assemblings:prop}
applied to $B_4$. Now the gluing map with matrix 
$\tiny{\matr 1{-1}{-1}0}$ matches $\theta^{(0)}$ to $\theta^{(-1)}$, and 
$d(\theta^{(0)},\theta^{(-1)})=1$, so the map can be realized via a self-assembling.
We deduce the upper bound $6+3+d(\alpha,\theta^{(-1)})$ on $c_3$,
and the proof is complete.
\end{proof}

\paragraph{Degenerate fillings of $B_6$}
Now we turn to the 
discussion of non-hyperbolic
fillings of $N_{1,-4}$. The following result will
be used in Propositions~\ref{no:Bi:for:Seifert:prop} 
and~\ref{cn:torus:bundle:prop} to show  that $B_6$ gives no
contribution to $c_8(M)$ if $M$ is a Seifert manifold or a torus bundle.

\begin{prop} \label{non:neg-curved:degenerate:prop}
Every slope on $B_6=N_{1,-4,\theta^{(-1)}}$ 
producing a non-hyperbolic manifold 
is degenerate.
\end{prop}

\begin{proof}
Corollary~\ref{B7:fill:cor} shows that the slopes that me must prove to be degenerate are
precisely $\infty, -3, -2, -1, -1/2, -1/3, 0, 1$.
The slopes $\infty, -1, 0$ are degenerate by 
Proposition~\ref{degenerate:prop}-(\ref{degenerate:contained:point}).
It follows from~(\ref{515:eqn}) and~(\ref{522:eqn}) 
that 

$$N_{1,-4,-3}=N_{1,-4,-1/3} = \big(D_2\times S^1\big)_{2,2}
\bigcup\nolimits_{\tiny{\matr 01{-1}{-1}}}
\big(D_2\times S^1\big)_{2,3}.$$
Lemma~\ref{c3:graph:lem}-(2) then implies that 
\begin{equation}\label{c3:ineq:used:twice}
\begin{array}{c}
c_3(N_{1,-4,-3})=c_3(N_{1,-4,-1/3})
\leqslant 6+0+0+0+1+0= 7 \\
<9=c(B_6)+d(-3,\theta^{(-1)})=c(B_6)+d(-1/3,\theta^{(-1)})
\end{array}
\end{equation}
so $-3$ and $-1/3$ are degenerate. For the slopes $-2$ and $-1/2$ we use 
the following equality, again deduced from~(\ref{515:eqn}) and~(\ref{523:eqn})
or~\cite[Theorem~1.3]{MaPe:chain}:
$$N_{1,-4,-2}=N_{1,-4,-1/2}=\big(S^2,(3,2),(3,2),(2,1),-1\big).$$
We can then apply 
Proposition~\ref{c3:atoroidal:prop}-(\ref{c3:atoroidal:three-fibres:point}) 
to see that 
$$c_3(N_{1,-4,-2})=c_3(N_{1,-4,-1/2})=
2+2+1+(-1)+1=5<8=c(B_6).$$
Therefore
$-2$ and $-1/2$ are also degenerate. We are left to deal with the slope $1$.
By~(\ref{526:eqn}) or~\cite[Theorem~1.3]{MaPe:chain}
we have $N_{1,-4,1}=T_A$ with
$A={\tiny{\matr {-3}{-1}10}}$. Now $\|A\|\leqslant 
d(\theta^{(-1)},A\theta^{(-1)})$ by Corollary~\ref{||:cor},
but $A\theta^{(-1)}=\theta^{(-3)}$, so $\|A\|\leqslant 2$, and
$c_1(T_A)\leqslant 7$ by Proposition~\ref{c1:torus:bundle:prop}.
Since $7<8=c(B_6)+d(1,\theta^{(-1)})$, we deduce that $1$ is 
also degenerate, and the proof is complete.
\end{proof}

\paragraph{Computation of $c_9$} 
The bricks $B_6$, $B_7$, $B_8$, $B_9$, and $B_{10}$ are
respectively $M2_2^1$, $M3_4^1$, $M4_1^2$, $N=M6_1^3$, and $N$
again (with appropriate markings; as above we use the
notation of~\cite{Ca-Hi-We} for cusped hyperbolic manifolds). 
Each of these  
manifolds can be obtained as an appropriate Dehn filling of $N$.
So for $i=6,\ldots,10$ we can describe $B_i$ as $N_{\alpha^{(i)}}$
where $\alpha^{(i)}$ is a triple of slopes and/or $\theta$-graphs.  Due to the
symmetries of $N$ the triples $\alpha^{(i)}$ are not unique, but we
choose the following definite ones:
$$\alpha^{(6)}=\big(1,-4,\theta^{(-1)}\big),\qquad \alpha^{(7)} =
\big(1,-5,\theta^{(-1)}\big), \qquad \alpha^{(8)} =
\big(1,\theta^{(-2)},\theta^{(-2)}\big),$$
$$\alpha^{(9)} = \big(\theta^{(-2)},\theta^{(-2)},\theta^{(-2)}\big), \qquad
\alpha^{(10)} = \big(\theta^{(-2)},\theta^{(-2)},\theta^{(-3)}\big).$$
The
proof that indeed $B_i=N_{\alpha^{(i)}}$ for $i=6,\ldots, 10$ uses the
explicit description of the $B_i$'s given in~\cite{MaPe}.

As opposed to what we did for $c_8$, we start for $c_9$ with the non-hyperbolic
fillings of $N$, because we have not proved yet well-behavedness on tori
of $n=9$. This property is readily implied by the next 
result, using Remark~\ref{few:cases:rem} and the
information already shown on $B_4$ in Proposition~\ref{c3:atoroidal:prop}.
The next result will also be used in Propositions~\ref{no:Bi:for:Seifert:prop} 
and~\ref{cn:torus:bundle:prop}
to show that $c_9$ equals $c_8$ on Seifert manifolds and torus bundles.

\begin{prop} \label{non:neg-curved:degenerate:prop2}
If $B\in\{B_7,\ldots,B_{10}\}$ and $\alpha$ is a set of slopes on 
$\partial B$ such that $B_\alpha$ is not hyperbolic then $\alpha$
is degenerate.
\end{prop}

\begin{proof}
Recall that the triples $\alpha$ such that $N_\alpha$ is not hyperbolic
were described in Theorem~\ref{short:chain:fill:teo}.
Since $c(B_i)>2$ for $i=7,\ldots,10$ we can now apply
Proposition~\ref{degenerate:prop}-(\ref{degenerate:superset:point}).
According to this result, to prove our proposition
it is sufficient to show degeneracy with respect to $B_i$ of the $\alpha$'s as follows:
\begin{itemize}
\item[(a)] $\alpha$ is a triple of slopes and $\theta$-graphs
obtained from $\alpha^{(i)}$ by removing some $\theta$-graphs 
and inserting slopes instead of them;
\item[(b)] $\alpha$ satisfies one of the conditions of Theorem~\ref{short:chain:fill:teo};
\item[(c)] it is impossible to attain one of the conditions of 
Theorem~\ref{short:chain:fill:teo}
with a smaller number of $\theta$-graph-to-slope replacements on $\alpha^{(i)}$.
\end{itemize}
To prove that all such $\alpha$'s are degenerate, we first establish upper bounds
on $c_9(N_\alpha)$ for all $\alpha$'s obtained as in (a)-(c) from $\alpha^{(9)}$.
Our inequalities use many facts
already established, together with the easy relation
\begin{equation}\label{special:matrix:inequality:eqn}
{\tiny\left\|\left(\begin{array}{cc} n & -1 \\ 1 & 0 \\ \end{array}\right)\right\|} 
\leqslant d\left(\theta^{(-1)},
{\tiny \left(\begin{array}{cc} n & -1 \\ 1 & 0 \\ \end{array}\right)}\theta^{(-1)}\right)=
d(\theta^{(-1)},\theta^{(n)})=|n+1|
\end{equation}
(already used with $n=-3$ in the proof of Proposition~\ref{non:neg-curved:degenerate:prop}).
To save space while helping the reader follow our arguments, we will
write under each equality or inequality sign the equation or statement it is implied by,
with ''\cite{MaPe:chain}´´ standing for \cite[Theorem~1.3]{MaPe:chain}.

\newcommand{\equalbecause}[1]{\mathop{=}\limits_{#1}}
\newcommand{\smallerbecause}[1]{\mathop{\leqslant}\limits_{#1}}

\begin{equation}\label{infty:ineq}
{\hspace{-.1cm}}\begin{array}{rcl}
c_9(N_{\infty,\theta^{(-2)},\theta^{(-2)}}) 
& \equalbecause{(\ref{521:eqn})} & 
c_9\big(T\times I,\theta^{(-2)}\times\{0\},\theta(-1,- 1/2,0)\times\{1\}\big)\\
& \smallerbecause{{\rm Prop.\ }\ref{easy:assemblings:prop}} & 
d(\theta^{(-2)},\theta(-1,-1/2,0))=2
\end{array}
\end{equation}
(here $\theta(-1,-1/2,0)$ denotes the $\theta$-graph containing the slopes
$-1$, $-1/2$, and $0$).
\begin{equation}\label{minusthree:ineq}
{\hspace{-.4cm}}\begin{array}{rcl}
c_9(N_{-3,\theta^{(-2)},\theta^{(-2)}}) 
& \equalbecause{(\ref{522:eqn})} & 
c_9\left((D_2\times S^1)_{2,\theta^{(0)}}
    \bigcup\nolimits_{\tiny{\matr 110{-1}}}
    (D_2\times S^1)_{2,\theta^{(0)}}\right) \\
& \smallerbecause{{\rm Lem.\ \ref{c3:graph:lem}(2)}} & 
6+0+0+0+0+0=6.
\end{array}
\end{equation}

\begin{equation}\label{minustwo:ineq}
c_9(N_{-2,\theta^{(-2)},\theta^{(-2)}}) 
\equalbecause{(\ref{523:eqn})}
c_9\big((D_2\times S^1)_{3/2,\theta^{(-1)},\theta^{(-1)}}\big)
\smallerbecause{{\rm Lem.\ \ref{c3:graph:lem}(1)}}
3+1+1=5.
\end{equation}

\begin{eqnarray*}
c_9(N_{-1,\theta^{(-2)},\theta^{(-2)}}) 
& \equalbecause{(\ref{524:eqn})} &
 c_9\big((D_2\times S^1)_{2,\theta^{(-2)},\theta^{(-2)}}\big)\\
& \smallerbecause{{\rm Prop.\ }\ref{easy:assemblings:prop}} & 
c_9\big((D_2\times S^1)_{2,\theta^{(-2)},\theta^{(-1)}}\big)+
d(\theta^{(-2)},\theta^{(-1)}) \\
& \smallerbecause{{\rm Lem.\ \ref{c3:graph:lem}(1)}} &
3+0+2+1=6.
\end{eqnarray*}

\begin{equation}\label{c9:zero:tt:ineq}
{\hspace{-.7cm}}\begin{array}{rcl}
c_9(N_{0, \theta^{(-2)},\theta^{(-2)}}) 
& \equalbecause{(\ref{525:eqn})} & 
c_9\left((D_2\times S^1)_{\theta^{(0)},\theta^{(0)}}
    \bigcup\nolimits_{\tiny{\matr 0{-1}11}}
    (D_2 \times S^1)_{2,3}\right) \\
& \smallerbecause{{\rm Lem.\ \ref{c3:graph:lem}(2)}} &
6+0+0+0+1+0=7.
\end{array}
\end{equation}

\begin{equation}\label{c9:11:t:ineq}
c_9(N_{1,1,\theta^{(-2)}}) 
\equalbecause{(\ref{526:eqn})} 
c_9\left((D_2\times S^1)_{\theta^{(-1)}}\Big/_{\tiny{\matr1{-1}{-1}0}}\right)
\smallerbecause{{\rm Lem.\ \ref{c3:graph:lem}(3)}} 9+0=9.
\end{equation}

\begin{eqnarray*}
c_9(N_{-4,-1/2,\theta^{(-2)}}) 
& \equalbecause{(\ref{511right:eqn})^{-1}} & 
c_9(N_{-3/2,0,\theta^{(-2)}}) \\
& \smallerbecause{{\rm Prop.\ }\ref{easy:assemblings:prop}} &
c_9(N_{0,\theta^{(-2)},\theta^{(-2)}})+d(-3/2,\theta^{(-2)})
\smallerbecause{(\ref{c9:zero:tt:ineq})} 7+0=7.
\end{eqnarray*}

\begin{eqnarray*}
c_9(N_{-5/2,-3/2,\theta^{(-2)}})
& \equalbecause{(\ref{512:eqn})} & 
c_9(N_{-5/2,-3,\theta^{(-2)}}) \\
& \smallerbecause{{\rm Prop.\ }\ref{easy:assemblings:prop}} & 
c_9(N_{\theta^{(-2)},-3,\theta^{(-2)}})+d(-5/2,\theta^{(-2)}) 
\smallerbecause{(\ref{minusthree:ineq})} 6+1= 7.
\end{eqnarray*}

$$c_9(N_{-5,-5,-1/2})
\equalbecause{(\ref{513:eqn})}
c_9(N_{1,1,-1/2})
\smallerbecause{{\rm Prop.\ }\ref{easy:assemblings:prop}}
c_9(N_{1,1,\theta^{(-2)}})+d(-1/2,\theta^{(-2)})
\smallerbecause{(\ref{c9:11:t:ineq})}
9+1=10.$$

\begin{equation}\label{citami:subito:eqn}
\begin{array}{rcl}
c_9(N_{-4,-4,-2/3})
& \equalbecause{(\ref{511right:eqn})} &
c_9(N_{-3/2,1,1}) \\
& \smallerbecause{{\rm Prop.\ }\ref{easy:assemblings:prop}} &
c_9(N_{1,1,\theta^{(-2)}})+d(-3/2,\theta^{(-2)})
\smallerbecause{(\ref{c9:11:t:ineq})}9+0=9.
\end{array}
\end{equation}

$$c_9(N_{-3/2,-7/3,-7/3})
\equalbecause{(\ref{511left:eqn})}
c_9(N_{-3/2,1,1})
\smallerbecause{(\ref{citami:subito:eqn})} 9.
$$

$$c_9(N_{-4,-3/2,-3/2})
\equalbecause{\rm \cite{MaPe:chain}}
c_9\left(T_{\tiny\matr {-3}{-1}10}\right)
\smallerbecause{{\rm Prop.\ }\ref{c1:torus:bundle:prop}}
5+\left\|{\tiny\matr {-3}{-1}10}\right\|
\smallerbecause{(\ref{special:matrix:inequality:eqn})}
5+2=7.$$

$$c_9(N_{-8/3,-3/2,-3/2})
\equalbecause{(\ref{511right:eqn})}
c_9(N_{-4,-1/3,1})
\equalbecause{(\ref{515:eqn})}
c_9(N_{1,-4,-3})
\smallerbecause{(\ref{c3:ineq:used:twice})} 7.$$

$$c_9(N_{-5/2,-5/2,-4/3})
\equalbecause{(\ref{512:eqn})} 
c_9(N_{-5/2,1,1})
\smallerbecause{{\rm Prop.\ }\ref{easy:assemblings:prop}}
c_9(N_{1,1,\theta^{(-2)}})+d(-5/2,\theta^{(-2)})
\smallerbecause{(\ref{c9:11:t:ineq})}9+1=10.$$

\begin{eqnarray*}
c_9(N_{-5/2,-5/3,-5/3})
& \equalbecause{\rm \cite{MaPe:chain}} & 
c_9\left((D_2\times S^1)_{2,2}
    \bigcup\nolimits_{\tiny{\matr 1{-1}01}}
    (D_2\times S^1)_{2,3}\right) \\
& \smallerbecause{{\rm Lem.\ \ref{c3:graph:lem}(2)}} &
6+0+0+0+1+1=8. 
\end{eqnarray*}

$$c_9(N_{1,2,2})
\equalbecause{\rm \cite{MaPe:chain}}
c_9\big(S^2,(2,1),(3,1),(7,1),-1\big)=7$$
(the last equality holds because $\big(S^2,(2,1),(3,1),(7,1),-1\big)$ is the
closed brick $E_2$ of $\calM^*$, and its complexity is $7$).

$$c_9(N_{1,2,3})
\equalbecause{\rm \cite{MaPe:chain}}
c_9\big(S^2,(2,1),(4,1),(5,1),-1\big)
\smallerbecause{{\rm Prop.\ \ref{c3:atoroidal:prop}-(\ref{c3:atoroidal:three-fibres:point})}}
1+3+4+(-1)+1=8.$$ 

$$c_9(N_{1,2,4})
\equalbecause{\rm \cite{MaPe:chain}}
c_9\big(S^2,(3,1),(3,1),(4,1),-1\big)
\smallerbecause{{\rm Prop.\ \ref{c3:atoroidal:prop}-(\ref{c3:atoroidal:three-fibres:point})}}
2+2+3+(-1)+1=7.$$

\begin{eqnarray*}
c_9(N_{1,2,5})
& \equalbecause{\rm \cite{MaPe:chain}} & 
c_9\left((D_2\times S^1)_{2,2}
    \bigcup\nolimits_{\tiny{\matr 0110}}
    (D_2\times S^1)_{2,3}\right) \\
& \smallerbecause{{\rm Lem.\ \ref{c3:graph:lem}(2)}} &
6+0+0+0+1+0=7.
\end{eqnarray*}

\begin{eqnarray*}
c_9(N_{1,3,3})
& \equalbecause{\rm \cite{MaPe:chain}} & 
c_9\left((D_2\times S^1)_{2,3}
    \bigcup\nolimits_{\tiny{\matr 120{-1}}}
    (D_2\times S^1)_{2,2}\right) \\
& \smallerbecause{{\rm Lem.\ \ref{c3:graph:lem}(2)}} &
6+0+1+0+0+1=8.
\end{eqnarray*}

\begin{equation}\label{222:ineq}
\begin{array}{rcl}
c_9(N_{2,2,2})
& \equalbecause{\rm \cite{MaPe:chain}} &
c_9\left((D_2\times S^1)_{2,3}
    \bigcup\nolimits_{\tiny{\matr 23{-1}{-2}}}
    (D_2\times S^1)_{2,2}\right) \\
& \smallerbecause{{\rm Lem.\ \ref{c3:graph:lem}(2)}} &
6+0+1+0+0+2=9.
\end{array}
\end{equation}

Using these inequalities, we prove now that the triples $\alpha$ 
obtained as explained in (a)-(c) from $\alpha^{(i)}$
are degenerate with respect to $B_i$ for $i=7,\ldots,10$.
We start with $i=7$. Since $\alpha^{(7)}=(1,-5,\theta^{(-1)})$
the only triples we have to consider are $(1,-5,n)$ for $n=-3,-2,-1,0,1,\infty$.
Now $-1,0,\infty$ are degenerate because they are contained in $\theta^{(-1)}$.
And the slopes $-3$, $-2$, and $1$ are degenerate because

\begin{eqnarray*}
c_9(N_{1,-5,-3})
& \smallerbecause{{\rm Prop.\ }\ref{easy:assemblings:prop}} &
c_9(N_{\theta^{(-2)},\theta^{(-2)},-3})+d(\theta^{(-2)},1)+d(\theta^{(-2)},-5)\\
& \smallerbecause{(\ref{minusthree:ineq})} & 6+1+2=9
< 10=9+1=c(B_7)+d(-3,\theta^{(-1)})
\end{eqnarray*}

\begin{eqnarray*}
c_9(N_{1,-5,-2})
& \smallerbecause{{\rm Prop.\ }\ref{easy:assemblings:prop}} &
c_9(N_{\theta^{(-2)},\theta^{(-2)},-2})+d(\theta^{(-2)},1)+d(\theta^{(-2)},-5)\\
& \smallerbecause{(\ref{minustwo:ineq})} & 5+1+2=8
<  9=9+0=c(B_7)+d(-2,\theta^{(-1)})
\end{eqnarray*}

\begin{eqnarray*}
c_9(N_{1,1,-5})
& \equalbecause{\rm \cite{MaPe:chain}} &
c_9\left(T_{\tiny{\matr 41{-1}0}}\right)
\smallerbecause{{\rm Prop.\ }\ref{c1:torus:bundle:prop}}
5+\|{\tiny{\matr 41{-1}0}}\|\\
& \smallerbecause{(\ref{special:matrix:inequality:eqn})} & 5+3=8 
< 9=9+0=c(B_7)+d(1,\theta^{(-1)}).
\end{eqnarray*}

To conclude for $i=8,9,10$ we now claim that for all the triples
$\alpha$ obtained as explained in (a)-(c) from $\alpha^{(9)}$ the following
inequality holds:
\begin{equation}\label{strong:upper:bound}
c_9(N_\alpha)<8+d(\alpha,\alpha^{(9)}).
\end{equation}
Before proving~(\ref{strong:upper:bound}), we show that it
is sufficient to conclude not only for $i=9$, which is obvious since
$c(B_9)=9$, but also for
$i=8$ and $i=10$. We start with $i=10$ and assume that 
$\alpha$ is obtained as in (a)-(c) from 
$\alpha^{(10)}$ by some $\theta$-graph-to-slope replacement. Let $\tilde\alpha$ be 
obtained by performing the same replacements on $\alpha^{(9)}$. 
Then~(\ref{strong:upper:bound}) applies to $\tilde\alpha$.
Using the fact that 
$d(\alpha^{(9)},\alpha^{(10)})=d(\theta^{(-2)},\theta^{(-3)})=1$,
it is also easy to see that either $d(\alpha,\tilde\alpha)=0$ and 
$d(\tilde\alpha,\alpha^{(9)})\leqslant d(\alpha,\alpha^{(10)})+1$,
or $d(\tilde\alpha,\alpha)=1$ and 
$d(\tilde\alpha,\alpha^{(9)})=d(\alpha,\alpha^{(10)})$. 
In particular, the following holds anyway:
\begin{equation}\label{usami:subito:bis:ineq}
d(\tilde\alpha,\alpha^{(9)})+d(\tilde\alpha,\alpha)
\leqslant d(\alpha,\alpha^{(10)})+1.
\end{equation}
Hence
\begin{eqnarray*}
c_9(N_\alpha)
& \smallerbecause{{\rm Prop.\ }\ref{easy:assemblings:prop}} &
c_9(N_{\tilde\alpha})+d(\tilde\alpha,\alpha)
\mathop{<}\limits_{(\ref{strong:upper:bound})}
8+d(\tilde\alpha,\alpha^{(9)})+d(\tilde\alpha,\alpha) \\
& \smallerbecause{(\ref{usami:subito:bis:ineq})} & 
9+d(\alpha,\alpha^{(10)})=c(B_{10})+d(\alpha,\alpha^{(10)}).
\end{eqnarray*}
For $i=8$ a very similar argument applies, except that 
to define $\tilde\alpha$ it is necessary in some cases 
to replace $\theta^{(-2)}$ by $1$, to ensure that one of the
properties (A), (B), or (C) is fulfilled. The key point is in any 
case that $d(\alpha^{(8)},\alpha^{(9)})=d(1,\theta^{(-2)})=1$.

Having shown that~(\ref{strong:upper:bound}) allows to conclude the proof, we then establish it.
Recall that
the $\alpha$'s to
consider are precisely those appearing in the left-hand side of the
inequalities~(\ref{infty:ineq}) to~(\ref{222:ineq}), and there is nothing to prove
whenever $d(\alpha,\alpha^{(9)})\geqslant 0$ and the
right-hand side is $7$ or less.
In the other cases we get the result by computing
$d(\alpha,\alpha^{(9)})$ and showing that it is greater than $-8$ plus the number on 
the right-hand side:
\begin{eqnarray*}
d((\infty,\theta^{(-2)},\theta^{(-2)}),\alpha^{(9)}) & = & -1 >2-8 \\
d((-2,\theta^{(-2)},\theta^{(-2)}),\alpha^{(9)}) & = & -1 >5-8 \\
d((-1,\theta^{(-2)},\theta^{(-2)}),\alpha^{(9)}) & = & -1 >6-8 \\
d((1,1,\theta^{(-2)}),\alpha^{(9)}) & = & 2 > 9 -8 \\
d((-5,-5,-1/2),\alpha^{(9)}) & = & 5 > 10 -8 \\
d((-4,-4,-2/3),\alpha^{(9)}) & = &  4 > 9 -8 \\
d((-7/3,-7/3-3/2,),\alpha^{(9)}) & = &  4 > 9 -8 \\
d((-5/2,-5/2,-4/3),\alpha^{(9)}) & = &  3 > 10 -8 \\
d((-5/2,-5/3,-5/3),\alpha^{(9)}) & = &  3 >  8-8
\end{eqnarray*}
\begin{eqnarray*}
d((1,2,2),\alpha^{(9)}) & = &  5 >  9-8 \\
d((1,2,3),\alpha^{(9)}) & = &  6 >  8-8 \\
d((1,3,3),\alpha^{(9)}) & = &  7 >  8-8 \\
d((2,2,2),\alpha^{(9)}) & = &  6 >  9-8.
\end{eqnarray*}
Our argument is now complete
\end{proof}

Having proved eventually that $9$ is well-behaved on tori, we are now ready to apply
Proposition~\ref{atoroidal:prop} to compute $c_9$ on the
negatively curved fillings of $N$. To do this we need to slightly extend the ``distance''
function $d$ introduced in Section~\ref{slope:triod:section}.  Recall
that $d(\alpha,\alpha')$ was originally defined when $\alpha$ is
either a $\theta$-graph or a slope on a given torus, and $\alpha'$ is a $\theta$-graph.
Now we also allow the case where both $\alpha$ and $\alpha'$ are
slopes, defining $d(\alpha,\alpha')$ to be 0 if they coincide and
$+\infty$ if they do not. The following result is an easy consequence of
Proposition~\ref{atoroidal:prop}

\begin{lemma} \label{atoroidal2:lem} 
  Let $M$ be hyperbolic. Then
  $c_9(M)$ equals the minimum of $c(B_i)+d(\alpha,\alpha^{(i)})$ over
  all realizations of $M$ as $N_\alpha$ and all $6\leqslant i \leqslant 10$.
\end{lemma}

Recall that at the end of Section~\ref{statement:section} we have
defined a certain function $h$ on triples of rational numbers.
The next result implies the second part of Theorem~\ref{hyp:Dehn:fill:teo}.

\begin{prop} 
Assume $M$ is irreducible, closed, atoroidal, and not a Seifert manifold.
 Then $c_9(M)$ equals the minimum of $h(p/q,r/s,t/u)$ over all
  realizations of $M$ as $N_{p/q,r/s,t/u}$.
\end{prop}

\begin{proof}
By Lemma~\ref{atoroidal2:lem}, it is
  sufficient to fix $p/q,r/s,t/u$ and prove that $h(p/q,r/s,t/u)$
  is $\min\{k_6,\ldots,k_{10}\}$ where
  $k_i=c(B_i)+d\big((p/q,r/s,t/u),\alpha^{(i)}\big)$.
  
  We first note that $d(p/q,\theta^{(-2)})=|p+2q,q|-1$, whence $k_9 =
  6+|p+2q,q|+|r+2s,s|+|t+2u,u|$. We will compute our minimum by
  systematically analyzing for $i\ne 9$ the triples on which $k_i$ is
  smaller than $k_9$.  We start with $k_{10}$ and note that it is
  smaller than $k_9$ precisely when $p/q<-2$ up to permutation, and
  $k_{10}=-1+k_9$ in this case.
  
  Concerning $k_8$, we note that it is finite precisely when $p/q=1$
  up to permutation.  Moreover $d(1,\theta^{(-2)})=1$, so $k_8$ is
  smaller than $k_9$ for $p/q=1$, and again the difference is $1$.
  
Now $k_7$ is finite precisely when $p/q=1$ and $r/s=-5$ up to permutation.
Moreover $d(1,\theta^{(-2)})=1$ and $d(-5,\theta^{(-2)})=2$,
so there is a $-3$ gain on the first two coordinates. And $d(\theta^{(-1)},\theta^{(-2)})=1$,
so there is a global gain whatever the last coordinate. More precisely, there is a further
$-1$ gain with respect to $k_9$ when $t/u>-1$, and a $+1$ loss if $t/u<-1$.
Summing up we have $k_7=-4+k_9$ if $p/q=1,r/s=-5,t/u>-1$ up to permutation, and
$k_7=-2+k_9$ if $p/q=1,r/s=-5,t/u<-1$ up to permutation.

We are left to analyze $k_6$, which is finite when $p/q=1$ and $r/s=-4$ up to
permutation. Since $d(1,\theta^{(-2)})=d(-4,\theta^{(-2)})=1$ we have a $-2$ gain on the
first two coordinates, together with a $-1$ gain coming from the fact that
$c(B_6)=8$. As above there is a further $-1$ gain if 
$t/u>-1$, and a $+1$ loss $t/u<-1$.
Therefore $k_7=-4+k_9$ if $p/q=1,r/s=-4,t/u>-1$ up to permutation, and
$k_7=-2+k_9$ if $p/q=1,r/s=-4,t/u<-1$ up to permutation.

This analysis easily implies that indeed $\min\{k_6,\ldots,k_{10}\}$ is 
$h(p/q,r/s,t/u)$.\end{proof}

\section{Complexity of Seifert manifolds and torus bundles}\label{Seifert:section}
The key tool used for the computation of $c_n$ in the atoroidal case,
namely Proposition~\ref{B_2:or:B_3:except:one:prop}, is typically not
valid for Seifert manifolds and torus bundles, so a different technique
has to be employed.

\paragraph{Seifert manifolds with compatible markings}
The bricks involved in the calculation of $c_3$ have markings containing fibres
of their product fibration. This motivates the following:

\begin{defn}
  \emph{The total space of a Seifert fibration has \emph{compatible}
		marked boundary if every $\theta$-graph of the marking
    contains a fibre (up to isotopy). }
\end{defn}

It turns out that the set of \emph{oriented} 
manifolds with markings compatible with a given Seifert fibration
is parameterized by $\matZ$. 
To show this, denote by $\theta_{(i)}$ the $\theta$-graph which contains
the slopes $0$, $1/i$, and $1/(i+1)$, and note that any 
oriented Seifert manifold with compatible marking can be expressed as
$$M = (F\timforsetil S^1)_{\frac{p_1}{q_1},\ldots,\frac{p_h}{q_h},
\theta_{(i_1)},\ldots,\theta_{(i_k)}}$$
with $p_j>q_j>0$ for all $j$.
Here $F$ is a surface with $h+k$ boundary components,
$F\timforsetil S^1$ denotes the unique oriented circle bundle over $F$,
and $k\geqslant 1$ otherwise $M$ is closed.
Moreover, slopes are expressed using homology bases as described 
at the beginning of Section~\ref{statement:section}, 
\emph{i.e.}~positive bases with first vector contained in a 
section of $F\timforsetil S^1$, and a fibre as second vector.
These requirements do not determine the bases uniquely if $k>1$, but
the general theory of Seifert bundles~\cite{FoMa} easily implies that
the $\matZ$-valued invariant
$$t(M)=i_1+\ldots+i_k$$
is well-defined and parameterizes compatible markings on $M$.

\begin{prop}\label{compatible:Seifert:prop}
Reversing orientation we get $t(-M) = -t(M)-h-k.$
\end{prop}
\begin{proof}
Under orientation reversal,
coordinates of slopes get mutiplied by $-1$, hence $\theta_{(i)}$ becomes
$\theta_{(-i-1)}$. Therefore
\begin{eqnarray*}
-M  & = & 
(F\timforsetil S^1)_{-\frac{p_1}{q_1},\ldots,-\frac{p_h}{q_h},
\theta_{(-i_1-1)},\ldots,\theta_{(-i_k-1)}} \\
 & = &
(F\timforsetil S^1)_{\frac{p_1}{p_1-q_1},\ldots, \frac{p_h}{p_h-q_h},\theta_{(-i_1-1-h)},
\ldots,\theta_{(-i_k-1)}}
\end{eqnarray*}
whence the conclusion.
\end{proof}

\paragraph{Small bricks have compatible marking}
The first bricks $B_0,\ldots,B_4$ have markings compatible
with the product fibration (but $B_5$ has not). 
We list their twisting numbers $t$ below.
For all bricks except $B_0$ there are two choices of $t$, depending on orientation.
\begin{equation} \label{seif:bricks:eqn}
\begin{tabular}{c}
$t(B_0) = -1, \quad t(B_1) \in \{-1,0\}, \quad t(B_2) \in \{-2,1\}$, \\
$t(B_3) \in \{-2,0\}, \quad t(B_4) \in \{-2,-1\}$.
\end{tabular}
\end{equation}
Let $S$ be the M\"obius strip. 
We will also need below the marked manifold 
$Z = S\timtil S^1$ with $t(Z)\in\{-1,0\}$.

Given two \emph{oriented} Seifert manifolds with compatible markings, 
we define now an assembling
of them to be \emph{compatible} if it is
orientation-reversing and fibre-preserving. The result of this assembling
is another oriented Seifert manifold with compatible marking. The same definition applies
to a self-assembling. In the following lemma, $N$ may be closed, and in this case $t(N)$
is the twisting number defined in Section~\ref{statement:section}.

\begin{lemma} \label{S_2:lemma}
If $N$ is obtained as a compatible assembling of $M$ and $M'$, then $t(N)=t(M)+t(M')+1$.
If $N$ is obtained as a compatible self-assembling of $M$, then $t(N)=t(M)+1$.
\end{lemma}

\begin{proof}
We prove the first assertion only, because a similar argument shows the second one.
Let us note that
$$t\left((F\timtil S^1)_{\frac{p_1}{q_1}, \ldots, \frac{p_h}{q_h}, \frac 1b, 
\theta_{(i_1)},\ldots,\theta_{(i_k)}}\right) = b+i_1+\ldots+i_k$$
whenever $p_j>q_j>0$, and this formula gives the correct value of the twisting
number even for $k=0$, \emph{i.e.} when the manifold is closed.
Hence we can write
$$M=(F\timtil S^1)_{\frac{p_1}{q_1}, \ldots, \frac{p_h}{q_h}, \frac 1{t(M)}, 
\theta_{(0)},\ldots,\theta_{(0)}}, \quad 
M'=(F'\timtil S^1)_{\frac{p_1'}{q_1'}, \ldots, \frac{p_h'}{q_h'}, \frac 1{t(M')+1}, 
\theta_{(-1)},\theta_{(0)}, \ldots,\theta_{(0)}}$$
and assume that one boundary torus of $M$ is glued to the boundary torus
of $M'$ marked by $\theta_{(-1)}$. The gluing matrix $A$ has determinant $-1$ and
maps $0$ to $0$ (because the assembling is compatible), so it maps
$\{1,\infty\}$ to $\{-1,\infty\}$ (because it maps
$\theta_{(0)}$ to $\theta_{(-1)}$). This implies that 
$A=\pm\tiny{\matr 100{-1}}$, \emph{i.e.}~$A$ matches the both the fibres and the
sections of the fibrations, whence
$$N=\big((F\cup F'\big)\timtil S^1)_{\frac{p_1}{q_1}, \ldots, \frac{p_h}{q_h}, 
\frac{p_1'}{q_1'}, \ldots, \frac{p_h'}{q_h'}, \frac 1{t(M)},\frac 1{t(M')+1}, 
\theta_{(0)},\ldots,\theta_{(0)}}$$
and the conclusion follows.
\end{proof}

We recall now that the Seifert parameters for a manifold with
non-empty boundary are almost always unique.  In the setting
of Seifert manifolds with compatible markings, the next lemma describes two 
non-uniqueness instances needed below. Recall that 
$D$ is the disc and $D_i$ is the disc with $i$ holes.

\begin{lemma} \label{two:fibrations:lemma}
The brick $B_2$ equals $(D_1\times S^1)_{2,\theta_{(-1)}}$, and
$Z$ equals $(D_2\times S^1)_{2,2,\theta_{(-1)}}$.
\end{lemma}

\begin{proof}
By~(\ref{seif:bricks:eqn}) we have $B_2=(D\times S^1)_{\theta_{(-2)}}$, but
$(D_1\times S^1)_2$ is homeomorphic to $D\times S^1$. Using the explicit
homeomorphism one easily gets the first assetion.

For the second assertion, recall that the fibrations
$(D_2\times S^1)_{2,2}$ and $S\timtil S^1$ are homeomorphic,
and the fibres have geometric intersection $1$ when the fibrations
are identified.
A $\theta$-graph is compatible with both fibrations when it contains
both the fibres, so there are two such $\theta$-graphs $\theta$ and $\theta'$.
Set 
$$\begin{array}{rclcrcl}
t_D & = & t((D_2\times S^1)_{2,2,\theta}), & &
t_D' & = & t((D_2\times S^1)_{2,2,{\theta'}}), \\
t_S & = & t((S\timtil S^1)_\theta), & & 
t_S' & = & t((S\timtil S^1)_{\theta'}).\end{array}
$$

Reversing orientation
$t_S$ gets replaced by $-t_S-1$, but the property that a $\theta$-graph contains
both fibres is insensitive to orientation, so $-t_S-1\in\{t_S,t_S'\}$
whence $t_S'=-t_S-1$. Moreover $\theta$ and $\theta'$ share two curves, 
so $t_S$ and $t_S'$ must
be consecutive, which implies that $\{t_S,t_S'\}=\{-1,0\}$.
A similar argument proves that $\{t_D,t_D'\}=\{-2,-1\}$. Therefore
$Z = (S\timtil S^1)_{\theta_{(t_S)}}$ for $t_S\in\{-1,0\}$ is also 
$(D_2\times S^1)_{2,2,\theta_{(t_D)}}$ for $t_D\in\{-2,-1\}$.
\end{proof}

The previous result implies the fact (also used in the proof) that
$Z$ has a marking compatible with both its original $S\timtil S^1$
and its alternative $(D_2\times S^1)_{2,2}$ fibrations.
Lemmas~\ref{S_2:lemma} and~\ref{two:fibrations:lemma} readily imply the following:

\begin{cor}\label{Z:ass:cor}
There is an assembling of one $B_4$ and two $B_2$'s realizing $Z$.
\end{cor}

We can now compute $c_3$ on Seifert manifolds.

\begin{teo}\label{c3:Seifert:teo}
Let $M=\big(F, (p_1, q_1), \ldots, (p_k, q_k), t\big)$ be a
closed genuine Seifert manifold
with parameters normalized so that
$k-\chi(F)>0$, $p_i>q_i>0$, $t\geqslant -k/2$.
If $M$ is not $\big(S^2,(2,1),(2,1),(j,1),-1\big)$ for $j\in\{2,3\}$ then:
$$c_3(M)=
\max\{0,t-1+\chi(F)\}+6(1-\chi(F))+\sum_{i=1}^k(|p_i,q_i|+2).$$
\end{teo}

\begin{proof}
Our proof is organized as a sequence of claims and constructions, and it uses the
Seifert $Z$ introduced above, fibred as $S\timtil S^1$.
The assumption that $M$ is not 
$\big(S^2,(2,1),(2,1),(j,1),-1\big)$ for $j\in\{2,3\}$
ensures that $M$ is not a closed brick of complexity up to $3$, so every
realization of $c_3(M)$ will consist of blocks $B_0,B_1,B_2,B_3,B_4$.

\emph{Claim 1. There exists an assembling of bricks realizing $c_3(M)$ and
such that:
\begin{itemize}
\item it involves only copies of $B_2$, $B_3$, and $B_4$;
\item each copy of $B_3$ is saturated (\emph{i.e.}~it is a union of fibres);
\item each copy of $B_4$ is either saturated or assembled to two copies of $B_2$
giving a saturated $Z$ as a result;
\item each copy of $B_2$ not contained in the $Z$'s just described is
a regular neighbourhood of an exceptional fibre of $M$, and conversely
each exceptional fibre is contained in one such $B_2$.
\end{itemize}}

To prove the claim we fix an assembling realizing $c_3(M)$ with a
minimal number of copies of $B_4$.  The only case where $B_0$'s cannot be dismissed
is when $M$ is a torus bundle, and this is not our case.
Lemma~\ref{no:B_1:lem} easily implies that our assembling cannot
involve $B_1$'s, so there are only $B_2$'s, $B_3$'s, and $B_4$'s.
Moreover, there is at least one $B_4$, otherwise $M$ is a lens space.

Let us concentrate on one $B_4$ and note that, by irreducibility of the
blocks involved, each component of $\partial B_4$ is either incompressible
in $M$ or the boundary of a solid torus $\Tsolid$.
Since $n=3$ is well-behaved on tori, minimality of the number of
$B_4$'s shows that such a $\Tsolid$ must be an assembling of
some $B_3$'s and one $B_2$. In particular, $\Tsolid$ 
lies outside $B_4$, and hence defines a Dehn filling of $B_4$.
Moreover, if the filled slope is $r/s$ with respect to the basis fixed on 
$H_1(\partial B_4)$, we have $r\geqslant 2$ (up to switching signs) by 
Proposition~\ref{c3:atoroidal:prop}-(\ref{c3:atoroidal:degenerate:point}).
Let $i\in\{0,1,2,3\}$ be the number of
incompressible components of $\partial B_4$.
Recall that incompressible tori are always fibred (because
$M$ is genuine, see~\cite{FoMa}).
The next discussion shows that we can modify the assembling so that
all $B_4$'s are as in the claim:

\begin{itemize}
\item If $i=3$ then $B_4$ is saturated and compatible (because its fibration is unique);
\item If $i=2$, the one solid torus incident to $B_4$
gives a $r/s$-Dehn filling $N$ of $B_4$ with $r\geqslant 2$, so
$N=(D_1,(r,s))$ has a unique Seifert fibration,
which restricts to the fibration of $B_4$. It follows that $B_4$ is saturated;
\item If $i=1$ the two $\Tsolid$'s incident to $B_4$ give
some saturated Dehn filling of $B_4$, again denoted by $N$. Let $r_1/s_1$ and
$r_2/s_2$ be the filled slopes,
with $r_1,r_2\geqslant 2$. If $r_1>2$ or $r_2>2$ then
$N = (D,(r_1,s_1),(r_2,s_2))$ has a unique fibration, so $B_4$ is saturated.
Assume $r_1=r_2=2$ and choose homology bases on $\partial(D_2\times S^1)$
so that $B_4$ is given by 
$\partial(D_2\times S^1)_{\theta_{(0)},\theta_{(0)},\theta_{(t)}}$ for some 
$t\in\{-2,-1\}$, and the filled components are those marked by 
$\theta_{(0)}$. Note that this change of bases affects $s_1$ and $s_2$:
the new value of $s_i$ is now intrinsically determined by the fact
that $|2,s_i|-1$ copies of $B_3$ (and one copy of $B_2$) realize the
filling $2/s_i$. Moreover
$$N = (D\times S^1)_{\frac 2{s_1}, \frac 2{s_2}, \theta_{(t)}} = 
(D\times S^1)_{2,2,\theta_{(t+(s_1-1)/2+(s_2-1)/2)}}.$$
Since $|2,2h+1|=|h|+1$, setting $s_i=2h_i+1$, we see that the original
assembling giving $N$ consists of one $B_4$, two $B_2$'s, and
$|h_1|+|h_2|$ copies of $B_3$. But $N$ can also be realized as one
$Z=(D\times S^1)_{2,2,\theta_{(t)}}$ with 
$|s_1/2+s_2/2-1|=|h_1+h_2|\leqslant |h_1|+|h_2|$ copies of $B_3$
(whence $h_1 h_2\geqslant 0$). Since
$(D\times S^1)_{2,2}$ admits two fibrations only, either the original
$B_4$ or the new $Z$ are saturated;
\item If $i=0$ then $M$ is atoroidal and $B_4$ is the centre of the assembling,
by Proposition~\ref{B_2:or:B_3:except:one:prop}.
\end{itemize}
Turning to the $B_2$'s not contained in the $Z$'s and to the
$B_3$'s, we note that they are saturated, because their boundary is.
The correspondence between exceptional fibres and $B_2$'s not contained
in the $Z$'s is a byproduct of the above construction.

\emph{A modified realization and its graph.}
A self-assembling which identifies two tori $T$ and $T'$
matches $\theta\subset T$ either with $\theta'\subset T'$ or with another $\theta$-graph
on $T'$ having distance $1$ from $\theta'$. In the second case
we can insert one copy of
$B_3$ and get a self-assembling which matches the $\theta$-graphs. Doing this
to an assembling giving $c_3(M)$ as in Claim 1, we get a
realization of $M$ as an assembling of blocks of type $B_2,B_3,B_4,Z$ such
that $\theta$-graphs are always matched. We call \emph{fake} the $B_3$'s we have inserted,
and we define a graph $\calG$
with the blocks of our assembling as vertices
and the gluings as edges. So the $B_4$'s are trivalent, the $B_3$'s are bivalent,
the $B_2$'s and $Z$'s are univalent.

\emph{Claim 2. Let $V_0,\ldots,V_{n+1}$ be the vertices of a simple (possibly
closed) path in $\calG$. Suppose $V_0,V_{n+1}$ are (possibly equal and)
of type $B_4$ or $Z$ and the $V_1,\ldots,V_n$ are of type $B_3$. Then $V_1,\ldots,V_n$
are compatible.}
Let $T_i$ be the torus along which $V_{i-1}$ is glued to $V_i$ and
let $\theta_i$ be the $\theta$-graph on $T_i$. Since the $T_i$'s are parallel
in $M$, we can identify them to a single fibred torus $T$. Now
$\theta_1,\ldots,\theta_{n+1}$ can be viewed as the vertices of a path
in $\Theta(T)$, because $d(\theta_i,\theta_{i+1})=1$.
Minimality shows that this path must be simple, so it is determined
by $\theta_1$ and $\theta_{n+1}$. Knowing
that $\theta_1$ and $\theta_{n+1}$ contain a fibre (because the blocks $B_4$ and $Z$
are compatible), it is now easy to show that
the other $\theta_i$'s do, and the claim is proved.

\emph{Claim 3. Let $V_0,\ldots,V_{n+1}$ be the vertices of a simple path in $\calG$.
Suppose $V_0$ is a $B_2$ block containing an exceptional fibre of type $(p,q)$.
Suppose $V_{n+1}$ is of type $B_4$ or $Z$
and $V_1,\ldots,V_n$ are of type $B_3$. Then $V_i$ is compatible for
$i=|p,q|,\ldots,n$ and it is not for $i=1,\ldots,|p,q|-1$.}
With notation as for the previous claim, we note that there is a preferred positive
basis $(a,b)$ for $H_1(T)$, such that $b$ is a fibre and $p\cdot a+q\cdot b$ is
a meridian of $V_0$. We use $(a,b)$ to identify $\calS(T)$ to
$\partial_\matQ(\matH^2)$ and note that $\theta_{n+1}$ contains the slope $0$, so
it is one of the vertices marked by solid dots in Fig.~\ref{compatible:fig}.
\begin{figure}
\begin{center}
\mettifig{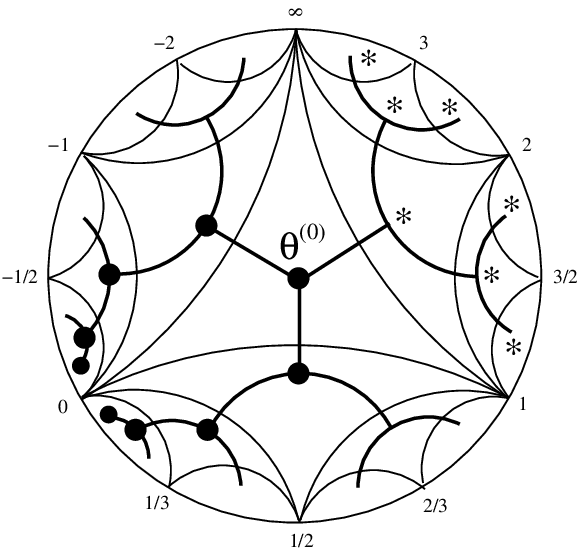}
\nota{The $\theta$-graph $\theta_{1}$ is one of the stars, and
$\theta_{n+1}$ is one of the solid dots.}
\label{compatible:fig}
\end{center}
\end{figure}
In addition $\theta_1$ contains $p/q$ and $1<p/q<\infty$, so
it is one of the vertices marked by stars. Using
Proposition~\ref{lines:and:|p,q|:prop} and recalling that $\theta^{(0)}$ is the
centre of the disc, we deduce that
$$n=d(\theta_1,\theta_{n+1})=d(\theta_1,\theta^{(0)})+d(\theta^{(0)},\theta_{n+1})=
|p,q|+d(\theta^{(0)},\theta_{n+1}).$$
The figure also shows that precisely the last
$1+d(\theta^{(0)},\theta_{n+1})$ $\theta$-graphs among $\theta_1,\ldots,\theta_{n+1}$ contain $0$,
whence the conclusion at once.

\emph{Claim 4. If $b_3^c$ is the number of non-fake compatible $B_3$'s in our assembling
then}
$$
c_3(M)=b_3^c+6(1-\chi(F))+\sum_{i=1}^k(|p_i,q_i|+2).
$$
Let us denote by $b_3,b_4,z,s$ the number of
non-fake $B_3$'s, of $B_4$'s, of $Z$'s, and of self-assemblings.
By assumption we have $c_3(M)=b_3+3b_4+3z+ 6s$,
and Claims~2 and~3 show that
$$b_3=b_3^c+\sum_{i=1}^k(|p_i,q_i|-1)=b_3^c-3k+\sum_{i=1}^k(|p_i,q_i|+2).$$
We now claim that $b_4=k-\chi(F)$ and $z=2-2s-\chi(F)$. These relations easily
imply the desired conclusion. Relation $\chi(F)=-b_4+k$ readily follows from
Lemma~\ref{S_2:lemma} and
the identities~(\ref{seif:bricks:eqn}). To show the other
relation consider the fibred space arising from the assembling of
the blocks of type $B_2,B_3,B_4$ only, without self-assemblings.
The base surface is then a sphere with some holes.
Filling $z$ of these holes and identifying the boundary circles of
$s$ pairs of them we get the closed surface $F$. So there are $z+2s$ holes
in all, and $\chi(F)=2-(z+2s)$, whence the desired relation.

\emph{Claim 5. Conclusion: $b_3^c=\max\{0,t-1+\chi(F)\}$.}
We first note that Claim 3 (and its proof) imply that by dismissing from our
assembling the $B_2$'s and the non-compatible $B_3$'s we get an $M'$ with $t(M')=t(M)$.
We also note that 
$z+b_4>0$, otherwise $M$ would be a lens space. 
Then we slightly change our perspective and ask ourselves
to what extent do the numbers $b_3^c,b_4,z,k$
actually determine $M'$.
Namely we discuss what $N$'s can be realized as follows:
\begin{itemize}
\item $N$ is connected and has $k$ boundary components;
\item $N$ is obtained via
fibre-preserving gluings from
$b_4$ copies of $B_4$, $z$ copies of $Z$,
$b_3^c$ copies of $B_3$, and some number $f\geqslant 0$
of additional copies of $B_3$, with $z+b_4>0$;
\item the graph of this gluing
remains connected after removing the $f$ additional copies of $B_3$.
\end{itemize}
These conditions easily imply that
there are $n=(z+2(b_3^c+f)+3b_4-k)/2$
gluings, precisely $s=(3b_4-z-k)/2-(b_4-1)$ of which are self-assemblings.
(Of course our notation for $s$ is consistent with above.)
Moreover $N$ has base surface 
with Euler characteristic $-b_4$ and no exceptional fibres. For the other Seifert
parameters we have some flexibility, because the gluing maps are not uniquely determined.
The flexibility reduces if we fix orientations on all the blocks, because
then the twisting number of each block is determined,
so $t(N)$ is computed via
Lemma~\ref{S_2:lemma}
If $s=0$ then the base surface is orientable, otherwise
we can realize both an orientable and a non-orientable one without
changing the other parameters.

Let us now discuss exactly what possible $t(N)$'s can be realized
once $b_3^c,b_4,z,k$ are given. Recall that, depending on orientation,
$t(B_4)$ is $-2$ or $-1$. Similarly $t(Z)$ is $-1$ or $0$ and $t(B_3)$ is
$-2$ or $0$. So:
$$-2b_4-z-2(b_3^c+f)+n\leqslant t(N) \leqslant -b_4+n$$
and all possibilities can be realized by appropriate choices of
orientations (because $z+b_4>0$). Noting that $f$ can take any value between
$0$ and $s$ and plugging in the values of $n$ and $s$
these inequalities easily imply that
\begin{equation}\label{chi:limits:eqn}
-b_4-b_3^c-1 \leqslant t(N) \leqslant b_4+b_3^c-k+1
\end{equation}
and all these values can be realized.

We are eventually ready to prove the desired equality, so we
assume that the parameters $b_3^c,b_4,z,k$ and $f$
arise as above from an assembling realizing $c_3(M)$.
Since $b_4-k=-\chi(F)$, the right-hand side of
equation~(\ref{chi:limits:eqn}) readily implies that
$b_3^c\geqslant\max\{0,t-1+\chi(F)\}$.
We prove equality by contradiction,
so we assume that $b_3^c>0$ and $e<b_4+b_3^c-k+1$.
Since $t\geqslant -k/2$ by assumption, we also have $t>-b_4-b_3^c-1$.
Now equation~(\ref{chi:limits:eqn})
and the discussion leading to it prove that
$(M',X)$ could be realized with $b_3^c-1$ copies of $B_3$
and all other parameters unchanged: a contradiction to minimality.
\end{proof}

\begin{prop}\label{no:Bi:for:Seifert:prop}
Let $M$ be a closed Seifert manifold. Then:
\begin{itemize}
\item the bricks $B_6,\ldots,B_{10}$ give no contribution to $c_8(M)$ and $c_9(M)$;
\item if $M$ is toroidal, the brick $B_5$ gives no contribution to $c_8(M)$.
\end{itemize}
\end{prop}

\begin{proof}
The first point is easy. Assume there is a hyperbolic brick $B\in\{B_6,\ldots,B_{10}\}$
in an assembling that realizes $c_8(M)$ or $c_9(M)$. Each boundary component of $B$
either bounds a solid torus outside $B$ or is incompressible in $M$. Let $M'$ be
the filling of $B$ given by the union of $B$ and the solid tori incident to it.
Then, by Propositions~\ref{non:neg-curved:degenerate:prop} 
and~\ref{non:neg-curved:degenerate:prop2}, 
either the filling is degenerate or
$M'$ is hyperbolic, so it is one of the blocks of the JSJ decomposition of $M$,
which is absurd.

Turning to the second point, we have
$$B_4=(D_2\times S^1)_{\theta^{(0)},\theta^{(0)},\theta^{(-1)}}\quad
{\rm and}\quad B_5=(D_2\times S^1)_{2,3,\theta}\quad {\rm with}\ 
\theta\supset\{-5/4,-6/5,-1\}.$$
Let us then assume by contradiction that there is a copy of 
$B_5$ in an assembling realizing $c_8(M)$.
Since $M$ is toroidal, the single boundary component of $B_5$ is incompressible.
Therefore $B_5$ is a union of fibres, and its induced Seifert structure, being
unique, is precisely that of $(D_2\times S^1)_{2,3,\theta}$.
Moreover, glued to $\partial B_5$, we 
have some $B_3$'s (possibly none of them) and then a brick $B$ equal to either
$B_4$ or $B_5$.

Consider first the case where $B=B_4$. Precisely as in the proof of 
Theorem~\ref{c3:Seifert:teo} we can see that either this $B_4$ is compatible
or its two free boundary components are glued to two $B_2$'s so 
to give a compatible $Z$. In both cases we see that the original $\theta$ on
$\partial B_5$ is modified along a simple path in $\calT$ until a compatible $\theta$-graph
$\theta'$ is reached, and the $d(\theta,\theta')$ edges in this
path correspond to the copies of $B_3$.
The sub-assembling consisting of $B_5$ and the neighbouring $B_3$'s
gives $(D_2\times S^1)_{2,3,\theta'}$
and contributes to $c_8(M)$ as 
$$c(B_5)+d(\theta,\theta')=8+d(\theta,\theta^{(-1)})
+d(\theta^{(-1)},\theta')=8+5+d(\theta^{(-1)},\theta')=13+d(\theta^{(-1)},\theta')$$
(identity $d(\theta,\theta')=d(\theta,\theta^{(-1)})
+d(\theta^{(-1)},\theta')$ is deduced from
Fig.~\ref{c8:fig}, because $\theta'$ contains $0$). 
Now the same $(D_2\times S^1)_{2,3,\theta'}$
can also be realized using $B_4$ with complexity
$$c(B_4)+d(\theta^{(0)},2)+d(\theta^{(0)},3)+
d(\theta^{(-1)},\theta')=3+0+1+d(\theta^{(-1)},\theta')=4+d(\theta^{(-1)},\theta').$$
This proves that $B_5$ can be dismissed with a gain of $9$ on complexity,
so the original assembling was not a realization of $c_8(M)$, a contradiction.

The case where $B=B_5$ is similar. Here $M$ is given by
$$\big(D_2\times S^1)_{2,3}\bigcup\nolimits_A \big(D_2\times S^1)_{2,3}.$$
But $M$ is Seifert,  so $A={\tiny\matr {\pm 1}0n1}$.
Since the two manifolds
glued along $A$ both have $\theta$ on the boundary, the number of 
$B_3$'s required to assemble them is 
$d(\theta,A\theta)$, so
$c_8(M)=2\cdot c(B_5)+d(\theta,A\theta)=16+d(\theta,A\theta)$. 
As above we can use $B_4$ instead of $B_5$
on both sides, deducing that 
$$c_3(M)\leqslant 2\cdot (3+0+1)+d(\theta^{(-1)},A\theta^{(-1)})
=8+d(\theta^{(-1)},A\theta^{(-1)}).$$

The conclusion now follows as soon as we show that 
$d(\theta^{(-1)},A\theta^{(-1)})<8+d(\theta,A\theta)$
for $A={\tiny\matr {\varepsilon}0n1}$, $\varepsilon=\pm 1$, $n\in\matZ$.
We begin with $\varepsilon=+1$ and note that $A=B^n$, where 
$B={\tiny\matr 1011}$. We consider the alternative picture in the
half-space model of $\matH^2$ of the Farey tessellation and its dual graph, 
as shown in Fig.~\ref{no:b5:fig}.
\begin{figure}
\begin{center}
\mettifig{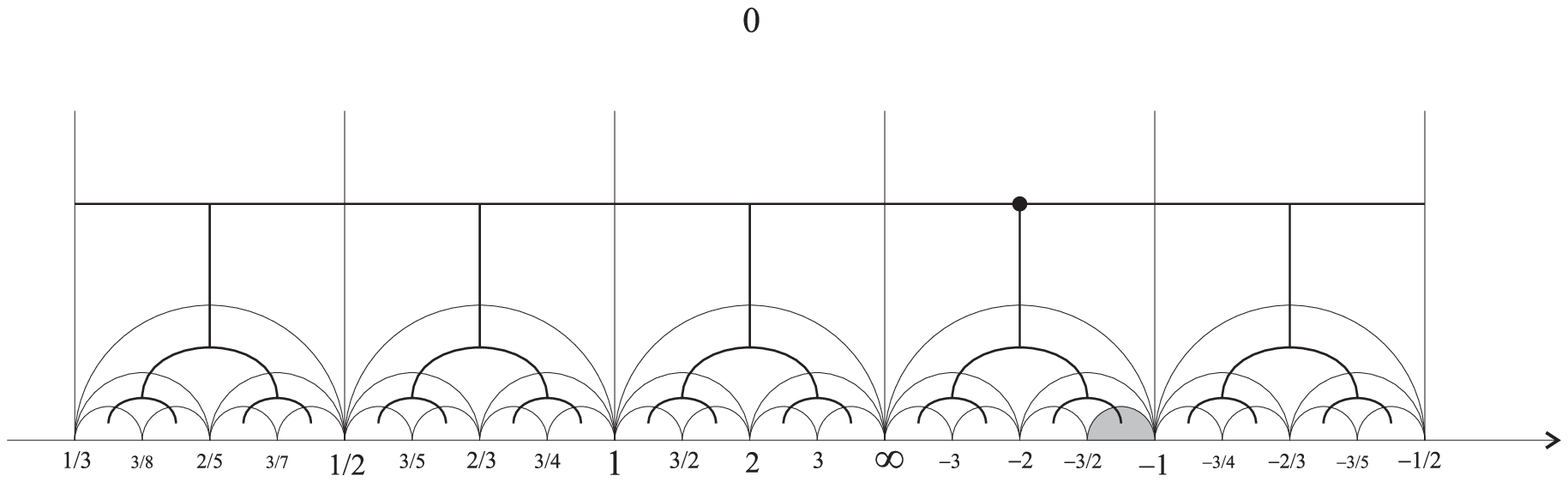, width=14cm}
\nota{An alternative picture of the Farey tessellation.} \label{no:b5:fig}
\end{center}
\end{figure}
(As above, the picture is only combinatorially, not geometrically, correct. Note
also that we have used a smaller and smaller character for $p/q$ as $|p,q|$ grows.)
Here $B$ acts on the tessellation of Fig.~\ref{no:b5:fig}
as the horizontal translation that moves each vertical half-line to the previous one.
Since $\theta^{(-1)}$ is the dot marked in the picture whereas $\theta$ lies somewhere
in the small shadowed region, we deduce that
$$d(\theta^{(-1)},A\theta^{(-1)})=|n|, \qquad d(\theta,A\theta)=
2d(\theta,\theta^{(-1)})+|n|=10+|n|$$
whence the desired inequality. When $A={\tiny\matr {-1}0n1}$ the action of $A$
on the tessellation is a horizontal translation followed by the reflection in the
vertical half-line through $0$ and $\infty$. It is then easy to see that 
$d(\theta^{(-1)},A\theta^{(-1)})=|n-1|$ and $d(\theta,A\theta)=10+|n-1|$, 
whence the conclusion.\end{proof}

The two previous results conclude the proof of Theorem~\ref{Seifert:bundle:teo}
for genuine Seifert manifolds, because
$c_n(M)$ was already computed in 
Propositions~\ref{c1:prop},
\ref{c3:atoroidal:prop}, and~\ref{c8:atoroidal:prop}
for $n\leqslant 8$ when $M$ is atoroidal and Seifert but not a member of
$\calM^*$.  For torus bundles, having already computed $c_1$ in
Proposition~\ref{c1:torus:bundle:prop}, we are left to show the following:

\begin{prop} \label{cn:torus:bundle:prop} 
$c_n(T_A) = c_1(T_A)$ for $n=2,\ldots,9$.
\end{prop}

\begin{proof} 
We begin by noting
the following topological fact:
if $H$ is an incompressible torus embedded in $T_A$ then, cutting
$T_A$ along $H$, we get either $T\times I$ or two copies of
$S\widetilde\times S^1$ (the twisted $S^1$-bundle over the M\"obius strip). And,
in the latter case, we have $A=\tiny{\matr{-1}0t{-1}}$ up to conjugation,
and $T_A=(K,+t)=(K,-t)$, where $K$ denotes the Klein bottle.

Let us then consider an assembling realizing $c_n(T_A)$, and let $\calH\subset T_A$ 
be the family of tori along which the gluings of the assembling are performed.
If all tori in $\calH$ are compressible, well-behavedness of $n\leqslant 9$ easily
implies that the assembling has a certain centre $B_i$ with
$i\geqslant 4$. However Propositions~\ref{c3:atoroidal:prop}, \ref{c8:atoroidal:prop}, 
\ref{non:neg-curved:degenerate:prop}, and~\ref{non:neg-curved:degenerate:prop2}
imply that a filling of $B_i$ giving $T_A$ is always degenerate, except
for $i=4$ and $A={\tiny\matr 01{-1}{-1}}$, but in this case we have 
$c_1(T_A)=c_3(T_A)=c(T_A)=6$, so $c_n(T_A)=6$ for $1\leqslant n\leqslant 9$.

Assume one $H\in\calH$ cuts $T_A$ into $T\times I$. Then a self-assembling is
performed along $H$ of a suitably marked $T\times I$.
Well-behavedness of $n\leqslant 9$ implies that we can realize this pair as
either $B_0$ or an assembling of $B_3$'s, so we have $c_n(T_A)=c_1(T_A)$.

Assume now one $H\in\calH$ cuts $T_A$ into two copies of 
$S\widetilde\times S^1=\big(D,(2,1),(2,1)\big)$.
Therefore $T_A = (K, \pm t)$ is Seifert, and every incompressible torus in $\calH$ is 
saturated.
Although $T_A$ is not genuine, this condition is enough in the proof of 
Theorem~\ref{c3:Seifert:teo} to conclude that $c_3(M) = 6 + \max\{0, |t|-1\} = 
\max\{6,||A||+5\} = c_1(M)$, as required.
\end{proof}

We conclude the section noting that a proof of Theorem~\ref{Mstar:teo}-(3)
is already contained in our results. In fact, Proposition~\ref{c3:atoroidal:prop}-(4)
proves that $c_3(C_{i,j})=1+i+j$ if $c^*(C_{i,j})>9$, whereas
$c^*(C_{i,j})=i+j$, Proposition~\ref{c8:atoroidal:prop}-(3) proves that
$c_8(E_k)=k+6$ if $c^*(E_k)>9$, whereas $c^*(E_k)=k+5$,
and Proposition~\ref{no:Bi:for:Seifert:prop} proves that
$c_9(C_{i,j})=c_3(C_{i,j})$ and $c_9(E_k)=c_8(E_k)$ under the same constraints.

\section{Tables of small manifolds}\label{tables:section}

We provide in this section a description, for $n=0,\ldots,9$,
of the set of closed orientable irreducible manifolds $M$ with
$c(M)=9$. A summary of our results was already given in 
Table~\ref{geom_compl:table}.
For $n\leqslant 6$ we give here complete lists, including
geometric information, while 
for $n=7,8,9$ we only concentrate on the most interesting geometric
features, addressing the reader to~\cite{website} for the lists.
The lists for $n\leqslant 6$ coincide with those of Matveev~\cite{Mat:book}.

Before turning to the data, we quickly sketch how we have 
found them.  Note first that for $M$ as above
one has $c(M)\leqslant 9$ if and only if $c_9(M)\leqslant 9$, and 
$c(M)=c_9(M)$ in this case. Different enumeration methods were used to
list manifolds of hyperbolic, geometric but non-hyperbolic, and graph
type. The hyperbolic case was already discussed in
Section~\ref{statement:section}.
When $M$ is geometric but non-hyperbolic, \emph{i.e.}~a Seifert manifold 
or a torus bundle,
Theorems~\ref{Seifert:bundle:teo} and~\ref{Mstar:teo} allow
to determine $c_9(M)$ directly. We have then written a computer program,
available from~\cite{website}, which lists
all such $M$'s with $c_9(M)\leqslant 9$, using also
Proposition~\ref{Seifert:generalities:prop} to avoid repetions.
To determine which graph manifolds $M$ satisfy 
$c_9(M)\leqslant 9$ we have used some estimates in the spirit of
Lemma~\ref{c3:graph:lem}, and the same computer program~\cite{website},
which does a careful analysis of the different
ways $M$ can be presented as a graph manifold.  Since this
paper is devoted to geometric manifolds, we do not insist on these
technicalities here.

In the rest of this section by `manifold' we mean a closed, connected, orientable,
irreducible 3-manifold. 
We recall now that we describe Seifert manifolds using the
filling parameters, not the orbital parameters (see the beginning of 
Section~\ref{statement:section}), and we proceed to a more accurate
explanation of our results. Below we will also use
the notation employed in Section~\ref{hyperbolic:complexity:section}
to describe graph manifolds (and an obvious extension of this notation).

\newcommand\tinymatr[4]{\tiny{\matr {#1} {#2} {#3} {#4}}}

\paragraph{Complexity up to 5}
It is well-known that the only manifolds of complexity $0$
are $S^3$, $L_{3,1}$, and $\matRP^3$. In complexity $1$ we have only the lens spaces 
$L_{4,1}$ and $L_{5,2}$. In complexity  $2$ the lens spaces 
$L_{5, 1}$, $L_{7, 2}$, and $L_{8, 3}$, and the Seifert space
$\big(S^2,(2,1),(2,1),(2,1),-1\big)$, which carries elliptic geometry
and belongs to $\calM^*$.
A similar phenomenon occurs in complexity $3$, with the lens spaces
$$L_{6, 1}\ \ L_{9, 2}\ \ L_{10, 3}\ \ L_{11, 3}\ \ L_{12, 5}\ \ L_{13, 5}$$
and the elliptic Seifert $\big(S^2,(2,1),(2,1),(3,1),-1\big)$
which belongs to $\calM^*$.
In complexity $4$ there are the lens spaces
$$\begin{array}{c}
L_{7, 1}\ \ L_{11, 2}\ \ L_{13, 3}\ \ L_{14, 3}\ \ L_{15, 4}\ \ 
L_{16, 7}\ \ L_{17, 5}\ \ L_{18, 5}\ \ L_{19, 7}\ \ L_{21, 8}
\end{array}$$
and the elliptic Seifert spaces
$$\begin{array}{cc}
\big(S^2,(2,1),(2,1),(2,1)\big), &
\big(S^2,(2,1),(2,1),(3,2),-1\big), \\
\big(S^2,(2,1),(3,1),(3,1),-1\big), &
\big(S^2,(2,1),(2,1),(4,1),-1\big),
\end{array}$$
the last two of which belong to $\calM^*$.

In complexity $5$ we have the lens spaces
$$\begin{array}{c}
L_{8, 1}\ \ L_{13, 2}\ \ L_{16, 3}\ \ L_{17, 3}\ \ L_{17, 4}\ \ 
L_{19, 4}\ \ L_{20, 9}\ \ L_{22, 5}\ \ L_{23, 5}\ \ L_{23, 7}\\ 
L_{24, 7}\ \ L_{25, 7}\ \ L_{25, 9}\ \ L_{26, 7}\ \ L_{27, 8}\ \ 
L_{29, 8}\ \ L_{29, 12}\ \ L_{30, 11}\ \ L_{31, 12}\ \ L_{34, 13}
\end{array}$$
and the following elliptic Seifert spaces
$$\begin{array}{c}
\big(S^2,(2, 1),(2, 1),(2, 1),1\big),\ \quad
\big(S^2,(2, 1),(2, 1),(3, 1)\big),\\
\big(S^2,(2, 1),(2, 1),(3, 2)\big),\ \quad 
\big(S^2,(2, 1),(2, 1),(4, 3),-1\big),\\
\big(S^2,(2, 1),(2, 1),(5, 2),-1\big),\ \quad 
\big(S^2,(2, 1),(2, 1),(5, 3),-1\big),\\
\big(S^2,(2, 1),(3, 1),(3, 2),-1\big),\ \quad
\big(S^2,(2, 1),(3, 2),(3, 2),-1\big),\\
\big(S^2,(2, 1),(3, 1),(4, 1),-1\big),\ \quad 
\big(S^2,(2, 1),(2, 1),(5, 1),-1\big),\\
\big(S^2,(2, 1),(3, 1),(5, 1),-1\big),
\end{array}$$
the last three of which belong to $\calM^*$ (the last one is
Poincar\'e's homology sphere).
Note that all the manifolds found so far are Seifert bundles 
over $S^2$ with at most three exceptional fibres.

\paragraph{Complexity 6}
The list of lens spaces in complexity $6$ consists of
$$\begin{array}{c}
L_{9, 1}\ \ L_{15, 2}\ \ L_{19, 3}\ \ L_{20, 3}\ \ L_{21, 4}\ \ L_{23, 4}\ \ 
L_{24, 5}\ \ L_{24, 11}\ \ L_{27, 5}\ \ L_{28, 5}\ \ L_{29, 9}\ \ L_{30, 7}\ \ 
L_{31, 7}\\ L_{31, 11}\ \ L_{32, 7}\ \  L_{33, 7}\ \ L_{33, 10}\ \ L_{34, 9}\ \ 
L_{35, 8}\ \ L_{36, 11}\ \  L_{37, 8}\ \ L_{37, 10}\ \ L_{39, 14}\ \ L_{39, 16}\ \ 
L_{40, 11}\\  L_{41, 11}\ \ L_{41, 12}\ \ L_{41, 16}\ \ L_{43, 12}\ \ L_{44, 13}\ \  
L_{45, 19}\ \ L_{46, 17}\ \ L_{47, 13}\ \ L_{49, 18}\ \ L_{50, 19}\ \ L_{55, 21}.
\end{array}$$
We then have many Seifert spaces fibred over $S^2$, with invariants and geometry
as described in table~\ref{s2_c6:table}.
We have grouped together the three members of $\calM^*$equals $7$
and the two Seifert spaces with $4$ exceptional fibres.

Besides those based on $S^2$, a few more Seifert manifolds have complexity $6$, as
described in  Table~\ref{seif_c6:table}, where
$T$ denotes the torus and $K$ the Klein bottle.
One should actually note that there is another such manifold, namely
$\big(K,0\big)$, which was already listed in Table~\ref{s2_c6:table}
as $\big(S^2,(2, 1),(2, 1),(2, 1),(2, 1),-2\big)$. Note also that a torus
bundle $T_A$ has $c_9(T_A)\leqslant 6$ if and only if $||A||\leqslant
1$, and in this case $T_A$ is one of the (flat or Nil) Seifert manifolds already
described.

\begin{table} \begin{center} \begin{tabular}{cc|cc} 
 invariants & geometry & invariants & geometry$\!\!$\\ \hline
$(2, 1),(2, 1),(2, 1),2$ & elliptic & 
$(2, 1),(2, 1),(3, 1),1$ & elliptic$\!\!$\\ 
$(2, 1),(2, 1),(3, 2),1$ & elliptic &
$(2, 1),(2, 1),(4, 1)$ & elliptic$\!\!$\\ 
$(2, 1),(2, 1),(4, 3)$ & elliptic &
$(2, 1),(2, 1),(5, 2)$ & elliptic$\!\!$\\ 
$(2, 1),(2, 1),(5, 3)$ & elliptic &
$(2, 1),(2, 1),(5, 4),-1$ & elliptic$\!\!$\\ 
$(2, 1),(2, 1),(7, 2),-1$ & elliptic &
$(2, 1),(2, 1),(7, 3),-1$ & elliptic$\!\!$\\ 
$(2, 1),(2, 1),(7, 4),-1$ & elliptic &
$(2, 1),(2, 1),(7, 5),-1$ & elliptic$\!\!$\\ 
$(2, 1),(2, 1),(8, 3),-1$ & elliptic &
$(2, 1),(2, 1),(8, 5),-1$ & elliptic$\!\!$\\ 
$(2, 1),(3, 1),(3, 1)$ & elliptic &
$(2, 1),(3, 1),(3, 2)$ & elliptic$\!\!$\\ 
$(2, 1),(3, 1),(4, 3),-1$ & elliptic &
$(2, 1),(3, 1),(5, 2),-1$ & elliptic$\!\!$\\ 
$(2, 1),(3, 1),(5, 3),-1$ & elliptic &
$(2, 1),(3, 2),(3, 2)$ & elliptic$\!\!$\\ 
$(2, 1),(3, 2),(4, 1),-1$ & elliptic &
$(2, 1),(3, 2),(4, 3),-1$ & elliptic$\!\!$\\ 
$(2, 1),(3, 2),(5, 2),-1$ & elliptic &
$(2, 1),(3, 2),(5, 3),-1$ & elliptic$\!\!$\\ 
$(3, 1),(3, 1),(3, 1),-1$ & flat &
$(3, 1),(3, 1),(3, 2),-1$ & Nil\\ 
$(3, 1),(3, 2),(3, 2),-1$ & Nil &
$(3, 2),(3, 2),(3, 2),-1$ & Nil\\  \hline
$(2, 1),(4, 1),(4, 1),-1$ & flat &
$(2, 1),(2, 1),(6, 1),-1$ & elliptic$\!\!$\\ 
$(2, 1),(3, 1),(6, 1),-1$ & flat & & \\ \hline
$\!\!(2, 1),(2, 1),(2, 1),(2, 1),-2\!$ & flat &
$\!(2, 1),(2, 1),(2, 1),(2, 1),-1\!$ & Nil
\end{tabular}\end{center}
\nota{Seifert manifolds of complexity $6$ based on $S^2$}\label{s2_c6:table}
\end{table}

\begin{table}[h] \begin{center} \begin{tabular}{c|c}
manifold & geometry \\ \hline 
$\big(\mathbb{P}^2,(2, 1),(2, 1),-1\big)$ & flat\\ 
$\big(\mathbb{P}^2,(2, 1) (2, 1)\big)$  & Nil\\ 
$\big(K,1\big)$ & Nil\\ 
$\big(T,0\big)$  & flat\\
$\big(T,1\big)$  & Nil
\end{tabular} \end{center}
\nota{Seifert manifolds of complexity $6$ not based on $S^2$}\label{seif_c6:table}
\end{table}

\paragraph{Complexity 7}
As mentioned above, we do not provide complete data for complexity higher than $6$,
addressing the reader to~\cite{website}. 
The numbers of manifolds per geometry are already in Table~\ref{geom_compl:table},
we only describe here the most interesting ones. In complexity $7$, besides
the lens spaces, there are 84 Seifert manifolds fibred over $S^2$ with 3 exceptional
fibres, with geometry either elliptic, or ${\rm SL}_2$, or Nil.
Three of these manifolds belong to $\calM^*$.
The other Seifert manifolds are as listed in 
Table~\ref{seif_c7:table}.
\begin{table} \begin{center} \begin{tabular}{c|c}
manifold & geometry \\ \hline 
$\big(S^2,(2, 1),(2, 1),(2, 1),(2, 1)\big)$ & Nil\\ 
$\big(S^2,(2, 1),(2, 1),(2, 1),(3, 1),-2\big)$ & ${\rm SL}_2$\\ 
$\big(S^2,(2, 1),(2, 1),(2, 1),(3, 1),-1\big)$ & ${\rm SL}_2$\\ 
$\big(S^2,(2, 1),(2, 1),(2, 1) (3, 2),-1\big)$  & ${\rm SL}_2$\\ 
$\big(\mathbb{P}^2,(2, 1),(2, 1),1\big)$ & Nil\\ 
$\big(\mathbb{P}^2,(2, 1),(3, 1),-1\big)$ & ${\rm SL}_2$\\ 
$\big(\mathbb{P}^2,(2, 1),(3, 1)\big)$ & ${\rm SL}_2$\\ 
$\big(\mathbb{P}^2,(2, 1) (3, 2)\big)$  & ${\rm SL}_2$\\ 
$\big(K,2\big)$ & Nil\\ 
$\big(T,2\big)$ & Nil
\end{tabular} \end{center}
\nota{Seifert manifolds of complexity $7$ except those over $S^2$
with $2$ or $3$ exceptional fibres}\label{seif_c7:table}
\end{table} 
In addition, we have the two Sol manifolds $T_A$ 
for $A={\tinymatr 3{-1}10},{\tinymatr {-3}1{-1}0}$, 
and 
$$\big(D,(2,1),(2,1)\big) \bigcup\nolimits_A \big(D,(2,1),(2,1)\big)$$
for $A=\tinymatr{1}{2}{-1}{-1},\tinymatr{1}{-1}{-1}{0},\tinymatr{0}{1}{-1}{0}$.
These manifolds
are fibred over the $1$-orbifold given by the interval with silvered endpoints,
with fibre a Klein bottle at the endpoints and a torus at the other points.
To conclude, we have the $4$ non-geometric manifolds 
$$\big(D,(2,1),(2,1)\big) \bigcup\nolimits_A \big(D,(2,1),(p,q)\big)$$
with $(p,q)\in\{(2,1),(3,2)\}$ and $A\in\{\tinymatr{0}{1}{-1}{-1},\tinymatr{0}{1}{1}{0}\}$.

\paragraph{Complexity 8}
After the 136 lens spaces, there are 226 Seifert
bundles over $S^2$ with three exceptional fibres and geometry
either elliptic, or ${\rm SL}_2$, or Nil. Of these 226 manifolds,
4 belong to $\calM^*$ and one is only discovered to have complexity $8$
using the brick $B_5$.
Next, there are 14 Seifert bundles over 
$S^2$ with 4 exceptional fibres and geometry either Nil or
$\mathbb{H}^2 \times S^1$, each in one case, or ${\rm SL}_2$, in the
other 12 cases. Then we have $14$ Seifert bundles
over $\matP^2$ with 2 exceptional fibres and exactly the same 
distribution of geometries.
The list of Seifert manifolds in complexity $8$ is completed
by $\big(K,3\big)$ and $\big(T,3\big)$, both carrying Nil geometry.

Next, we have the Sol manifolds $T_A$ with $A={\tinymatr 4{-1}10},{\tinymatr {-4}1{-1}0}$ and
$$\big(D,(2,1),(2,1)\big) \bigcup\nolimits_A \big(D,(2,1),(2,1)\big)$$
for $A=\tinymatr{1}{3}{0}{-1},\tinymatr{2}{3}{-1}{-1},
\tinymatr{2}{-1}{-1}{0},\tinymatr{1}{2}{1}{1},
\tinymatr{1}{2}{0}{1},\tinymatr{1}{-1}{1}{0},\tinymatr{1}{-1}{0}{-1}$. 

The list of manifolds of complexity $8$ is completed by $35$ non-geometric ones,
all of the form
$$\big(D,(2,1),(p,q)\big) \bigcup\nolimits_A \big(D,(r,s),(t,u)\big).$$

\paragraph{Complexity 9}
We have here, besides the $272$ lens spaces, Seifert bundles of 
many sorts. First, we have 586 over $S^2$ with 3 fibres, and geometry
${\rm SL}_2$, elliptic and (sporadically) Nil. Four of these
586 belong to $\calM^*$ and 3 require $B_6$ in the assembling
realizing the complexity.
Next, we have $41$ Seifert bundles over $S^2$ with 
4 fibres, always with ${\rm SL}_2$ geometry except for one Nil case,
and two ${\rm SL}_2$ bundles over $S^2$ with 5 fibres. Turning to
Seifert bundles over $\matP^2$, we have 34 with 2 fibres (all 
${\rm SL}_2$ except one Nil), and 2 with 3 fibres and geometry
${\rm SL}_2$. The list of Seifert bundles of complexity $9$ is completed by
the Nil manifolds $\big(K,4\big)$ and $\big(T,4\big)$ and the 
${\rm SL}_2$ manifolds 
$$\big(K,(2,1),0\big),\ \ \big(K,(2,1),1\big),
\ \ \big(T,(2,1),0\big),\ \ \big(T,(2,1),1\big).$$
There are then 6 Sol manifolds fibred over the circle and 17 fibred over the interval.
The 4 hyperbolic manifolds have been described in Section~\ref{statement:section}.

We have then non-geometric manifolds, 168 of them with JSJ decomposition
given by two blocks of the form $\big(D,(p,q),(r,s)\big)$, where $(p,q)$
is always $(2,1)$ for one of the blocks. Next we have
$$\big(D,(2,1),(2,1),(2,1)\big) \bigcup\nolimits_A \big(D,(2,1),(2,1)\big)$$
with $A=\tinymatr{1}{1}{-1}{0},\tinymatr{0}{1}{1}{0},\tinymatr{1}{1}{1}{0}$. Denoting
by $S$ the M\"obius strip and by $D_1$ the annulus we also have
$$\big(S,(2,1)\big) \bigcup\nolimits_A \big(D,(2,1),(2,1)\big)$$
with $A=\tinymatr{2}{1}{-1}{0},\tinymatr{1}{1}{1}{0},\tinymatr{1}{1}{-1}{0}$,
$$\big(D,(2,1),(2,1)\big) \bigcup\nolimits_A \big(D_1,(2,1)\big)
\bigcup\nolimits_B \big(D,(2,1),(2,1)\big)$$
with $(A,B)=
\left(\tinymatr{0}{1}{-1}{-1},\tinymatr{1}{1}{-1}{0}\right),
\left(\tinymatr{0}{1}{1}{0},\tinymatr{0}{1}{1}{0}\right),
\left(\tinymatr{0}{1}{1}{0},\tinymatr{1}{1}{-1}{0}\right)$, and
$$\big(D_1,(2,1)\big)_A$$
with $A=
\tinymatr{0}{1}{1}{0},
\tinymatr{0}{-1}{-1}{0},
\tinymatr{1}{1}{0}{-1},
\tinymatr{-1}{-1}{0}{1},
\tinymatr{-1}{-1}{1}{2},
\tinymatr{1}{1}{-1}{-2},
\tinymatr{1}{2}{0}{-1},
\tinymatr{-1}{-1}{0}{1}$.

\vspace{1cm}

{\small 
\noindent
\hspace{7.85cm}
\vbox{\noindent Dipartimento di Matematica Applicata\\ 
Via Bonanno Pisano, 25B\\ 
56126 Pisa, Italy\\
martelli@mail.dm.unipi.it\\
petronio@dm.unipi.it}}

\end{document}